\documentclass{amsart}
 \usepackage{amssymb,amsmath,color}
 \usepackage{graphicx}
 \usepackage{pdfsync,epsfig,subfigure}

\numberwithin{equation}{section}
\theoremstyle{plain}
\newcommand{\red}{\color{black}}

\newcommand{\clt}{central limit theorem}

\newcommand{\garch}{{\rm GARCH}$(1,1)$}

\newcommand{\sta}{St\u aric\u a}
\newcommand{\ex}{{\rm e}\,}

\newcommand{\asy}{asymptotic}

\newcommand{\Imath}{I}
\newcommand{\Jmath}{J}
\newcommand{\ts}{time series}
\newcommand{\tsa}{\ts\ analysis}

\definecolor{darkblue}{rgb}{.1, 0.1,.8}
\definecolor{darkgreen}{rgb}{0,0.8,0.2}
\definecolor{darkred}{rgb}{.8, .1,.1}

\newtheorem{thm}{Theorem}
\newtheorem{theorem}[thm]{Theorem}
\newtheorem{lemma}[thm]{Lemma}

\newcommand{\bbc}{{\mathbb C}}
\newtheorem{proposition}[thm]{Proposition}
\newtheorem{definition}[thm]{Definition}
\newtheorem{corollary}[thm]{Corollary}
\newtheorem{example}[thm]{Example}
\newtheorem{exercise}[thm]{Exercise}
\newtheorem{remark}[thm]{Remark}
\newtheorem{fig}[thm]{Figure}
\newtheorem{tab}[thm]{Table}

\newcommand{\bth}{\begin{theorem}}
\newcommand{\ethe}{\end{theorem}}
\newcommand{\sv}{stochastic volatility}
\newcommand{\bre}{\begin{remark}\em }
\newcommand{\ere}{\end{remark}}
\newcommand{\fre}{frequenc}
\newcommand{\ble}{\begin{lemma}}
\newcommand{\ele}{\end{lemma}}

\newcommand{\pp}{point process}
\newcommand{\bde}{\begin{definition}}
\newcommand{\ede}{\end{definition}}
\newcommand{\bco}{\begin{corollary}}
\newcommand{\eco}{\end{corollary}}

\newcommand{\bpr}{\begin{proposition}}
\newcommand{\epr}{\end{proposition}}

\newcommand{\bexer}{\begin{exercise}}
\newcommand{\eexer}{\end{exercise}}

\newcommand{\bexam}{\begin{example}}
\newcommand{\eexam}{\end{example}}

\newcommand{\bfi}{\begin{fig}}
\newcommand{\efi}{\end{fig}}

\newcommand{\btab}{\begin{tab}}
\newcommand{\etab}{\end{tab}}

\newcommand{\per}{periodogram}

\newcommand{\lhs}{left-hand side}
\newcommand{\fidi}{finite-dimensional distribution}
\newcommand{\rv}{random variable}

\newcommand{\var}{{\rm var}}

\newcommand{\cov}{{\rm cov}}

\newcommand{\rhs}{right-hand side}
\newcommand{\df}{distribution function}

\newcommand{\beao}{\begin{eqnarray*}}
\newcommand{\eeao}{\end{eqnarray*}\noindent}

\newcommand{\beam}{\begin{eqnarray}}
\newcommand{\eeam}{\end{eqnarray}\noindent}

\newcommand{\beqq}{\begin{equation}}
\newcommand{\eeqq}{\end{equation}\noindent}

\newcommand{\bce}{\begin{center}}
\newcommand{\ece}{\end{center}}

\newcommand{\barr}{\begin{array}}
\newcommand{\earr}{\end{array}}

\newcommand{\stp}{\stackrel{P}{\rightarrow}}
\newcommand{\std}{\stackrel{d}{\rightarrow}}

\newcommand{\stv}{\stackrel{v}{\rightarrow}}

\newcommand{\eqd}{\stackrel{d}{=}}

\newcommand{\w}{\omega}

\newcommand{\nto}{n\to\infty}
\newcommand{\kto}{k\to\infty}
\newcommand{\xto}{x\to\infty}

\newcommand{\ov}{\overline}
\newcommand{\wt}{\widetilde}
\newcommand{\wh}{\widehat}
\newcommand{\vep}{\varepsilon}

\newcommand{\la}{\lambda}

\newcommand{\regvary}{regularly varying}

\newcommand{\regvar}{regular variation}

\newcommand{\bbr}{{\mathbb R}}

\newcommand{\bbb}{{\mathcal B}}
\newcommand{\bbz}{{\mathbb Z}}

\newcommand{\con}{convergence}

\newcommand{\st}{such that}

\newcommand{\wrt}{with respect to}

\newcommand{\fct}{function}

\newcommand{\ds}{distribution}

\newcommand{\rep}{representation}
\newcommand{\cmt}{continuous mapping theorem}
\newcommand{\seq}{sequence}

\newcommand{\pro}{probabilit}

\newcommand{\ms}{measure}

\begin{document}
\title{The integrated periodogram  of a dependent extremal event
  sequence  }
\thanks{This paper was written when Yuwei Zhao was a PhD student at the
Department of Mathematics of the University of Copenhagen. He would
like to thank his Department for generous financial support. Major
parts of this paper were written in 2013 when
both authors spent sabbatical periods
at the Department of Statistics at Columbia University and the Forschungsinstitut
f\"ur Mathematik of ETH Z\"urich. Both authors take pleasure to thank 
their hosts for financial support and a constructive scientific
atmosphere. The research of Thomas Mikosch is partly supported by the
Danish Research Council Grant DFF-4002-000435.}
\author[T. Mikosch]{Thomas Mikosch}
\author[Y. Zhao]{Yuwei Zhao}
\address{Thomas Mikosch, University of Copenhagen, Department of Mathematics,
Universitetsparken 5,
DK-2100 Copenhagen\\ Denmark} \email{mikosch@math.ku.dk}
\address{Yuwei Zhao, Ulm University, Institute of Mathematical
  Finance, Helmholtzstrasse 18, D-89081 Ulm\\ Germany }\email{yuwei.zhao@uni-ulm.de}

% \begin{keyword}
%   Extreme value theory \sep functional central limit theorem \sep stationary
%   bootstrap \sep goodness-of-fit test \sep  spectral analysis
% \end{keyword}

\begin{abstract}
We investigate the \asy\ properties of the integrated   \per\ 
calculated from a \seq\ of indicator \fct s of dependent extremal events.
An event in Euclidean space is extreme if it occurs far away from the
origin.
We use a \regvar\ condition on the underlying stationary \seq\ to make these
notions precise. Our main result is a \fct al \clt\ for the integrated 
\per\ of the indicator \fct s of dependent extremal events. The limiting process is a continuous Gaussian process whose covariance structure 
is in general unfamiliar, but in the iid case a 
Brownian bridge appears. {\red In the general case, we propose a
stationary bootstrap procedure for approximating the \ds\ of the limiting
process. The developed theory can be used to construct 
classical goodness-of-fit tests such as the Grenander-Rosenblatt and Cram\'er-von Mises
tests which are based only on the extremes in the sample. We apply
the test statistics to simulated and real-life data.}
% \keywords{Extreme value theory \and functional limit theorem \and
%   stationary bootstrap \and goodness-of-fit test \and spectral
%   analysis }
% \subclass{60G70 \and 60F17 \and 62F40 \and 62F03 \and 62M15 \and 62M10 }
\end{abstract}
\maketitle
\section{Introduction} \label{sec:intr}%\setcounter{equation}{0}
\subsection{Regularly varying sequences}
We consider a strictly stationary $\bbr^d$-valued \seq\ $(X_t)$ for
some $d\ge 1$ with a generic element $X$ and assume that its \fidi s
are \regvary .
This means that for every $h\ge 1$, there exists a non-null Radon
\ms\ $\mu_h$  on the Borel $\sigma$-field $\ov \bbb_0^{dh}$ of
$\ov\bbr_0^{dh}=\ov \bbr^{dh}\backslash \{0\}$, $\ov \bbr=\{-\infty,\infty\}$, 
\st
\beam\label{eq:1a}
\dfrac{P(x^{-1} (X_1,\ldots,X_h)\in\cdot ) }{P(|X|>x)}\stv \mu_h(\cdot)\,,
\eeam
where $\stv$ denotes vague \con\ in $\ov\bbb_0^{dh}$; cf. Resnick
\cite{resnick:1987,resnick:2007}, Kallenberg \cite{kallenberg:1983}. 
The limiting \ms\ $\mu_h$ necessarily has the property $\mu_h(t
\cdot)=t^{-\alpha} \mu_h(\cdot)$, $t>0$, for some $\alpha\ge 0$, the
index
of \regvar . In what follows, we assume that $\alpha>0$. Relation
\eqref{eq:1a} is equivalent to the sequential definition
\beam\label{eq:2}
n\,P(a_n^{-1} (X_1,\ldots,X_h)\in\cdot ) \stv \mu_h(\cdot)\,,\quad \nto\,,
\eeam 
where $(a_n)$ is chosen \st\ $P(|X|>a_n)\sim n^{-1}$ as $\nto$.
We will say that the \seq\ $(X_t)$ and any of the vectors
$(X_1,\ldots,X_h)$, $h\ge 1$, are {\em \regvary\ with index $\alpha$}.
\par
Examples of \regvary\ strictly stationary \seq s are linear and \sv\ processes 
with iid \regvary\ noise, GARCH processes, infinite variance stable
processes and max-stable processes with Fr\'echet marginals. These
examples are discussed e.g. in 
Davis et al.
\cite{davis:mikosch:2009,davis:mikosch:cribben:2012,davis:mikosch:zhao:2013},
Mikosch and Zhao \cite{mikosch:zhao:2012}.
\subsection{The extremogram}\label{subsec:extremogram}
Consider a $\mu_1$-continuity Borel set 
$D_0=A\subset \ov \bbr_0^d$ bounded away from
zero and \st\ $\mu_1(A)>0$.
Then the sets 
$D_h=A\times \overline{\mathbb{R}}^{d(h-1)}\times A$ are
bounded away from zero as well and are continuity
sets \wrt\ the corresponding limiting \ms s $\mu_{h+1}$, $h\ge 1$.  We conclude
from \eqref{eq:2} that the limits  
\beam\label{eq:rv}
\gamma_A(h)=\lim_{\nto} n\,P(a_n^{-1}X_0\in A\,, 
a_n^{-1}X_h\in A)= \mu_{h+1}(D_h)\,, \quad h\ge 0\,,
\eeam
exist. For $t\in\bbz$, it is not difficult to see that
\beao
n\,\cov (I_{\{a_n^{-1} X_t\in A\}},I_{\{a_n^{-1}X_{t+h}\in A\}})&\sim& 
n E I_{\{a_n^{-1} X_t\in A\,,a_n^{-1}X_{t+h}\in A\}}\\
&=& n\,P(a_n^{-1}X_0\in A\,, 
a_n^{-1}X_h\in A)\\& \to& \gamma_A(h)\,,\quad \nto
\,.
\eeao
Hence $\gamma_A$ constitutes the covariance  \fct\ of a stationary
process.
We refer to $\gamma_A$ as the {\em extremogram relative to the set $A$.}
We will also consider the  {\em standardized extremogram} 
given as the limiting \seq 
\beao
\rho_A(h)= \lim_{\nto}P(a_n^{-1} X_h\in A\mid a_n^{-1} X_0\in
A)=\dfrac{\mu_{h+1}(D_h)}{\mu_1(D_0)}\,,\quad h\ge 0\,.
\eeao The quantities $\rho_A(h)$ have an intuitive
interpretation as limiting conditional \pro ies. Moreover, $\rho_A$
is the autocorrelation \fct\ of a stationary process. 
The quantities $\rho_A(h)$
are generalizations of  the upper tail dependence coefficient of a 
two-dimensional vector $(Y_1,Y_2)$ with identical marginals given as the limit $\lim_{\xto}
P(Y_2>x\mid Y_1>x)$. 
\par
The extremogram was introduced in Davis and Mikosch
\cite{davis:mikosch:2009} as a \ms\ of serial extremal dependence in a
strictly stationary \seq . There and in Davis et al. 
\cite{davis:mikosch:cribben:2012,davis:mikosch:zhao:2013} 
various aspects of the estimation of the extremogram were discussed,
including \asy\ theory and the use of the stationary bootstrap for
the construction of confidence bands. 
\subsection{The sample extremogram}
Natural estimators of the extremograms $\gamma_A$ and $\rho_A$ 
are given by
their respective sample analogs 
\beao
\widetilde{\gamma}_A(h) & =& \frac{m}{n} \sum_{t=1}^{n-h}
\widetilde{I}_t\widetilde{I}_{t+h}\quad \mbox{and}\quad
\widetilde{\rho}_A(h)
=\dfrac{\widetilde{\gamma}_A(h)}{\widetilde{\gamma}_A(0)}\,,\quad h\ge
0\,.
\eeao
Here $m=m_n$ is any integer \seq\ satisfying the conditions
$m_n\to\infty$ and $m_n/n=o(1)$ and 
\beao
I_t=I_{\{a_m^{-1}X_t\in A\}}\,,\quad  \wt I_t= I_t-
p_0\,,\quad\mbox{and}\quad  p_0= E I_t=
P(a_m^{-1} X\in A)\,, t\in\bbz\,.
\eeao
It is shown in Davis and Mikosch \cite{davis:mikosch:2009} that the
conditions $m_n\to\infty$ and $m_n/n=o(1)$ are needed 
for the validity of the \asy\ properties
$E\wt \gamma_A(h)\to \gamma_A(h)$ and $\var(\wt \gamma_A(h))\to 0$ as
$\nto$. Moreover, under a mixing condition, the \fidi s of $\wt
\gamma_A$ and $\wt \rho_A$ satisfy a \clt\ with normalization
$(n/m)^{1/2}$; cf. Lemma~\ref{thm:acf} below.
\subsection{Spectral density and periodogram}\label{subsec:spectral}
Since $\gamma_A$ and $\rho_A$ are the autocovariance and autocorrelation
\fct s of a stationary
process, respectively, it is possible to enter the corresponding \fre
y domain. If $\gamma_A$ is square summable one can define the
spectral densities 
\beao
h_A(\la)=\sum_{h\in\bbz} \gamma_A(h)\, \ex^{-ih\la }\quad
\mbox{and}\quad f_A(\la)=\sum_{h\in \bbz} \rho_A(h)\, \ex^{-ih\la }\,,\;
\la\in [0,\pi]=\Pi\,.
\eeao
A natural estimator of the spectral density is the \per . Since the
sample autocovariances $\wt \gamma_A(h)$ are derived from the
triangular array of the stationary \seq s $(\wt I_t)$, an analog of
the classical \per\ for $h_A$ is given by  
\beao%\label{eq:peri1}
I_{n,A}(\lambda) & = & \frac{m}{n} \Big|\sum_{t=1}^n \widetilde{I}_t\,
\ex^{-i\,t\,\lambda} \Big|^2=\widetilde{\gamma}_A (0)+2 \sum_{h=1}^{n-1}\widetilde{\gamma}_A(h)\cos(h\lambda) \,, \quad \lambda \in \Pi\,,
\eeao
and the \per\ for the standardized spectral density $f_A$ is obtained as the
scaled  \per\ $I_{n,A}/\wt \gamma_A(0)$. Mikosch and Zhao
\cite{mikosch:zhao:2012} showed under mixing conditions that the {\em
  extremal \per\ ordinates}  $I_{n,A}(\la)$
share various of the classical properties of the \per\ ordinates 
for a 
stationary \seq\ (cf. Brockwell and Davis \cite{brockwell:davis:1991}): consistency in the mean, \con\ in \ds\ to
independent exponential \rv s with expectation $h_A(\la_j)$ at distinct fixed \fre ies $\la_j\in
(0,\pi)$ and at
distinct  Fourier \fre ies $\w_n(j)=2\pi j/n\in (0,\pi)$ provided
these \fre ies converge to a limit $\la_j\in (0,\pi)$ as $\nto$. The
latter property ensures that weighted versions of the
\per\ $I_{n,A}$ at fixed \fre ies $\la\in (0,\pi)$ converge in mean
square to $h_A(\la)$.
\par
For practical purposes, one will mostly work with the \per\ at
the Fourier \fre ies $\w_n(j)\in (0,\pi)$. Then 
\beao
I_{n,A}(\w_n(j))=
  \frac{m}{n} \Big|\sum_{t=1}^n I_t\,
\ex^{-i\,t\,\w_n(j)} \Big|^2\,,
\eeao
i.e., centering of the indicator \fct s $I_t$ is not needed.
However, for proving \asy\ theory it will be convenient to work 
with the extremal \per\ $I_{n,A}$ based on the centered quantities $\wt I_t$,
$t=1,\ldots,n$. 

\subsection{The integrated \per }
The integrated \per\ of a stationary \seq\ has a long history in \tsa
, starting with classical work of Grenander and Rosenblatt
\cite{grenander:rosenblatt:1984},  and was extensively used in the
monographs Hannan \cite{hannan:1960}, Priestley \cite{priestley:1981},
Brockwell and Davis \cite{brockwell:davis:1991}, to name a few references.
Dahlhaus \cite{dahlhaus:1988} discovered a close relationship of the
integrated \per , considered as a process indexed by \fct s, and
empirical process theory. Under entropy conditions, he proved 
uniform \con\ results over suitable classes of index \fct s; see also
the survey paper Dahlhaus and Polonik \cite{dahlhaus:polonik:2002}.
These papers gave some general theoretical background for 
various \per\ based techniques such as Whittle estimation of the
parameters of a FARIMA process and goodness of fit tests for linear
processes as mentioned in Grenander and Rosenblatt
\cite{grenander:rosenblatt:1984} and Priestley \cite{priestley:1981}.
\par 
In this paper, we will consider the integrated \per
\beam\label{eq:intper}
J_{n,A}(g)&=&\int_\Pi I_{n,A}(\lambda)\,g(\lambda) \,d\la = 
c_0(g)\,\widetilde{\gamma}_A(0)+2
\sum_{h=1}^{n-1}c_h(g)\,\widetilde{\gamma}_A(h)\,,
\eeam
and its standardized version
\beao
% \label{intper1}
J_{n,A}^{\circ}(g)&=&\dfrac{1}{\wt \gamma_A(0)}\int_\Pi I_{n,A}(\lambda)\,g(\lambda) \,d\la =  c_0(g) +2
\sum_{h=1}^{n-1}c_h(g)\,\widetilde{\rho}_A(h)\,,\nonumber
\eeao
where $g$ is non-negative and square integrable \wrt\ Lebesgue \ms\ on
$\Pi$ (we write $g\in L_+^2(\Pi)$)
with corresponding Fourier coefficients
\beao
c_h(g)= \int_\Pi \cos(h\lambda)\,g(\lambda) \,d\la\,,\quad h\in\bbz\,.
\eeao
We will understand $J_{n,A}(g)$ and $J_{n,A}^{\circ}(g)$ as  
natural estimators of 
\beam\label{eq:intpera}
J_A(g)=\int_\Pi h_A(\lambda)\, g(\la)\,d\la &= & c_0(g)\,\gamma_A(0)+2
\sum_{h=1}^{\infty} c_h(g)\,\gamma_A(h)\,,\\
%\label{eq:intpera1}
J_A^{\circ}(g)=\int_\Pi f_A(\lambda)\, g(\la)\,d\la &= & c_0(g) +2
\sum_{h=1}^{\infty} c_h(g)\,\rho_A(h)\,,\nonumber
\eeam
respectively.
The latter identities hold if $\sum_{h=0}^\infty\gamma_A(h)<\infty$, a
condition we assume throughout this paper; {\red see also
  Remark~\ref{rem:2} below.}
\par
The main results of this paper (see Section~\ref{sec:fclt}) are \fct al \clt s for the integrated \per\ $J_{n,A}$  with $g=h I_{[0,\cdot]}$ 
for a sufficiently smooth \fct\ $h$ on $\Pi$. 
The limit processes are 
Gaussian whose covariance structure strongly depends on 
the limit \ms s $(\mu_h)$. The rate of \con\ in these results is typically slower than $\sqrt{n}$.  However, in the case of an iid \seq , 
the limiting process  is a Brownian bridge and the \con\ rates are 
much faster than in the case 
of a dependent \seq . These results differ from classical theory for 
the \per\ of a stationary \seq\ $(X_t)$ (see e.g. Dahlhaus \cite{dahlhaus:1988}, Kl\"uppelberg and Mikosch \cite{kluppelberg:mikosch:1996}), 
where the limiting process is completely determined by 
the covariance structure of $(X_t)$. 
The methods of proof combine 
classical techniques of weak \con\ and strong mixing 
(e.g. Billingsley \cite{billingsley:1999}) with extreme value theory 
for dependent \seq s (e.g. Davis and Mikosch \cite{davis:mikosch:2009}).
The proofs are rather technical due to the fact that the \seq s of 
indicator \fct s $(I_t)$ have triangular structure: they change in
dependence on the threshold $a_m$.
\par
{\red As in classical \tsa , the \fct al central limit theory 
for the integrated \per\ can be used to construct \asy\ goodness-of-fit
tests such as the Grenander-Rosenblatt and Cram\'er-von Mises
tests. In contrast to their classical counterparts, these tests are based only
on the extremal part of the underlying sample, i.e., we test whether the
extremes of the sample are in agreement with the null hypothesis
about a given type of \ts\ model. Such tests may be useful, for example, for
distinguishing between a GARCH and a \sv\ model fitted 
to a return \ts . The aforementioned two types of models may have similar 
autocorrelation structure for the data, their absolute values and
squares, so their spectral properties are very similar as well, 
while their extremograms are rather distinct: the extremogram
$\gamma_A$ relative to the set $A=(1,\infty)$ decays exponentially fast for
GARCH and for the simple \sv\ model $\gamma_A$ vanishes at all
positive lags; see Davis and Mikosch \cite{davis:mikosch:2009}.}
\par
The paper is organized as follows.
We start in Section~\ref{sec:2} with 
some moment calculations and we also introduce the relevant mixing conditions
and central limit theory for the sample extremogram. In 
Section~\ref{sec:3} we provide a result about the 
mean square consistency of the integrated \per ; the proof is given in
Section~\ref{sec:lemma}.
The main results (Theorems~\ref{thm:main} and \ref{thm:mdep2})  are
\fct al central limit theorems for the integrated \per . 
They are given in Section~\ref{sec:fclt}; the corresponding
proofs are provided in Sections~\ref{sec:proof1} and
\ref{sec:proof2}. The covariance structure of the limiting Gaussian
processes in Theorem~\ref{thm:main} is rather complicated. Therefore
in Section~\ref{sec:bootstrap} we
supplement the \asy\ theory  
by consistency results for the
stationary bootstrap applied to the integrated \per\ of extremal events
in a strictly stationary \seq . 
The corresponding proofs are given in 
Section~\ref{sec:proof3}. 
In Section~\ref{sec:simulation} we indicate
how the integrated \per\ works for simulated and real-life data.

\section{Preliminaries}%\setcounter{equation}{0}
\label{sec:2}
\subsection{Some moment calculations}
Recall the notation and conditions of 
Section~\ref{sec:intr}. We write
\beao
p_0=P(a_m^{-1} X_0 \in A)\quad\mbox{and}\quad
p_h= P(a_m^{-1}X_0 \in A, a_m^{-1}X_h \in A)\,,\quad h\ge 1\,,
\eeao
where as above, $m_n\to\infty$ and $m_n/n=o(1)$ as $\nto$.
For integers $s,t,u,v\ge 0$, we set
\beao
\Gamma(s,t,u,v) &=& E \widetilde{I}_s\widetilde{I}_t \widetilde{I}_u
\widetilde{I}_v\,,\\ 
 \Gamma(s,t,u)& =& E\widetilde{I}_s\widetilde{I}_t \widetilde{I}_u\,,\\
\Gamma(s,t)&=& E \widetilde{I}_s \widetilde{I}_t=p_{|s-t|}-p_0^2\,.
\end{eqnarray*}
\par
We will often have to calculate variances and
covariances of the sample extremogram $\wt \gamma_A$. 
We provide some of these formulas for further use. 
\ble\label{lem:varcov} Let $(X_t)$ be  a strictly stationary \seq .
Then, for $1\le h \le n-1$, 
\beao
%\label{eq:var}
 (n/m)^2E \wt \gamma_A^2(h)
%\lefteqn{E \sum_{t=1}^{n-l}\widetilde{I}_t(A)\widetilde{I}_{t+l}(A)\sum_{s=1}^{n-l}\widetilde{I}_s(A)\widetilde{I}_{s+l}(A)= E\Big(\sum_{t=1}^{n-l}\widetilde{I}_t(A)\widetilde{I}_{t+l}(A) \Big)^2}\\
&= & (n-h) E(\wt I_0 \wt I_h)^2+
2\sum_{t=1}^{n-h-1}(n-h-t)\Gamma(0,h,t,t+h)
\eeao
and for 
$1\le h<h+u\le n-1$,
\beao%\label{eq:cov}
&&(n/m)^2 \,E \wt \gamma_A(h)\wt \gamma_A (h+u)\\
 & =& (n-h-u)\Gamma(0,h,0,h+u)\\&& + \sum_{t=1}^{n-h-u-1}(n-h-u-t)
\Gamma(0,h,t,t+h+u)\\ 
& &  + \sum_{t=1}^{n-h-1}\min(n-h-u,n-h-t)\Gamma(0,h+u,t,t+h)\,.
\eeao
\ele
\subsection{Mixing conditions}
The following two mixing conditions were introduced in   Davis and
Mikosch \cite{davis:mikosch:2009} for a strongly mixing
$\bbr^d$-valued  \seq\ $(X_t)$ with rate \fct\ $(\xi_h)$. 
\subsubsection*{ Condition {\rm (M)} }
There exist integer \seq s $m=m_n\to\infty$ and  
  $r_n\to\infty$  
\st\ $m_n/n\to 0$, $r_n/m_n\to 0$ and  
\beam\label{eq:m1}  
\lim_{\nto}m_n \sum_{h=r_n}^\infty \xi_h&=&0\,,
\eeam
Moreover, an anti-clustering condition holds:
\beam  
\label{eq:m2}  
\lim_{k\to\infty}\limsup_{\nto} \sum_{h=k}^{r_n} P(|X_h|>\epsilon  
\,a_m \mid|X_0|>\epsilon\, a_m)&=&0\,,\quad \epsilon >0\,.  
\eeam
\subsubsection*{Condition {\rm (M1)}}
Assume (M) and that
the \seq s  $(m_n)$, $(r_n)$, $k_n=[n/m_n]$ from (M) also satisfy  
the growth conditions $k_n \xi_{r_n} \to 0$, and $m_n=o(n^{1/3})$.   
\bre
The condition $m_n=o(n^{1/3})$ in (M1) can be replaced 
by  $\frac{m_nr_n^3}{n} \to 0$ and
$
\frac{m_n^4}{n} \sum_{j=r_n}^{m_n}\xi_j \to 0$
which is often much weaker. 
\ere
Condition \eqref{eq:m1} is easily satisfied if 
the mixing rate $(\xi_h)$ is geometric, i.e., exponentially 
decaying to zero. Under mild conditions, the popular classes of 
ARMA, max-stable, 
GARCH and \sv\ processes are strongly mixing with geometric rate;
cf. Davis et al. 
\cite{davis:mikosch:2009,davis:mikosch:cribben:2012,davis:mikosch:zhao:2013,mikosch:zhao:2012} for discussions 
of these examples. 
Condition \eqref{eq:m2} is similar to (2.8) in Davis
and Hsing \cite{davis:hsing:1995}. It serves the purpose of
establishing the convergence of a sequence of point processes to a
limiting cluster \pp . This condition is much weaker than the anti-clustering
condition $D'(\epsilon a_n)$ of Leadbetter; cf. Section 5.3.2 in
Embrechts et al. \cite{embrechts:kluppelberg:mikosch:1997}.
\par
The mixing rate $(\xi_h)$ in conditions (M) and (M1) is useful for
finding bounds on the moments $\Gamma(s,t,u,v)$ introduced above. 
{\em In what follows, $c$ will denote any (possibly different) 
constants whose value is not of interest.}
\ble
\label{lem:nb}
Let $(X_t)$ be a strongly mixing sequence with mixing rate
$(\xi_h)$. Then for integers $h,l,u\ge 1$ and for some constants $c>0$
{\red which do not depend on $n$},
\begin{eqnarray}
&& |\Gamma(0,h,h+l,h+l+u)| \le c\,\min(\xi_h,\xi_u)\,,\label{eq:nb1}\\
&&|\Gamma(0,h,h+l,h+l+u)
-(p_h-p_0^2)(p_u-p_0^2)|\le c\,\xi_l\,,\label{eq:nb2}\\
\label{eq:nb3}
&& |\Gamma(0,h,h+l)| \le c\,\min(\xi_h,\xi_l)\,,\\ 
\label{eq:nb4}
&&|\Gamma(0,h)| \le  \xi_h\,.
\end{eqnarray}
\ele
The proof of Lemma \ref{lem:nb} follows by a direct application of 
Theorem 17.2.1 in Ibragimov and Linnik \cite{ibragimov:linnik:1971}.
Relation \eqref{eq:nb1} combined with \eqref{eq:m1} will ensure 
that sums of $\Gamma(0,h,h+l,h+l+u)$ are \asy ally negligible if
$h$ or $u$ exceed~$r_n$.

\subsection{Central limit theory for the sample extremogram} \label{sec:clt}
In this section we recall a \clt\ for the extremogram from Davis and
Mikosch \cite{davis:mikosch:2009}, Section 3. 
\ble \label{thm:acf}
Assume that $(X_t)$ is an  $\mathbb{R}^d$-valued 
strictly stationary \regvary\ sequence with index $\alpha >0$
and that the Borel set $A$ satisfies the
conditions of Section~\ref{subsec:extremogram}.
If the mixing conditions 
{\rm (M)}, {\rm (M1)} hold and $\sum_{l=1}^\infty
\gamma_A(l)<\infty$ then for $h\ge 0$,
\beam
\label{eq:acf00}\widetilde{\gamma}_A(h)&\overset{P}{\to}& \gamma_A(h)\,,\\
 \label{eq:acf11}(n/m)^{1/2} \big(
\widetilde{\gamma}_A(i)-E \widetilde{\gamma}_A(i)\big)_{i=0,\ldots,h} 
&\std& (Z_i)_{i=0,\ldots,h} \,,
\eeam
where $(Z_i)_{i=0,\ldots,h}$ is Gaussian with mean zero and
covariance matrix $\Sigma_h=(\sigma_{ij})_{i,j=0,\ldots,h}$
given by
\beao
\sigma_{ij} & =& \gamma_A(i,j)+  \sum_{l=1}^{\infty}\big[
\gamma_A(i,l,l+j)+\gamma_A(j,l,l+i)\big]\,, \quad i,j=0,\ldots,h\,,
\eeao
and for $u,s,t\ge 0$,
\beao
\gamma_A(u,s,t)&=& \lim_{n\to\infty }n\,P(a_n^{-1} X_0\in A\,,a_n^{-1}
X_u\in A, a_n^{-1} X_s\in A, a_n^{-1} X_{t}\in A)\,,
\eeao
with the convention that 
$\gamma_A(u,t)=\gamma_A(u,u,t)$. Moreover, we have for $h\ge 1$ 
\beam
\wt \rho_A(h) &\stp& \rho_A(h)\,,\label{eq:opu1}\\
(n/m)^{1/2} \Big(\wt \rho_A(i)-
\frac{p_i}{p_0}\Big)_{i=1,\ldots,h}&\std& \frac 1 {\gamma_A(0)} \Big(
Z_i- \rho_A(i) Z_0\Big)_{i=1,\ldots,h}\,.\label{eq:opu2} 
\eeam
\ele
\begin{proof}
The proof of \eqref{eq:acf00}  was given 
in Section 3 of Davis and Mikosch 
\cite{davis:mikosch:2009}. There we can also find 
the proof of \eqref{eq:acf11} in a more general context.
Here we will calculate the
covariance matrix $\Sigma_h$ explicitly. The expressions 
for $\sigma_{ii}$, $i\ge 0$, were derived in 
Davis and Mikosch \cite{davis:mikosch:2009} for $i=0$ and
$i\ge 1$ in Theorem 3.1 and Lemma~5.2, respectively. 
We notice that $\gamma_A(i,l,l+j)\le
\gamma_A(l)$ and therefore the infinite series in $\sigma_{ij}$ are
finite. 
\par
For $i\ne j$, similar calculations as for 
Lemma \ref{lem:varcov} yield for $k\ge 1$ and \mbox{$r_n/m_n\to 0$,}
\beao
\lefteqn{\frac{m}{n}\cov\Big(\sum_{t=1}^n \widetilde{I}_t\widetilde{I}_{t+i},\sum_{s=1}^n \widetilde{I}_s\widetilde{I}_{s+j}\Big)}\\
 &=&
m\, \Gamma(0,0,i,j) + m\,\sum_{l=1}^{n}\Big[ (1-l/n)
\big[\Gamma(0,i,l,l+j)\\
&& \quad+\Gamma(0,j,l,l+i)\big]-
(p_i-p_0^2)(p_j-p_0^2)\Big]\\
&=&m\, \Gamma(0,0,i,j) +m\Big(\sum_{l=1}^k +\sum_{l=k+1}^{r_n}+\sum_{l=r_n+1}^{n}\Big)\Big[ (1-l/n)
\big[\Gamma(0,i,l,l+j)\\
&& \quad+\Gamma(0,j,l,l+i)\big]-
(p_i-p_0^2)(p_j-p_0^2)\Big]\\
&=& Q_1+Q_2+Q_3+Q_4\,.
\eeao
By \regvar , for fixed $k\ge 1$ as $\nto$, 
\beao
Q_1+Q_2 \to\gamma_A(i,j)+ \sum_{l=1}^{k}\big[
\gamma_A(i,l,l+j)+\gamma_A(j,l,l+i)\big]\,,
\eeao
and the \rhs\ converges to $\sigma_{ij}$ as
$\kto$.
By \eqref{eq:m2}, we have
\beao
\lim_{\kto}\limsup_{\nto} |Q_3|=0\,.
\eeao
Using \eqref{eq:nb2} and \eqref{eq:m1}, we also have
$
|Q_4|\le c m_n\sum_{l=r_n+1}^\infty \xi_l\to 0$ as $\nto$.
This proves \eqref{eq:acf00} and \eqref{eq:acf11}. Relations
\eqref{eq:opu1} and \eqref{eq:opu2} follow by a continuous mapping
argument, observing that for $1\le i\le h$,
\beao
\Big(\frac nm\Big)^{1/2} 
\Big(\wt\rho_A(i) - p_i/p_0\Big)&=&\Big(\frac nm\Big)^{1/2}
    \dfrac{\wt\gamma_A(i)-E\wt\gamma_A(i)}{\wt \gamma_A(0)} \\
&& \quad- E\wt \gamma_A(i)
\dfrac{(n/m)^{1/2}
\big(\wt \gamma_A(0)-E\wt\gamma_A(0))}{\wt\gamma_A(0) E \wt
  \gamma_A(0)} +
o_P(1)\\
&\std & \frac{1}{\gamma_A(0)} \Big(Z_i- \rho_A(i) Z_0\Big)\,. 
\eeao
\qed
\end{proof}
\bre\label{rem:2}{\red 
The summability condition on $\gamma_A$ which we assume in the
previous lemma and 
throughout this paper is satisfied for a large variety of \regvary\
\ts\ models; see the calculation of $\gamma_A$ in
Davis et
al. \cite{davis:mikosch:2009,davis:mikosch:cribben:2012,davis:mikosch:zhao:2013}.
For example, finite order ARMA models with \regvary\ iid noise and GARCH models
have exponentially decaying extremogram, and the simple \sv\ model with
log-normal volatility process has vanishing extremogram at all 
positive lags. Formulas for $\gamma_A$ also exist for infinite
variance stable and max-stable processes with Fr\'echet
marginals. Also for these processes the summability condition on
$\gamma_A$ may hold, depending on the specification of the process. }
\ere
Recall that a strictly stationary process $(X_t)$ is $\eta$-dependent
for some integer $\eta\ge 0$ if $(X_t)_{t\le 0}$ and $(X_t)_{t>\eta}$
are independent. For such a process we observe that 
$\sigma_{hh}=0$ for $h>\eta$ and 
hence \eqref{eq:acf11}
collapses into $(n/m)^{0.5} \widetilde{\gamma}_A(h)\stp 0$ for $h>
\eta$. In particular, for an iid \seq\ $(X_t)$, $Z_h=0$ a.s. for 
$h\ge 1$, while $(n/m)^{0.5} \widetilde{\gamma}_A(0)\std Z_0$ and
$Z_0$ is $N(0,\gamma_A(0))$ distributed.
\par
In these cases, the rate of convergence in
\eqref{eq:acf11} can be improved.  
\ble \label{thm:mdep1}
Assume that $(X_t)$ is an $\mathbb{R}^d$-valued $\eta$-dependent
\regvary\ strictly stationary sequence with index $\alpha >0$ for some
$\eta\ge 0$, and 
the Borel set $A$ satisfies the conditions of
Section~\ref{subsec:extremogram}. Additionally, assume that for
$j\ge i>\eta$ and $1\le t\le \eta- (j-i)$, the following limits exist:
\begin{eqnarray} \label{eq:mdep1}
\lefteqn{\overline{\gamma}_A(t,i,t+j)}\nonumber\\&=&\lim_{\red m\to\infty} m^2P(a_m^{-1}X_0\in A, a_m^{-1}X_t\in A,a_m^{-1}X_i\in A,a_m^{-1}X_{t+j}\in A)\,.\nonumber\\
\end{eqnarray} 
Then for $h\ge 1$,
$
n^{0.5}
\big(\widetilde{\gamma}_A(\eta +i)\big)_{i=1,\ldots,h} 
\std (Z_i)_{i=1,\ldots,h}\,,
$
where $(Z_i)_{i=1,\ldots,h}$ is Gaussian
$N(0,\overline{\Sigma}_h)$ whose
covariance matrix 
$\overline{\Sigma}_h=(\sigma_{ij})_{i,j=1,\ldots,h}$ is given by
\beam\label{opa}
 \sigma_{ij}&= &
\gamma_A(0)\gamma_A(j-i)+ \sum_{t=1}^{\eta -(j-i)}\big[
\ov \gamma_{\red A}(t,i,t+j)
+\ov \gamma_{\red A}(t,j,t+i)\big]\,,\\&&1\le i\le j\,.\nonumber  
\eeam
\ele
\bre
Condition \eqref{eq:mdep1} is an additional \asy\ independence condition.
Indeed, \regvar\ of $(X_t)$ only implies that the limits 
\beao
\lim_{\red m\to \infty} m P(a_m^{-1}X_0\in A, a_m^{-1}X_t\in A,a_m^{-1}X_i\in A,a_m^{-1}X_{t+j}\in A)
\eeao 
exist and are finite. Then \eqref{eq:mdep1} implies that the latter
limits must be zero. In Example~\ref{exam:1} we consider some simple cases
when  \eqref{eq:mdep1} is satisfied.
\ere
\bre  Assume $j-i>\eta$. Then, by $\eta$-dependence, 
$\gamma_A(j-i)=0$ and the index set in \eqref{opa} is empty. 
Hence $\sigma_{ij}=0$ for $j-i>\eta$. In particular, if $(X_t)$ is
iid, $\sigma_{ij}=0$ for $i\ne j$ and $\sigma_{ii}=\gamma_A^2(0)$.
\ere
\begin{proof}
We start by calculating the
asymptotic covariances. Assume $j\ge i>{\red \eta}$. Then, using the
independence of $I_0$ and $(I_jI_i,I_iI_tI_{t+j},I_jI_tI_{t+i})$ for
  $t>\eta$ and of $I_{t+j}$ and $I_0I_t I_i$ for  $t\le \eta$ and
  $t\ge \eta - (j-i)$, we obtain
\begin{eqnarray*}
\lefteqn{\cov\big(n^{0.5}\widetilde{\gamma}_A(i),n^{0.5}\widetilde{\gamma}_A(j)\big)}\\&=& m^2
E\widetilde{\Imath}_0^2E\widetilde{\Imath}_i\widetilde{\Imath}_j+  
m^2  \sum_{t=1}^{\eta}\big[
 E\widetilde{I}_0\widetilde{I}_i\widetilde{I}_{t}\widetilde{I}_{t+j}
+E\widetilde{I}_0\widetilde{I}_j\widetilde{I}_{t}\widetilde{I}_{t+i}\big]
 +o(1)\\
&=& \gamma_A(0)\gamma_A(j-i)+m^2  \sum_{t=1}^{\eta-(j-i)}\big[
 E\widetilde{I}_0\widetilde{I}_i\widetilde{I}_{t}\widetilde{I}_{t+j}
+E\widetilde{I}_0\widetilde{I}_j\widetilde{I}_{t}\widetilde{I}_{t+i}\big]
+o(1)\\
&\to & \gamma_A(0)\gamma_A(j-i) + \sum_{t=1}^{\eta-(j-i)}\big[ \ov
\gamma_A(t,i,t+j)+\ov\gamma_A(t,j,t+i)\big] \,,\quad \nto\,.
\end{eqnarray*}
In the last step we used \eqref{eq:mdep1}.
This completes the calculation of $\overline{\Sigma}_h$. Furthermore,
we observe that for $h\ge 1$,
\beam\label{eq:mdep22}
n^{0.5}\big(
\widetilde{\gamma}_A(i)\big)_{i=\eta +1,\ldots,\eta+h}
&=& (m/n^{0.5}) \sum_{t=1}^n 
\big( \wt I_t \wt I_{t+i}\big)_{i=\eta +1,\ldots,\eta+h}+
o_P(1)\,.\nonumber\\
\eeam
The vector \seq\ $( \wt I_t \wt I_{t+i})_{i=\eta+1,\ldots,\eta+h}$,
$t=1,2,\ldots$, is strictly stationary and $(h+\eta)$-dependent.
Now an application of the \clt\ for strongly mixing triangular arrays in Rio \cite{rio:1995} and the Cram\'er-Wold device 
to \eqref{eq:mdep22} conclude the proof. 
\qed
\end{proof}
The following examples fulfill the conditions of Lemma~\ref{thm:mdep1}.
\bexam \label{exam:1}\rm
An iid regularly varying sequence $(X_t)$ is $0$-dependent, and thus
 \eqref{eq:mdep1} holds. Its
 limiting covariance matrix $\overline{\Sigma}_h$ is a diagonal matrix
 with entries $\gamma_A^2(0)=(\mu_1(A))^2$ on the main diagonal. 

We consider the stochastic volatility model $X_t=\sigma_tV_t$ where
 $(\sigma_t)$ is independent of $(V_t)$, $(\sigma_t)$ is a positive
 $\eta$-dependent strictly stationary sequence and $(V_t)$ is a regularly
 varying iid sequence with index $\alpha>0$; see Davis and Mikosch \cite{davis:mikosch:2009b}. Assume that
 $E\sigma^{\alpha+\varepsilon}<\infty$ for some $\varepsilon>0$. In
 this case, $(X_t)$ is $\eta$-dependent, strictly stationary and
 regularly varying with index $\alpha$. We
 will show that \eqref{eq:mdep1} holds with
 $\overline{\gamma}_A(u,s,t)=0$ for $0<u<s<t$. 
Since $A$ is bounded away from zero,
 there exists a $\delta>0$ such that
\begin{eqnarray*}
 &&\overline{\gamma}_A(u,s,t)\\&\le &\limsup_{\red m\to\infty}m^2P(a_m^{-1}\min(|X_0|,|X_u|,|X_s|,|X_t|)>\delta )\\ 
 & \le&
\limsup_{\red m\to\infty}m^2P(a_m^{-1}\max(\sigma_0,\sigma_u,\sigma_s,\sigma_t)\min(|V_0|,|V_u|,|V_s|,|V_t|)>\delta)\\
 & \le&\limsup_{\red m\to\infty} 4m^2P(a_m^{-1}\sigma_0\min(|V_0|,|V_u|,|V_s|,|V_t|)>\delta)\\
 & \le& \limsup_{\red m\to\infty}cm^2(E\sigma^{\alpha})^4(P(|V_0|>a_m\delta))^4=0 \,,
\end{eqnarray*}
where we used that $P(\sigma_0|V_0|>a_m)\sim
E\sigma^{\alpha}P(|V_0|>a_m\delta)$
by virtue of Breiman's lemma; see \cite{breiman:1965}.
\eexam
{\red In the iid case, the limiting quantities 
$Z_h$, $h\ge 1$, in Lemma~\ref{thm:acf} vanish. The same observation
can be made in the case of a strictly stationary \seq\ with \asy\ (extremal)
independence in the following sense:\\[1mm]
{\em Condition {\rm (AI)}}: 
Assume there exist \seq s  $m=m_n\to\infty$ and $r_n\to\infty$ \st\
$m=o(n)$ and $r_n=o(m)$ as $\nto$ and the following conditions are
satisfied for any Borel set $A\subset \bbr^d$ bounded away from zero
and the axes \st\ $\mu_1(\partial A)=0$:
\begin{enumerate}
\item
$\lim _{n\to\infty}m^2\,p_h$ exists and is finite for $h\ge 1$,
\item
$\lim_{n\to\infty} m^2
\sup_{1\le
  i<j\le r_n}P(a_m^{-1}X_0\in A\,,a_m^{-1} X_i\in
 A\,,a_m^{-1} X_j\in A)= 0\,,$
\item{\small
$\lim_{\nto} r_n m^2\,\sup_{1\le
  i<j<t\le r_n}P(a_m^{-1}X_0\in A\,,a_m^{-1} X_i\in
 A\,,a_m^{-1} X_j\in A,a_m^{-1} X_t\in A)= 0\,.$}
\end{enumerate}
\bexam\rm  
We consider the \sv\ model from Example~\ref{exam:1} but we
drop the condition of $\eta$-dependence. Conditions 
(AI.2) and (AI.3) are verified in the same way as in Example~\ref{exam:1}.
We also observe that for some constant $c>0$,
\beao
m^2\,p_h&\sim& c\,\dfrac{P(a_m^{-1}X_0\in A, a_m^{-1} X_h\in
  A)}{P(\min (V_1,V_2)>a_m)}\\
&=& c\,\dfrac{P(a_m^{-1}{\rm diag}(\sigma_0,\sigma_h) (V_1,V_2)'\in
  A\times A)}{P(\min (V_1,V_2)>a_m)}\,.
\eeao  
An application of a Breiman-type result 
for \regvary\ vectors on cones due 
to Janssen and Drees \cite{janssen:drees:2014} ensures the existence
and finiteness of the limits $\lim_{m\to\infty} m^2p_h$ for $h\ge 1$. This is (AI.1).
\eexam
\ble \label{thm:acf2}
Assume that $(X_t)$ is an  $\mathbb{R}^d$-valued 
strongly mixing strictly stationary \regvary\ sequence with index $\alpha >0$
and that the Borel set $A$ satisfies the
conditions of Section~\ref{subsec:extremogram}. We also assume the
\asy\ independence condition {\rm (AI)} and the mixing condition
\beam\label{eq:mix3}
\lim_{n\to \infty}m^2\sum_{h=r_n}^n \xi_h=0\,.
\eeam
Then
\beam
\widetilde{\gamma}_A(h)&\overset{P}{\to}& 0\,,\qquad h\ge 1\,, \nonumber\\
n^{0.5} \big(
\widetilde{\gamma}_A(i)-E\wt \gamma_A(i)\big)_{i=1,\ldots,h} 
&\std& (Z_i)_{i=1,\ldots,h} \,,\label{eq:cafg}
\eeam
where $(Z_i)_{i=1,\ldots,h}$ are independent Gaussian with mean zero and
variances
\beao
\var(Z_i)=\lim_{m\to\infty} m^2\,p_i\,,\quad i\ge 1\,.
\eeao
\ele
\begin{proof} 
We will apply the central limit theorem in Rio
  \cite{rio:1995} for 
strongly mixing triangular arrays to the \lhs\ in \eqref{eq:cafg}.
For this reason, we have to calculate the \asy\ covariance matrix of
the left-hand vector. 
We observe that for fixed $j>i\ge 1$, in view of
the mixing condition \eqref{eq:mix3} as $\nto$,
\beam\label{eqwm}
\lefteqn{\cov(n^{0.5}\wt \gamma_A(i),n^{0.5}\wt \gamma_A(j))}\nonumber\\&=&
m^2 \cov(\wt I_0\wt I_j,\wt I_0\wt I_i)+
m^2\sum_{t=1}^n \big[\cov (\wt I_0\wt I_i, \wt I_t\wt I_{t+j})+
\cov (\wt I_0\wt I_j, \wt I_t\wt I_{t+i})\big]+o(1)\nonumber\\
&=&m^2 \cov(\wt I_0\wt I_j,\wt I_0\wt I_i)+
m^2\sum_{t=1}^{r_n} \big[\cov (\wt I_0\wt I_i, \wt I_t\wt I_{t+j})+
\cov (\wt I_0\wt I_j, \wt I_t\wt I_{t+i})\big]+o(1)\,.\nonumber\\
\eeam
Condition (AI) implies that 
$m^2 \cov(\wt I_0\wt I_j,\wt I_0\wt I_i)\to 0$ as $m\to\infty$. The
same argument also implies that the first $j$ 
summands in \eqref{eqwm} 
vanish as $\nto$. Therefore it suffices to consider
\beao
m^2\sum_{t=j+1}^{r_n} \big[\cov (\wt I_0\wt I_i, \wt I_t\wt I_{t+j})+
\cov (\wt I_0\wt I_j, \wt I_t\wt I_{t+i})\big]\,.
\eeao
In the latter sum, the indices $0,i,t,t+j$ are distinct and the same
observation applies to $0,j,t,t+i$. Direct 
calculation with condition (AI)
shows that this sum is \asy ally negligible. This implies that the
covariance matrix of the limiting vector is diagonal. 
The calculation of the \asy\ variances is similar by observing that as $\nto$,
\beao
\lefteqn{\var(n^{0.5}\wt\gamma_A(h))}\\
&=&m^2 \var(\wt I_0\wt I_h)+
2\,m^2\sum_{t=1}^n \cov (\wt I_0\wt I_i, \wt I_t\wt I_{t+i})+o(1)=m^2 p_h +o(1)\,.
\eeao
\end{proof}
\bre\label{rem:null}
Although $\wt\gamma_A(h)\stp 0$, $h\ge 1$, it is in general not possible to
avoid centering in \eqref{eq:cafg}. However, under (AI.1), 
$n^{0.5}E\wt\gamma_A(h)\to 0$ if $n/m^2=o(1)$ as $\nto$, and the
latter condition can even be weakened if $m^2(p_h-p_0^2)\to 0$ as $m\to \infty$.
\ere
}

\subsection{Mean square consistency of the integrated periodogram}\label{sec:3}
Recall the definitions of $J_{n,A}(g)$ and $J_A(g)$ for $g\in
L^2_+(\Pi)$ from \eqref{eq:intper}  and \eqref{eq:intpera},
respectively.
The following elementary result deals with the \con\ of  the first and second
moments of $J_{n,A}(g)$ for a given \fct\ $g$. 
\ble \label{thm:app}
Consider an $\mathbb{R}^d$-valued 
strictly stationary \regvary\ sequence $(X_t)$ with index $\alpha
>0$. Assume that
the Borel set $A\subset \ov \bbr_0^d$ satisfies the
conditions of Section~\ref{subsec:extremogram}, $\sum_{l=1}^\infty
\gamma_A(l)<\infty$ and {\rm (M)} holds. Then the following \asy\
relations hold for  $g\in L^2_+(\Pi)$.
\begin{enumerate}
\item[\rm 1.]
$EJ_{n,A}(g)\to J_A(g)$ as $\nto$.
\item[\rm 2.] 
If in addition, $m \log^2 n /n=O(1)$ as $\nto$, and there exists a constant $c>0$ \st
\begin{eqnarray} \label{eq:fcoef}
 |c_h(g)|\le c/h\,,\quad h\ge 1\,,
\end{eqnarray} 
then $E(J_{n,A}(g)- J_A(g))^2\to 0$ and 
  $J_{n,A}^{\circ}(g)\overset{P}{\to} J_A^\circ (g)$ as $\nto$.
\end{enumerate}
\ele
The proof of the lemma is given in Section~\ref{sec:lemma}.
\bre 
Condition \eqref{eq:fcoef} holds under mild smoothness conditions on
$g$, e.g. if
$g$ is Lipschitz or has bounded variation on $\Pi$; see Theorem 4.7
on p. 46
and Theorem 4.12 on p. 47 
in Zygmund~\cite{zygmund:2002}. 
\ere

\section{Functional central limit theorem for the integrated \per }\label{sec:fclt}%\setcounter{equation}{0}
Recall the definition of the spectral density $h_A$ from
Section~\ref{subsec:spectral}.
In this section, we assume that the weight \fct\ $g$ is a non-negative 
continuous \fct . Abusing notation, 
we define the empirical spectral
\df\ with weight \fct\ $g$ by
\beam\label{eq:cont}
J_{n,A}(x)=J_{n,A}(g I_{[0, x]})= \int_{0}^x  I_{n,A}(\la)\, g(\la)\,
d\la\,,\quad
x\in \Pi \,.
\eeam
Under the conditions of Lemma~\ref{thm:app}, again abusing notation,
we have
\beao
J_{n,A}(x)\stp J_A(x)=J_A(g I_{[0, x]} )=\int_0^x h_A(\la)\, g(\la)\,d\la\,,\quad x\in
\Pi \,.
\eeao
In view of the monotonicity and continuity of the \fct s 
$J_{n,A}$ and $J_A$ we also have \beam\label{eq:convp}
\sup_{x\in\Pi} |J_{n,A}(x)-J_A(x)|\stp 0\,.
\eeam
Our next goal is to complement  this consistency result by a \fct al 
\clt \ of the type 
$
(n/m)^{0.5} (J_{n,A}-J_A) \std G\,,
$
in $\bbc(\Pi)$, the space of continuous \fct s on $\Pi$ equipped 
with the uniform topology,   for a suitable Gaussian limit process $G$.
\par
However, this result is unlikely to hold in general, due to 
\asy\ bias problems. It is mentioned in Davis and Mikosch
\cite{davis:mikosch:2009} in relation with the \clt\ for the sample
extremogram (see Lemma~\ref{thm:acf} above) that the pre-\asy\
centerings $E\wt \gamma_A(i)= ((n-i)/n)m(p_i-p_0^2)$ can in general not be replaced  
by their limits $\gamma_A(i)$ due to the failure of the relation $(n/m)^{0.5}|
m(p_i-p_0^2)-\gamma_A(i)|\to 0$ as $\nto$. Therefore we will
equip the empirical  spectral \df\ $J_{n,A}$ with the pre-\asy\ 
centering $E J_{n,A}$. It follows from Lemma~\ref{thm:app}
that under (M), $E J_{n,A}(x)\to J_A(x)$ for every $x\in \Pi$, and
again using monotonicity of $E J_{n,A}$ and $J_A $, we have 
$\sup_{x\in \Pi}|E J_{n,A}(x)- J_A(x)|\to 0$. 
\par
We observe that
\beao
J_{n,A}(x) &=& \psi_0(x) \,\wt \gamma_A(0)
+ 2 \sum_{h=1}^{n-1} \psi_h(x) \, \wt \gamma_A(h)\,,\\
J_{n,A}^\circ(x)&=& \psi_0(x) + 2\sum_{h=1}^{n-1} \psi_h(x)\,\wt \rho_A(h)\,,
\eeao
where
$
\psi_h(x)= \int_{0}^x \cos(h\la) \,g(\la)\,d\la\,,\quad
x\in \Pi\,.
$
We also consider a Riemann sum  approximation of the coefficients
$\psi_h(x)$ at the Fourier \fre ies  $\w_n(i)=2i\pi/n \in \Pi$ given by 
\beao
\widehat{\psi}_h(x)&=& \dfrac{2\pi}{n} \sum_{i=1}^{x_n}
g(\w_n(i)) \cos(h\w_n(i))\,,\quad x\in \Pi\,,
\eeao 
where $x_n=[nx/2\pi ]$.
The corresponding analogs of $J_{n,A}$ and $J_{n,A}^\circ$ are then given by
\beao%\label{eq:discr}
\widehat{J}_{n,A}(x)&= &\widehat{\psi}_0(x)\wt \gamma_A(0) + 2
\sum_{h=1}^{n-1} \widehat{\psi}_h(x)  \,\wt \gamma_A(h)\,,\\
\widehat{J}_{n,A}^\circ(x)&= &\widehat{\psi}_0(x) + 2
\sum_{h=1}^{n-1} \widehat{\psi}_h(x)  \,\wt \rho_A(h)\,,
\eeao
\par
Now we are ready to formulate the main result of this paper.
\bth\label{thm:main}
Assume that $(X_t)$ is an  $\mathbb{R}^d$-valued 
strictly stationary \regvary\ sequence with index $\alpha >0$
and the Borel set $A\subset \ov \bbr_0^d$ is bounded away from zero, 
$\mu_1(\partial A)=0$ and $\mu_1(A)>0$. Let
$g$ be a non-negative $\beta$-H\"{o}lder continuous \fct\ with $\beta\in (3/4,1]$.
If the conditions %{\rm (M)}, 
{\rm (M1)} and $\sum_{l=1}^\infty
\gamma_A(l)<\infty$ hold then in $\bbc(\Pi)$,
\beam \label{eq:fclt1}
(n/m)^{0.5} (J_{n,A} - EJ_{n,A}) &\std&  G\,,\quad \nto\,,\\
 \label{eq:fclt2}
(n/m)^{0.5} (\widehat{J}_{n,A}-E\widehat{J}_{n,A}) &\std & G\,, \quad
\nto\,,
\eeam
where the limit process is given by the infinite series
\beam
G=\psi_0  Z_0 +2 \sum_{h=1}^{\infty} \psi_h\, Z_h\,,
\eeam
which converges in \ds\ in $\bbc(\Pi)$,
$(Z_h)$ is a mean zero Gaussian \seq\ such that $(Z_0,\ldots,Z_h)$ has
the covariance  
matrix $(\Sigma_h)$, $h\ge 0$, given in Lemma~\ref{thm:acf}. Moreover,
the following limit relations hold
\beam
 \label{eq:fclt01}
(n/m)^{0.5} \big(J_{n,A}^{\circ} - EJ_{n,A}/(mp_0)\big) &\std&  G^\circ\,,\quad \nto\,,\\
 \label{eq:fclt02}
(n/m)^{0.5} \big(\widehat{J}_{n,A}^{\circ}-E\widehat{J}_{n,A}/(mp_0)\big) &\std & G^\circ \,, \quad
\nto\,,
\eeam\label{eq:g1}
where the limit process is given by the infinite series
\beao
G^\circ= \dfrac{ 2}{ \gamma_A(0)}\sum_{h=1}^\infty \psi_h (Z_h
-\rho_A(h) Z_0)\,.
\eeao
\ethe
The proof of this result is given in Section~\ref{sec:proof1}.
\bre \label{rem:main1} For practical purposes, the discretized version 
$\wh J_{n,A}$ will be preferred to $J_{n,A}$ since it does not involve the
calculation of integrals. 
Moreover, since $\sum_{t=1}^n \ex^{i \w_n(j)t}=0$ 
for $ \w_n(j)\in (0,\pi)$, centering of the indicators $I_t$ with
the unknown parameter  $p_0$
in the \per\ ordinates $I_{n,A}(\w_n(j))=(m/n)|\sum_{t=1}^n I_t \ex^{i\w_n(j) t}|^2$ is not needed.
\ere
For an $\eta$-dependent \seq\ $(X_t)$, we know that $Z_h=0$ a.s. for
$h>\eta$. Then we conclude from Theorem~\ref{thm:main} and Lemma~\ref{thm:acf} that the
limit process $G$ collapses into
$
G=  \psi_0 Z_0+2\sum_{h=1}^{\eta}\psi_h Z_h\,.
$
However, taking into account  
Lemma~\ref{thm:mdep1}, a more sophisticated result with a different
\con\ rate can be derived. The 
corresponding result for $J_{n,A}^\circ$ is similar and  therefore omitted.
\bth\label{thm:mdep2}
Assume that $(X_t)$ is an  $\mathbb{R}^d$-valued strongly mixing strictly stationary
$\eta$-dependent \regvary\ sequence with index $\alpha >0$ for some
$\eta\ge 0$
and the Borel set $A\subset\ov \bbr_0^d$ is bounded away from zero, 
$\mu_1(\partial A)=0$ and $\mu_1(A)>0$.
Also assume that the
limits in \eqref{eq:mdep1} exist. Let
$g$ be a non-negative
$\beta$-H\"{o}lder continuous \fct\ with $\beta\in (3/4,1]$.
Then the relations
\beao %\label{eq:iid21}
\sqrt{n}\big(J_{n,A} -\psi_0\widetilde{\gamma}_A(0)-2 \sum_{h=1}^{\eta}
\psi_h\widetilde{\gamma}_A(h)\big) &\std&  \overline{G}\,,\nonumber\\
 %\label{eq:iid22}
\sqrt{n} \big(\widehat{J}_{n,A}-\widehat{\psi}_0\widetilde{\gamma}_A(0) -2 \sum_{h=1}^{\eta}
\widehat{\psi}_h\widetilde{\gamma}_A(h)\big) &\std & \overline{G}\,,
\eeao
hold in $\bbc(\Pi)$, 
where the limit process is given by the a.s. converging infinite series 
\begin{eqnarray*}
 \overline{G}=2 \sum_{h=1}^{\infty} \psi_{\eta+h}\, Z_h\,,
\end{eqnarray*} and $(Z_h)$ is a mean zero Gaussian \seq\ such that
$(Z_1,\ldots,Z_h)$ has covariance matrix
$\overline{\Sigma}_h$, $h\ge 1$, given in Lemma~\ref{thm:mdep1}.
\ethe
The proof is given in Section~\ref{sec:proof2}.
\bexam\label{exam:bb}\rm
Assume that $(X_t)$ is an iid regularly varying sequence with index
 $\alpha>0$. Then $(Z_h)$ is an iid mean zero Gaussian sequence 
with $\var(Z)=\gamma_A^2(0)=(\mu_1(A))^2$. 
If we choose the \fct\ $g\equiv 1$ we obtain
\beao
\psi_h(x)= \int _0^x \cos(h\la) d\la = \dfrac{\sin(hx)}{h}\,,\quad h\ge 0\,,\quad x\in \Pi\,,
\eeao
and 
\beao
\overline{G}(x)=2 \sum_{h=1}^{\infty}  \dfrac{\sin(hx)}{h}\,Z_h\,,\quad x\in\Pi\,.
\eeao
We notice that $\ov G$ is a series \rep\ of a Brownian bridge; 
see Hida \cite{hida:1980}.
\eexam 
{\red In the case of \asy\ (extremal) independence a result similar to
  Theorem~\ref{thm:mdep2} holds.
\bth\label{thm:mdep2a}
Assume that $(X_t)$ is an  $\mathbb{R}^d$-valued strictly stationary  
 \regvary\ sequence with index $\alpha >0$ 
and the Borel set $A\subset\ov \bbr_0^d$ is bounded away from zero and
the axes, 
$\mu_1(\partial A)=0$ and $\mu_1(A)>0$. Also assume the mixing
condition \eqref{eq:mix3} and the \asy\ independence condition {\rm (AI).}
Let
$g$ be a non-negative
$\beta$-H\"{o}lder continuous \fct\ with $\beta\in (3/4,1]$.
Then the relations
\beao %\label{eq:iid21}
\sqrt{n}\big((J_{n,A} -EJ_{n,A})-\psi_0(\widetilde{\gamma}_A(0)-E\wt \gamma_A(0)) \big) &\std&  \wh{G}\,,\nonumber\\
 %\label{eq:iid22}
\sqrt{n} \big((\widehat{J}_{n,A}-E\widehat{J}_{n,A})-\widehat{\psi}_0(\widetilde{\gamma}_A(0)-E\widetilde{\gamma}_A(0))\big) &\std & \wh{G}\,,
\eeao
hold in $\bbc(\Pi)$, 
where the limit process is given by the a.s. converging infinite series 
\begin{eqnarray*}
 \wh {G}=2 \sum_{h=1}^{\infty} \psi_{h}\, Z_h\,,
\end{eqnarray*} and $(Z_h)$ is a \seq\ of independent mean zero
Gaussian variables with variances
$\var(Z_h)=\lim_{m\to\infty}m^2\,p_h$, $h\ge  1$. 
\ethe
The proof is based on Lemma~\ref{thm:acf2} and tightness arguments
which are similar to the proofs of Theorem~\ref{thm:main} and
\ref{thm:mdep2}. We omit further details. In view of Remark~\ref{rem:null}, 
centering in Theorem~\ref{thm:mdep2a} can be avoided if $n/m^2=o(1)$
as $\nto$.}
\par
\begin{figure}
\begin{minipage}[ht]{0.33\linewidth}  
\centering  
\includegraphics[width=\textwidth]{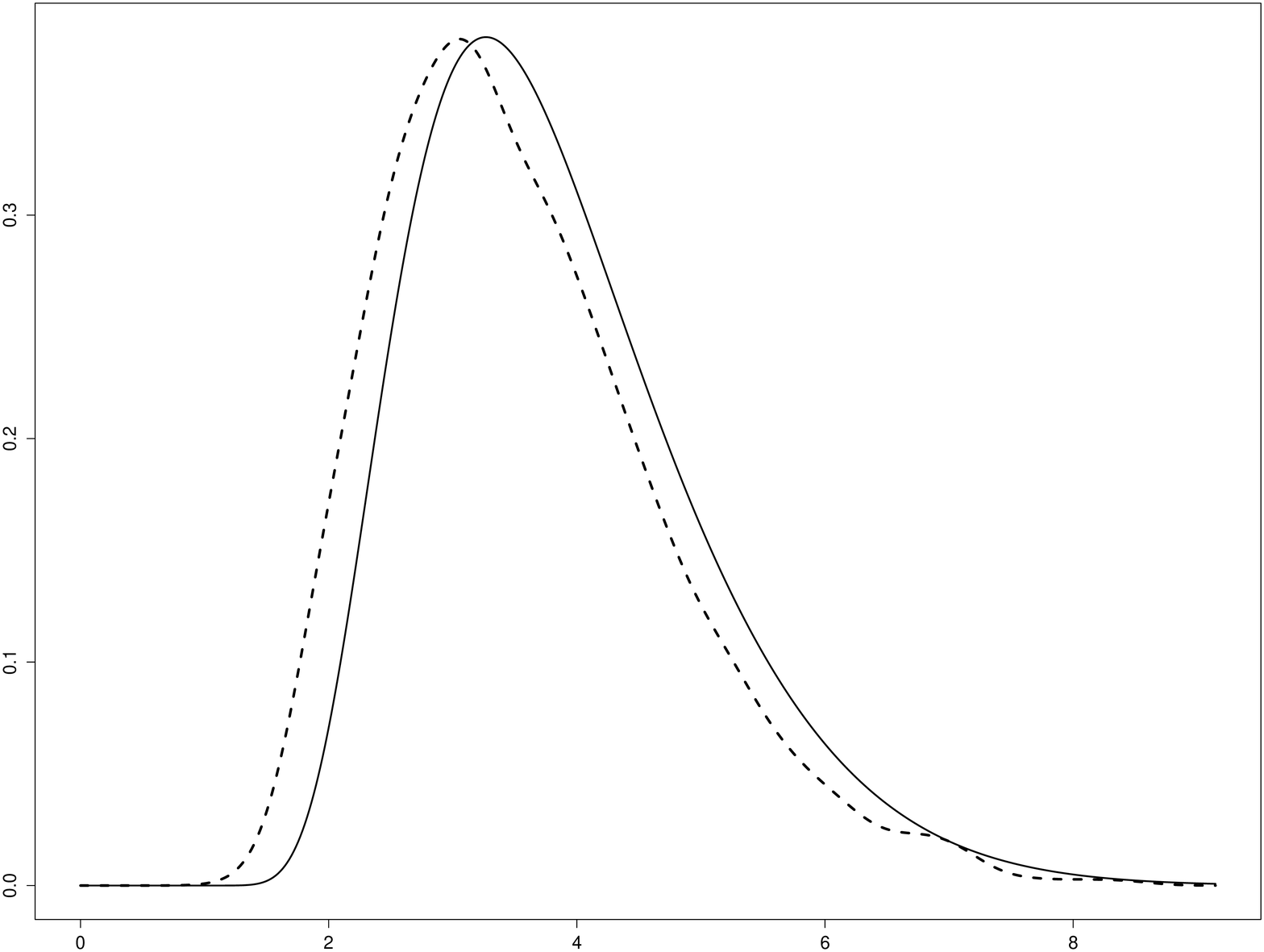}  
\end{minipage}%
\begin{minipage}[ht]{0.33\linewidth}  
\centering  
\includegraphics[width=\textwidth]{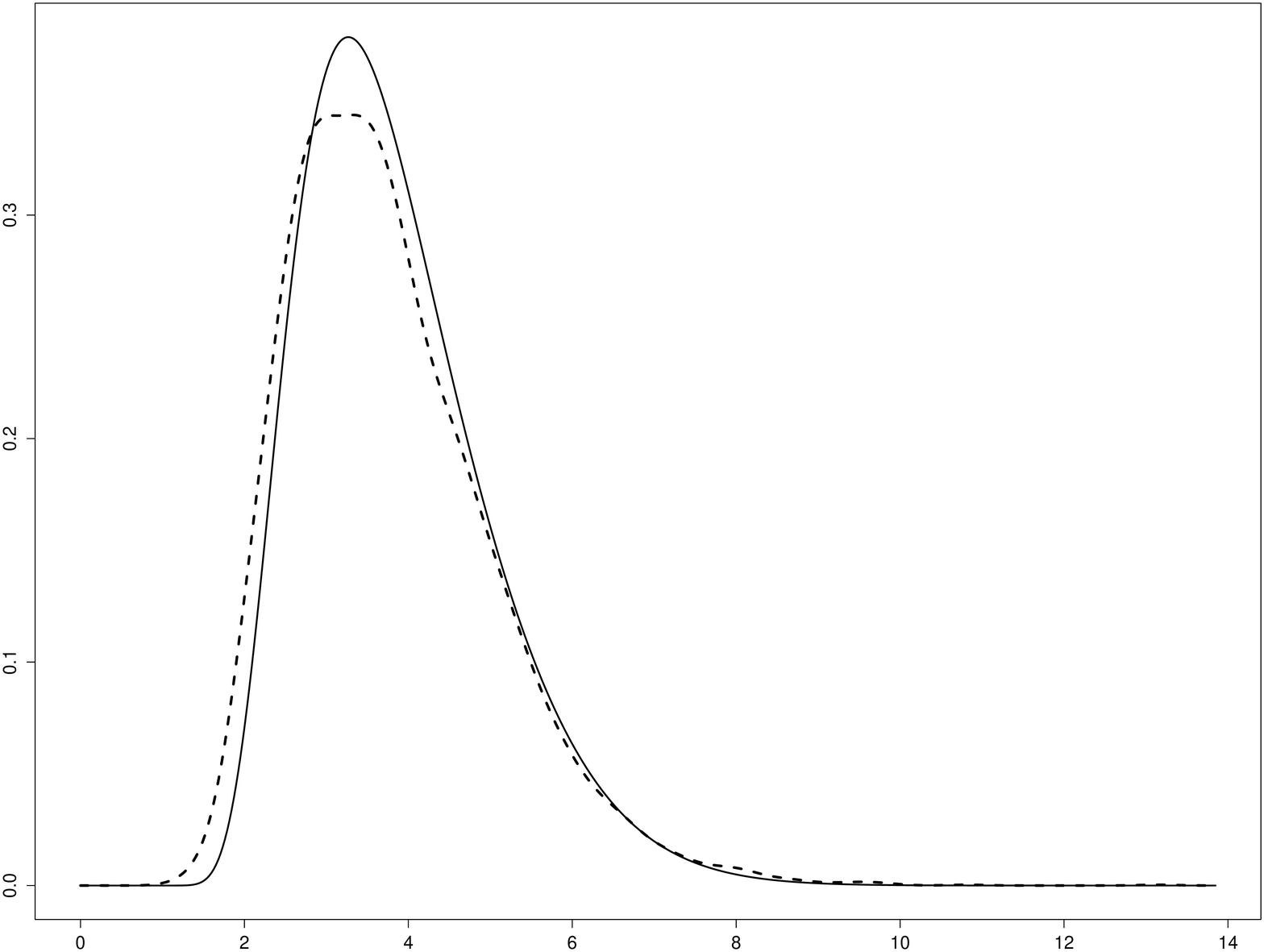}  
\end{minipage}%
\begin{minipage}[ht]{0.33\linewidth}  
\centering  
\includegraphics[width=\textwidth]{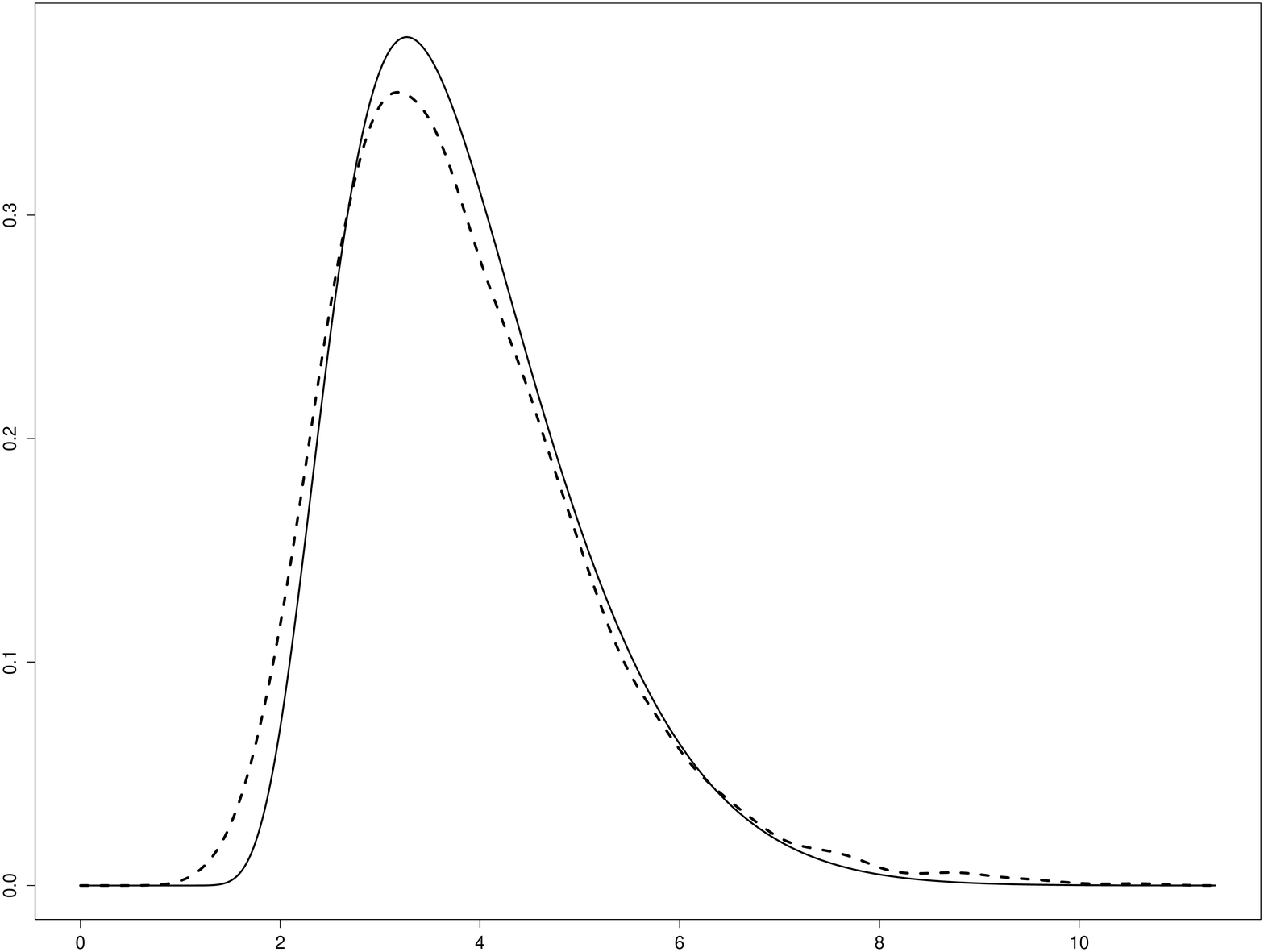}  
\end{minipage}%
\\
\begin{minipage}[ht]{0.33\linewidth}  
\centering  
\includegraphics[width=\textwidth]{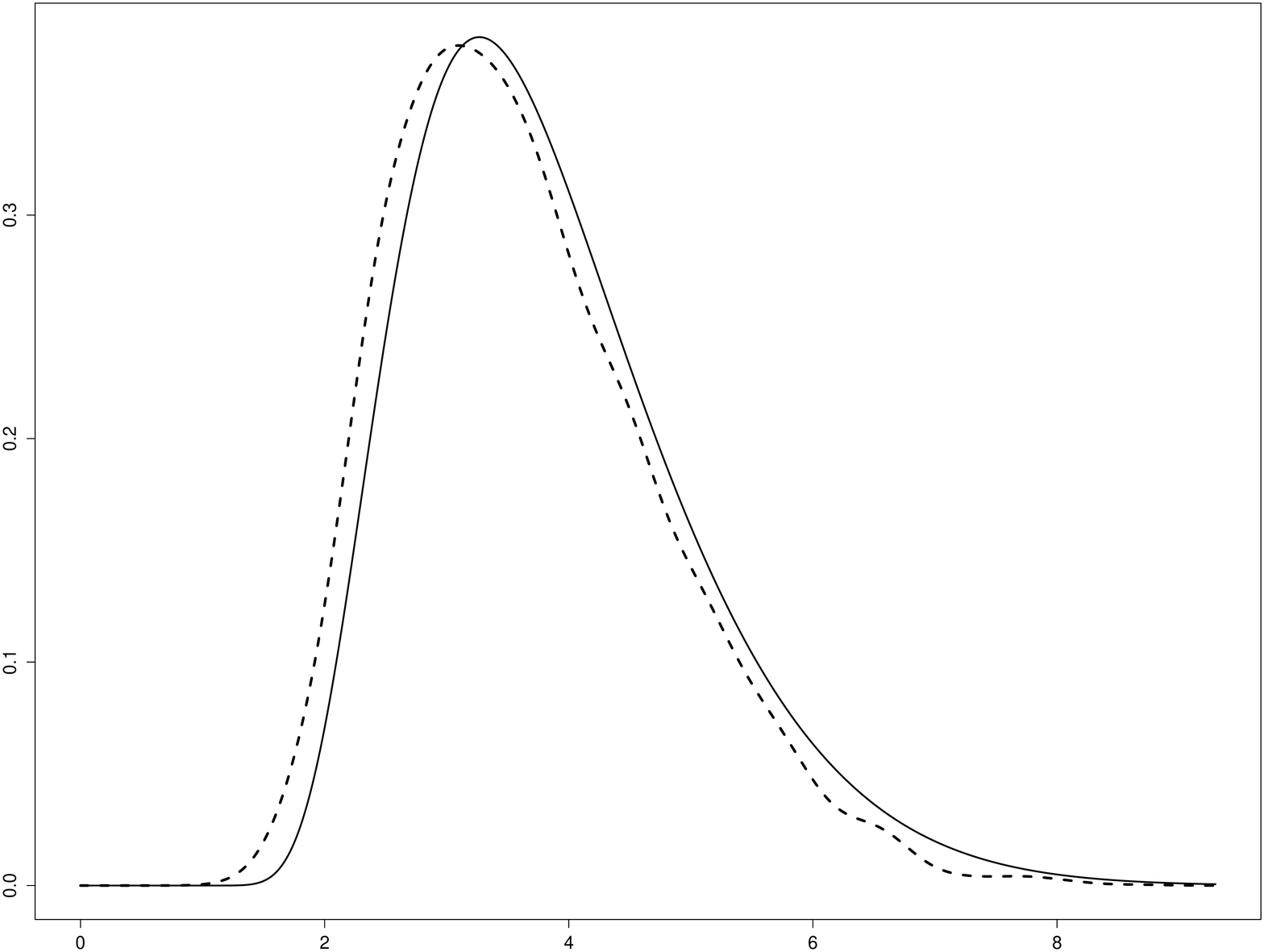}  
\end{minipage}%
\begin{minipage}[ht]{0.33\linewidth}  
\centering  
\includegraphics[width=\textwidth]{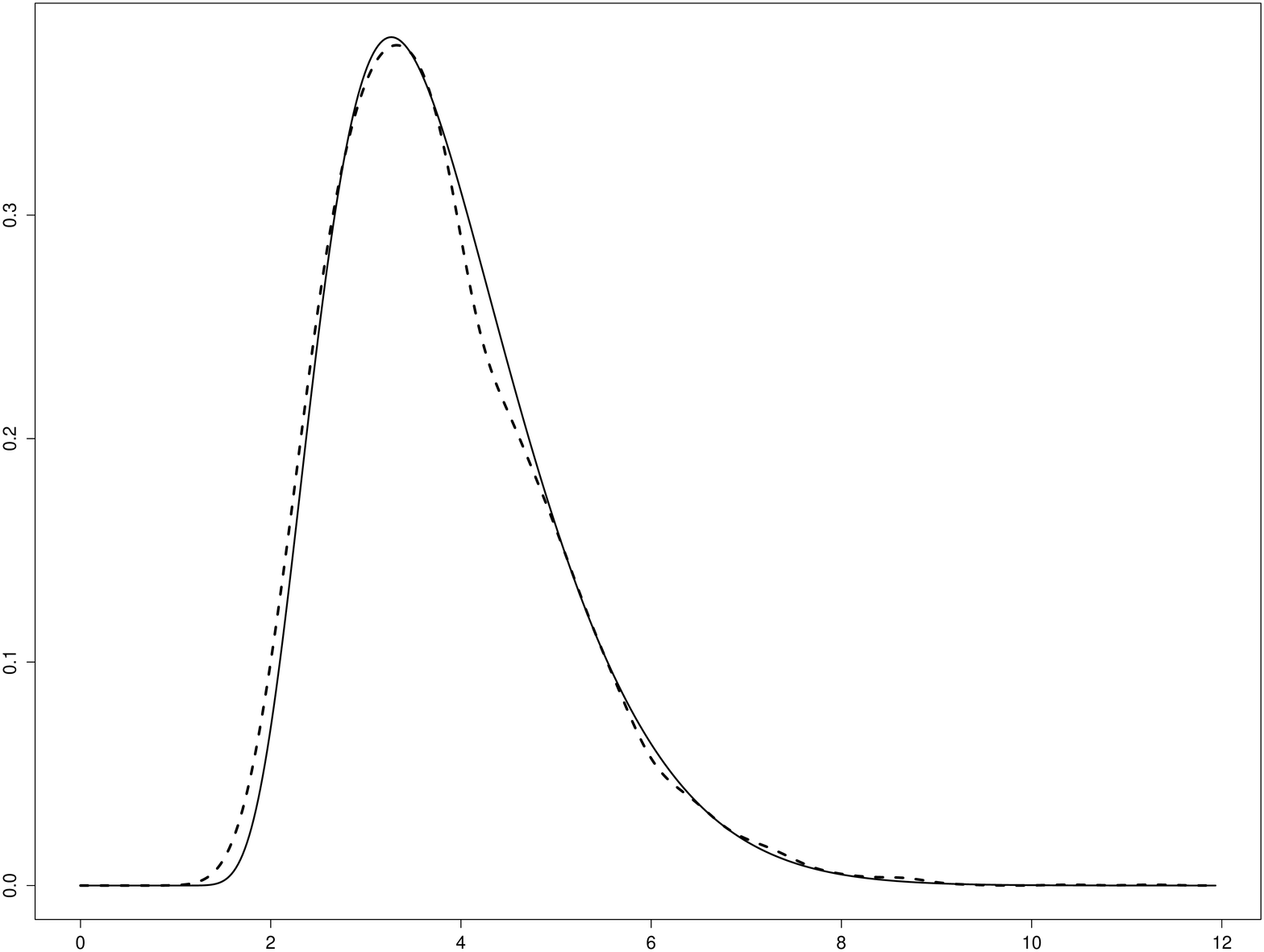}  
\end{minipage}%
\begin{minipage}[ht]{0.33\linewidth}  
\centering  
\includegraphics[width=\textwidth]{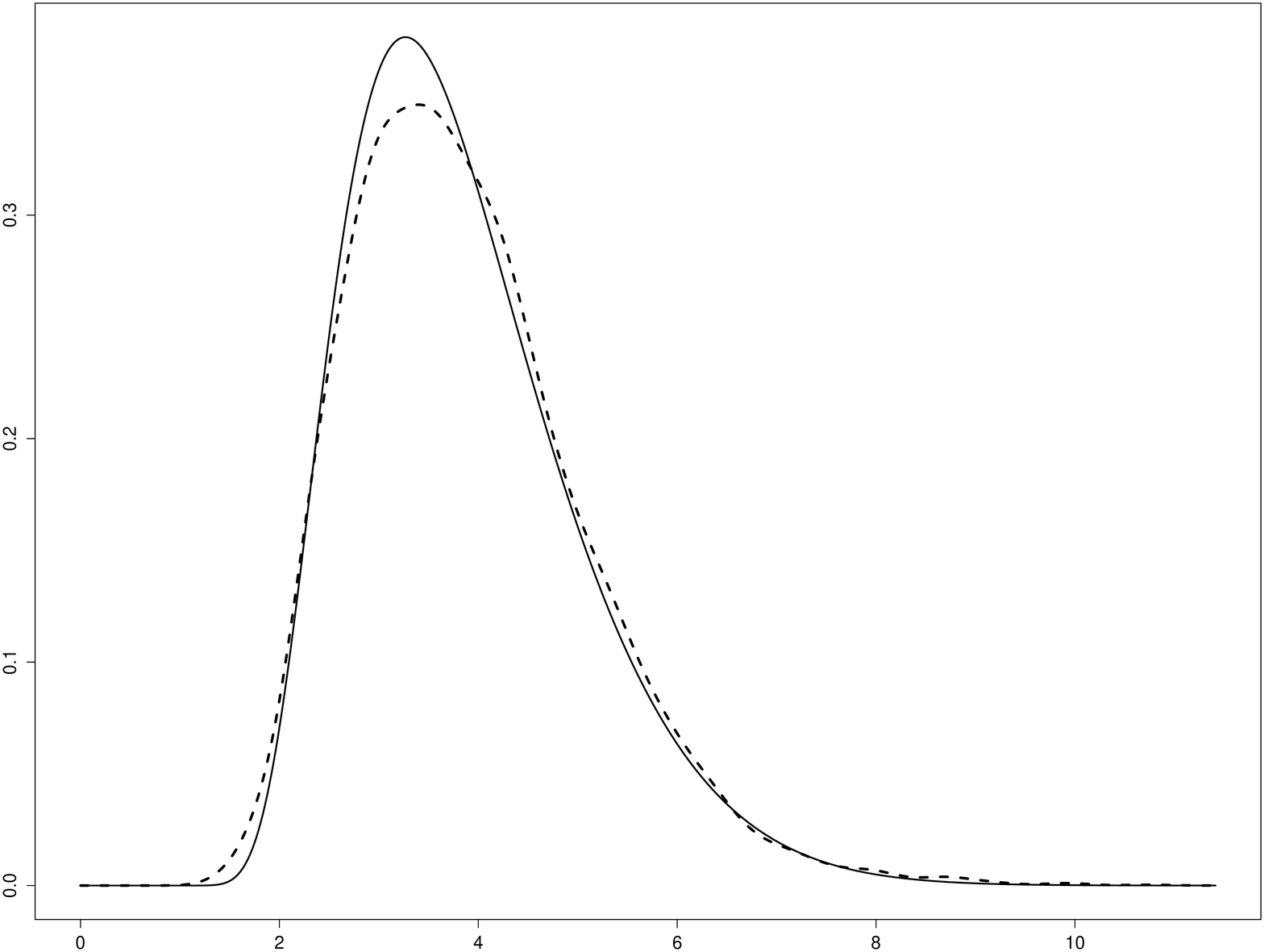}  
\end{minipage}%
\caption{\red Density of the \lhs\ in \eqref{eq:pp1} with $\eta=0$ (dotted line) and 
its limit $\sup_{x \in \Pi}|\ov G(x)|$ (solid line) for $g\equiv 1$. We choose the set 
$A=(1,\infty)$, different thresholds $a_m$ with $p_0=P(X>a_m)$ and
different sample sizes $n$. The underlying \seq\ $(X_t)$ is iid
$t$-distributed with $\alpha=3$ degrees of freedom.
 The sample sizes
are chosen as $n=2,000$ in the first row and $n=10,000$ in the second
row. The thresholds $a_m$ are chosen such that $p_0=0.1$ in the first
column, $p_0=0.05$ in the second column and $p_0=0.03$ in the third column.}
\label{fig:01}
\end{figure}
\begin{figure}
\begin{minipage}[ht]{0.33\linewidth}  
\centering  
\includegraphics[width=\textwidth]{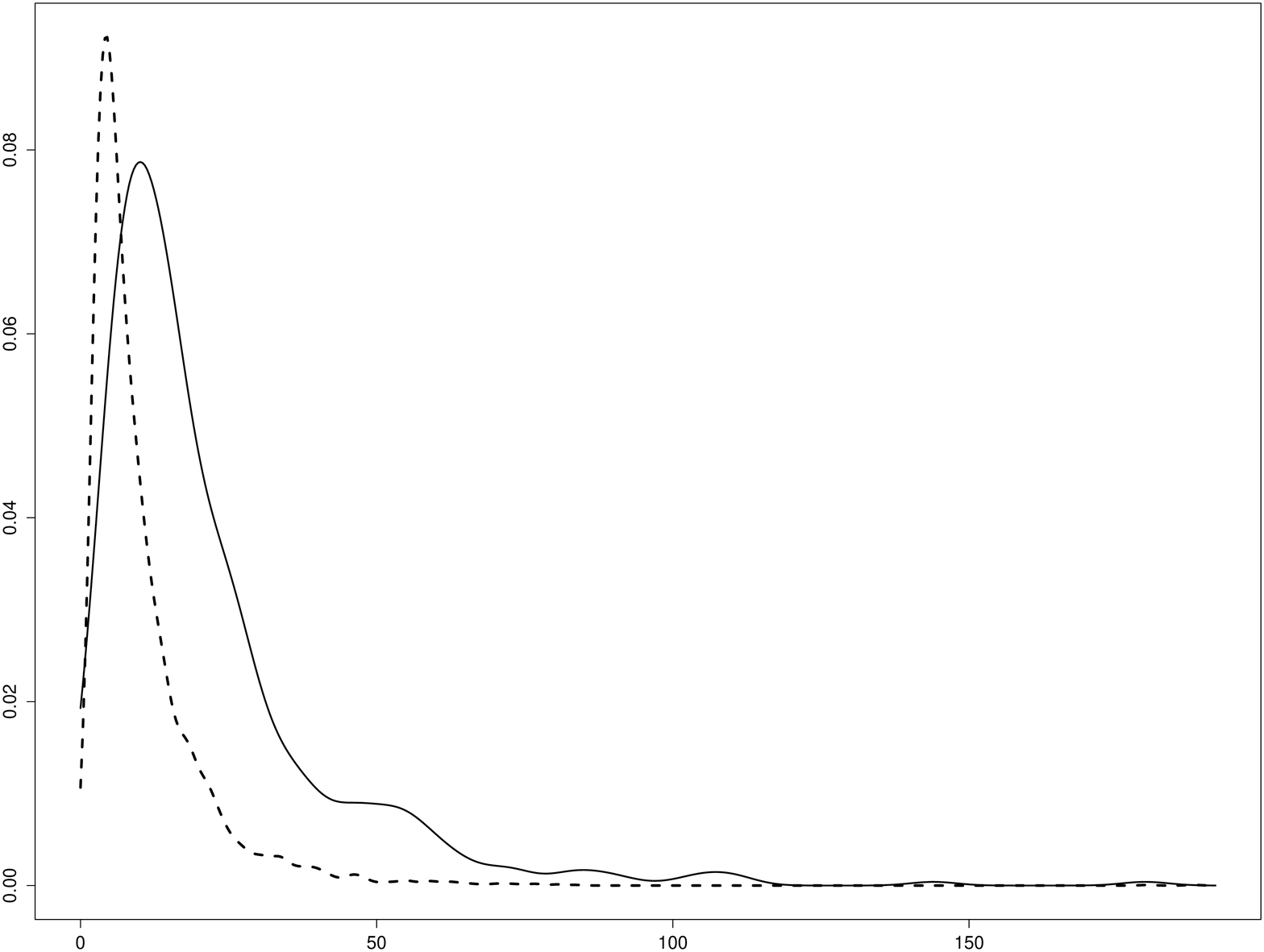}  
\end{minipage}%
\begin{minipage}[ht]{0.33\linewidth}  
\centering  
\includegraphics[width=\textwidth]{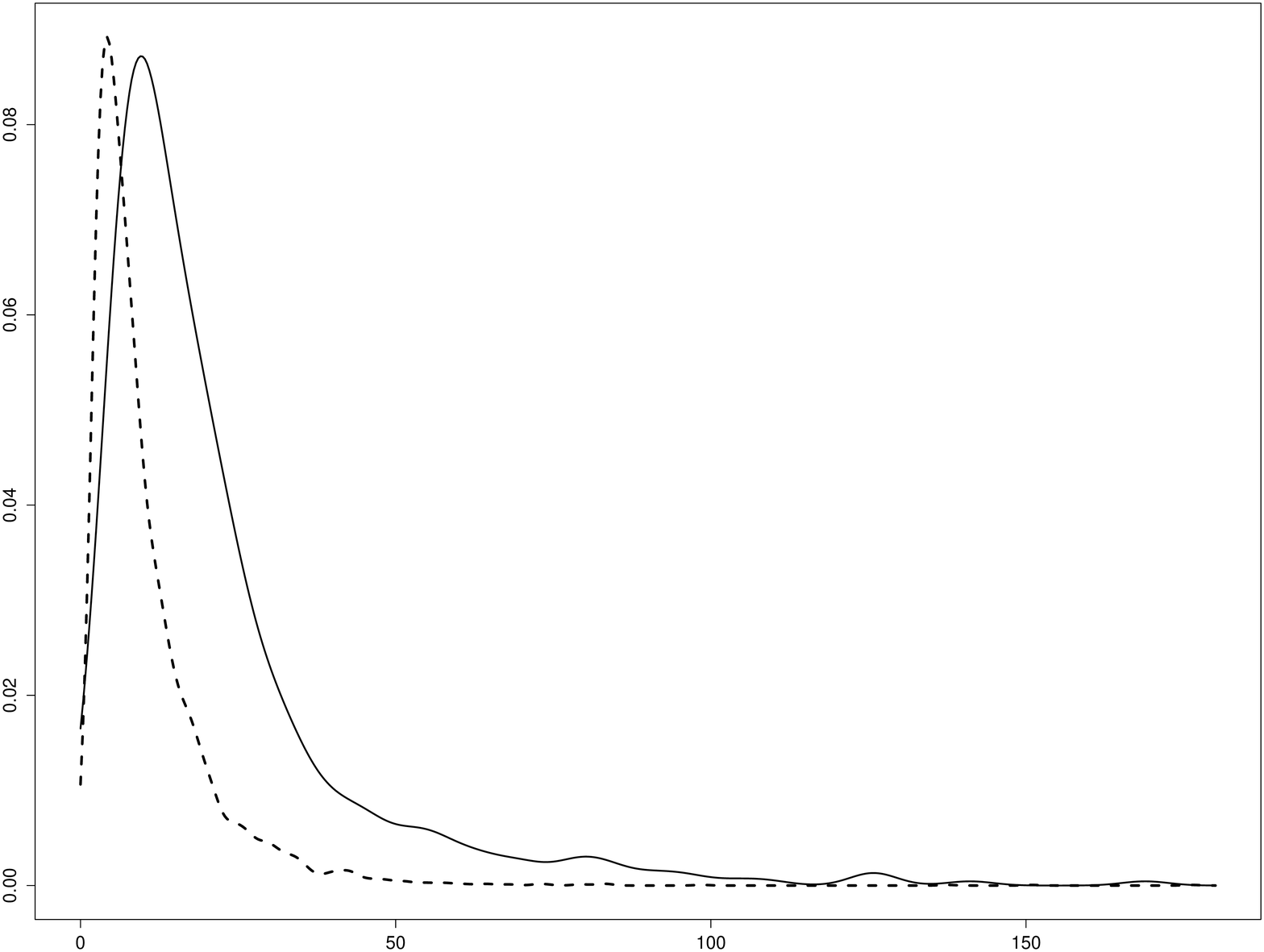}  
\end{minipage}%
\begin{minipage}[ht]{0.33\linewidth}  
\centering  
\includegraphics[width=\textwidth]{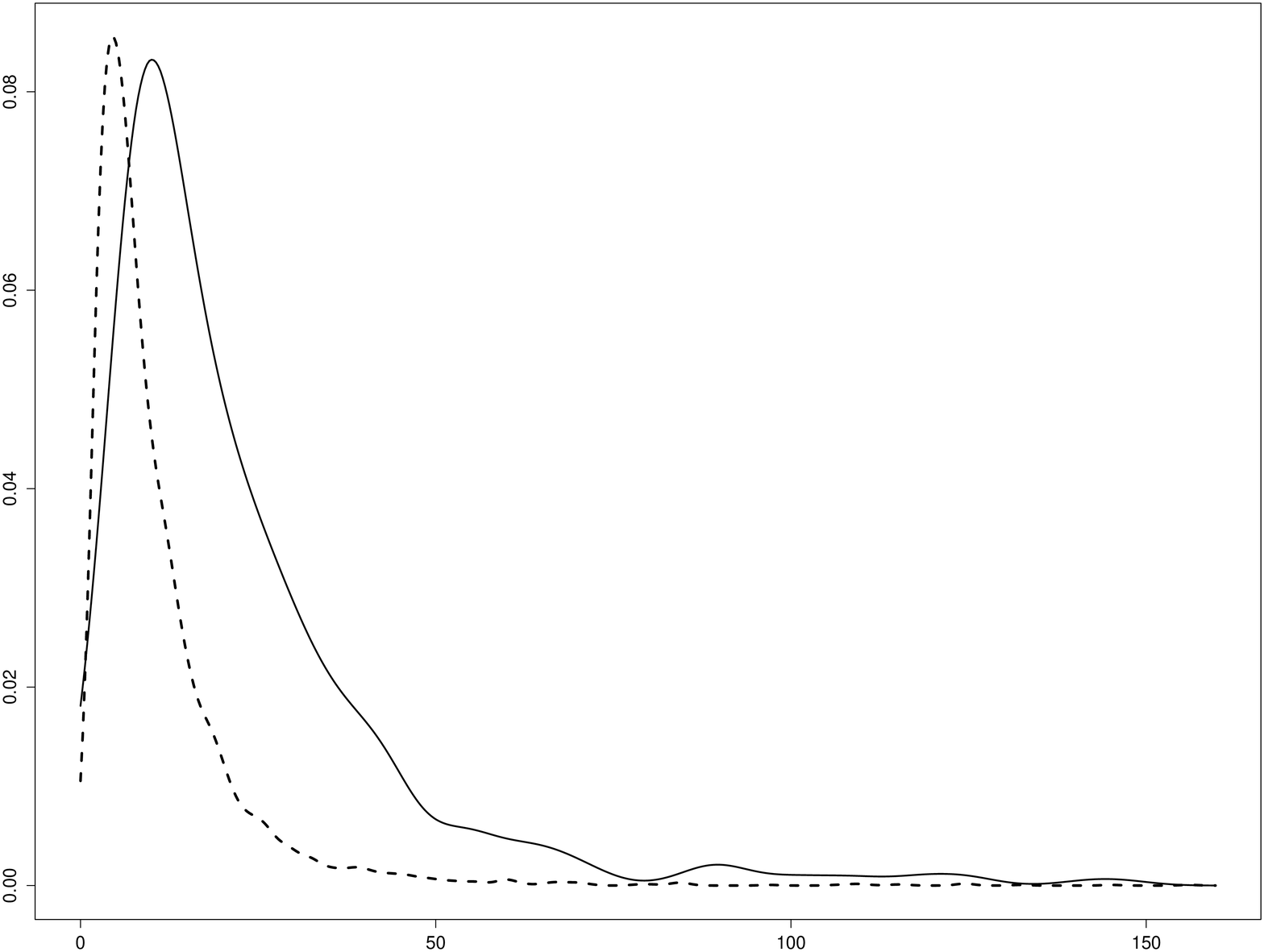}  
\end{minipage}%
\\
\begin{minipage}[ht]{0.33\linewidth}  
\centering  
\includegraphics[width=\textwidth]{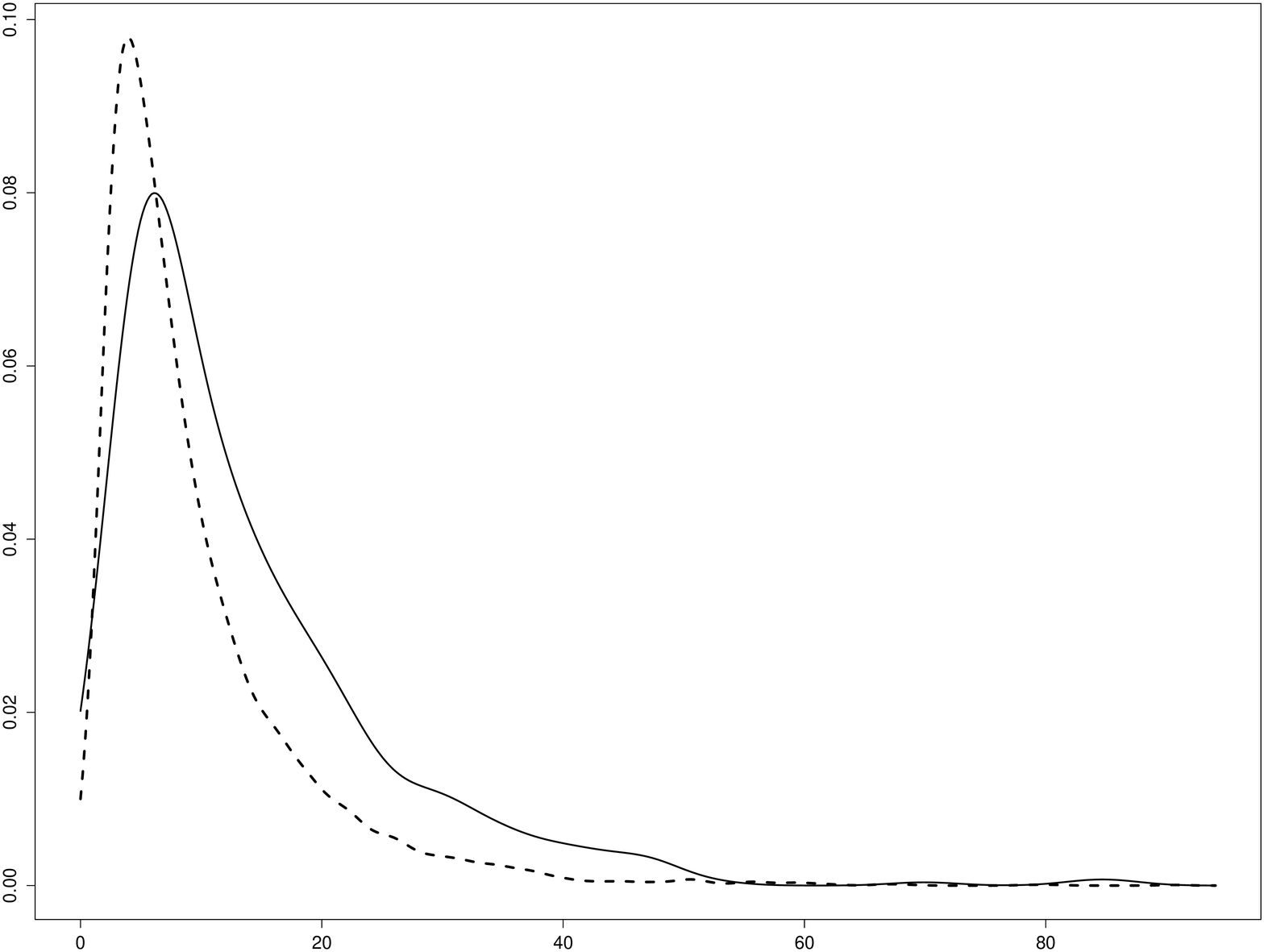}  
\end{minipage}%
\begin{minipage}[ht]{0.33\linewidth}  
\centering  
\includegraphics[width=\textwidth]{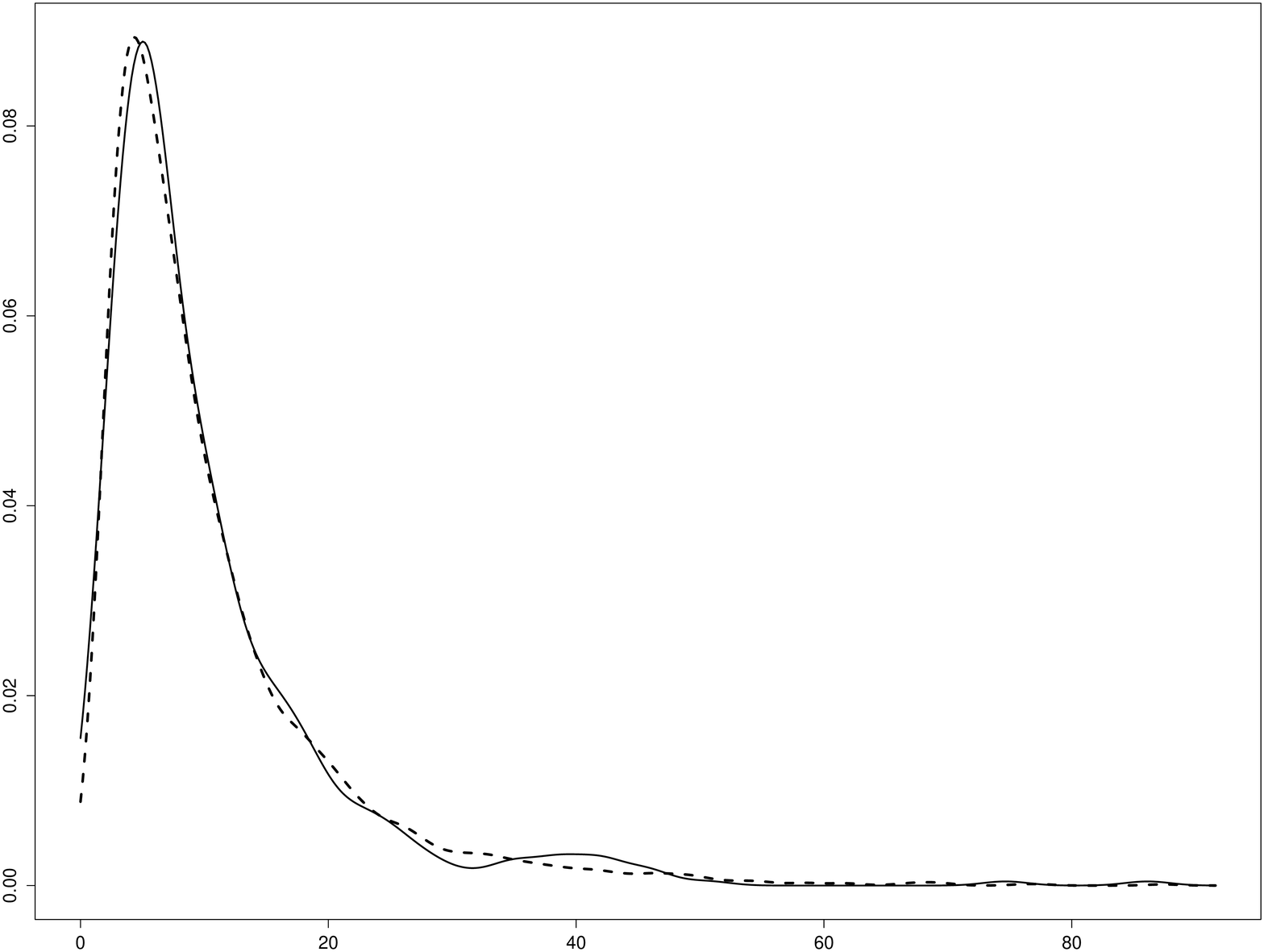}  
\end{minipage}%
\begin{minipage}[ht]{0.33\linewidth}  
\centering  
\includegraphics[width=\textwidth]{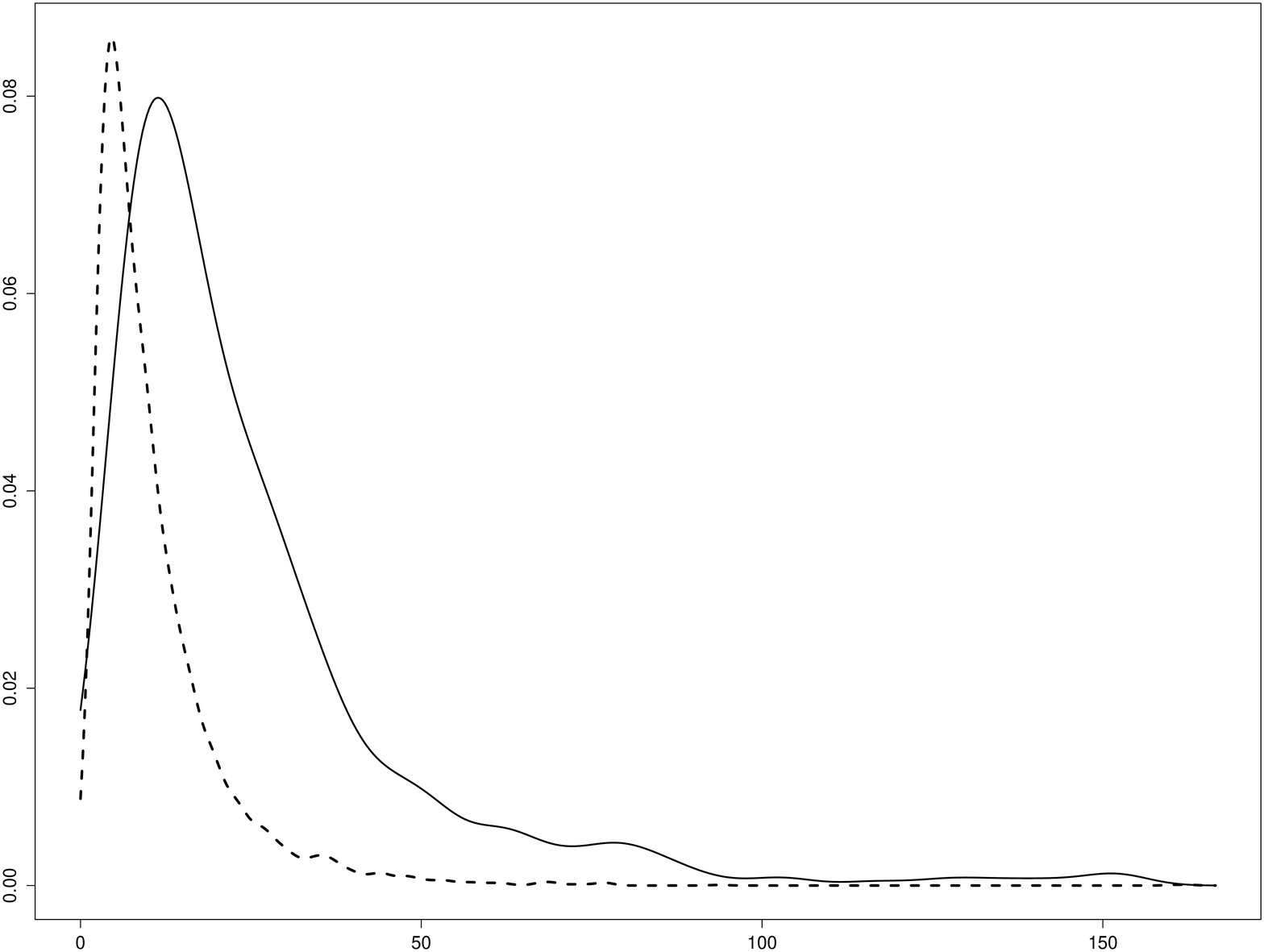}  
\end{minipage}%
\caption{\red Density of the \lhs\ in \eqref{eq:ppr1} with $\eta=0$ (dotted
  line) and its limit
 $\int_{x\in \Pi} \ov G^2(x)\,dx$ (solid line) for $g\equiv
1$. We choose the same setting as in Figure~\ref{fig:02}.}
\label{fig:02}
\end{figure}

As in classical limit theory for the empirical spectral \ds\ (see Gre\-nan\-der and Rosenblatt \cite{grenander:rosenblatt:1984},
Dahlhaus  \cite{dahlhaus:1988}), an application of the \cmt\ to
Theorems~\ref{thm:main} and \ref{thm:mdep2} yields limit theory for  \fct
als of the integrated \per . These \fct als can be used 
for testing the goodness of fit of the spectral density of the \ts\ model
underlying the data, under the null hypothesis that the model is
correct. From Theorem~\ref{thm:main} we get the following limit
results
for the corresponding test statistics.
\begin{itemize}
\item {\em Grenander-Rosenblatt test:}
\beam\label{eq:g-r}
(n/m)^{0.5}\sup_{x\in \Pi}\Big|J_{n,A}(x)-EJ_{n,A}(x)\Big|&\overset{d}{\to} & \sup_{x \in \Pi}|G(x)|\,.
\eeam
\item {\em $\omega^2$- or Cram\'{e}r-von Mises test:}
\begin{eqnarray*}
 (n/m)\int_{x\in \Pi}\Big({J}_{n,A}(x)-E{J}_{n,A}(x)\Big)^2\,dx \overset{d}{\to} \int_{x\in \Pi} G^2(x)\,dx\,.
\end{eqnarray*}
\end{itemize}
If $(X_t)$ is an $\eta$-dependent \seq\ satisfying the conditions of Theorem~\ref{thm:mdep2}, 
the corresponding limit results read as follows:
\begin{itemize}
\item {\em Grenander-Rosenblatt test:}
\beam\label{eq:pp1}
\sqrt{n}\sup_{x\in
  \Pi}\Big|{J}_{n,A}(x)-\psi_0(x)\widetilde{\gamma}_A(0)-2\sum_{h=1}^{\eta}
\psi_h(x) \widetilde{\gamma}_A(h)\Big|\overset{d}{\to}  \sup_{x \in \Pi}|\ov G(x)|\,.\nonumber\\
\eeam
\item $\omega^2$-statistic or Cram\'{e}r-von Mises test:{\small
\beam\label{eq:ppr1}
 n\int_{x\in \Pi}\big({J}_{n,A}(x) -\psi_0(x)\widetilde{\gamma}_A(0) -2\sum_{h=1}^{\eta}
\psi_h(x) \widetilde{\gamma}_A(h)\big)^2\,dx \overset{d}{\to} \int_{x\in \Pi} \ov G^2(x)\,dx\,.\nonumber\\
\eeam }
\end{itemize}
{\red In Figures~\ref{fig:01} and \ref{fig:02} 
we show the estimated densities of the test statistics 
in \eqref{eq:pp1} and \eqref{eq:ppr1} for $n=2,000$ and $n=10,000$,
for different thresholds $a_m$ and $ g\equiv 1$. We compare the
estimated densities with their
corresponding limits. The samples are iid
$t$-distributed with $\alpha=3$ degrees of freedom.  
We mention that
the density of  $\sup_{x \in \Pi}|\ov G(x)|$ is given by $ 4\pi^{-2}
\sum_{j=1}^{\infty} (-1)^{j+1}x\exp\big(-j^2x^2/\pi^2\big)$, $x>0$; 
see Shorack and Wellner \cite{shorack:wellner:1986}.
We use the identity in law $\int_{x\in \Pi} \ov G^2(x)\,dx\eqd
\sum_{j=1}^{\infty}(2/j^2)N_j^2 $ for an iid standard normal \seq\ 
$(N_j)$ (see \cite{shorack:wellner:1986}) for the simulation of the limiting
density of the $\omega^2$-statistic.
\par
Not surprisingly, these graphs show that one needs rather large sample sizes
to make the tests reliable. The Grenander-Rosenblatt statistic
shows a better overall behavior in comparison with the
$\omega^2$-statistic. The \ds\ of the former statistic is close to its
limit for a variety of thresholds like $p_0=0.1, 0.05$ and even for $p_0=0.03$.
In contrast, the $\omega^2$-statistic is rather sensitive to the
choice of threshold and sample size; the best overall approximation
is achieved for $n=10,000$ and $p_0=0.05$. For applications, one would
need to focus on the quality of the approximation of high/low quantiles of 
the test statistics by the limiting quantiles. This task is not
addressed in this paper. 
}
%\begin{figure}[htbp]
%\begin{center}
%\includegraphics[width=0.5\textwidth]{Fig2.eps}
%\caption{Density of the \lhs\ in \eqref{eq:pp1} with $\eta=0$ (dotted line) and 
%its limit $\sup_{x \in \Pi}|\ov G(x)|$ (solid line) for $g\equiv 1$. We choose the set 
%$A=(1,\infty)$, the threshold $a_m$ \st\ $p_0=P(X>a_m)=0.05$ and the sample size $n=10,000$. The underlying \seq\ $(X_t)$ is iid
%$t$-distributed with $\alpha=3$ degrees of freedom.
%The density of  $\sup_{x \in \Pi}|\ov G(x)|$ is given by  $ 4\pi^{-2}
% \sum_{j=1}^{\infty} (-1)^{j+1}x\exp\big(-j^2x^2/\pi^2\big)$, $x>0$; 
%see Shorack and Wellner \cite{shorack:wellner:1986}.}\label{fig:2}
%\end{center}
%\end{figure}

%\begin{figure}[t]
%\begin{center}
%\includegraphics[width=0.5\textwidth]{Fig3.eps}
%\caption{Density of the \lhs\ in \eqref{eq:ppr1} with $\eta=0$ (dotted
%  line) and its limit
% $\int_{x\in \Pi} \ov G^2(x)\,dx$ (solid line) for $g\equiv
%1$. We choose the same setting as in Figure~\ref{fig:2}.
%We use the identity in law $\int_{x\in \Pi} \ov G^2(x)\,dx\eqd
%\sum_{j=1}^{\infty}(2/j^2)Z_j^2 $ for an iid standard normal \seq\ 
%$(Z_j)$ (see \cite{shorack:wellner:1986}) for the simulation of the limiting
%\rv .}\label{fig:3}
%\end{center}
%\end{figure}

\section{The bootstrapped integrated periodogram}\label{sec:bootstrap}
With a few exceptions, the limit processes $G$ and $\overline{G}$ in
Theorem~\ref{thm:main} and ~\ref{thm:mdep2} have an unfamiliar
dependence structure and then it is impossible to give confidence
bands for the test statistics mentioned in the previous section. 
One faces a similar problem when dealing with the sample extremograms
whose asymptotic covariance matrix is a complicated function of the
measures $\mu_h$ in \eqref{eq:2}. Davis et
al.~\cite{davis:mikosch:cribben:2012} proposed to apply the stationary
bootstrap for constructing confidence bands for the sample
extremogram. The stationary bootstrap can also be used for the
integrated periodogram, as we will show below. 

\subsection{Stationary bootstrap}\label{sec:bp}
The stationary bootstrap was introduced by Politis and
Romano~\cite{politis:romano:1994} as an alternative block bootstrap
method. First, we describe this procedure for a strictly stationary
sequence $(Y_t)$. Given a sample $Y_1,\ldots,Y_n$, consider the
bootstrapped sequence
\begin{eqnarray}\label{eq:bp0}
Y_{K_1},\ldots, Y_{K_1+L_1-1},\ldots , Y_{K_N},\ldots ,Y_{K_N+L_N-1},\ldots\,,
\end{eqnarray}
where $(Y_i)$, $(K_i)$, $(L_i)$ are independent sequences, $(K_i)$ is
an iid sequence of random variables uniformly distributed on
$\{1,\ldots,n\} $, $(L_i)$ is an iid sequence of geometrically 
distributed random variables with distribution
$P(L_1=i)=\theta(1-\theta)^{i-1}$, $i=1,2,\ldots,$ for some
$\theta=\theta_n\in (0,1)$ such that $\theta_n\to 0$ as $\nto$ and 
$N=N_n=\inf\{i\ge 1: \sum_{j=1}^i L_j\ge n \}$. If any element $Y_t$
in \eqref{eq:bp0} has an index $t>n$, we replace it by $Y_{t\mod
  n}$. As a matter of fact, $(Y_t)_{t\ge 1}$ constitutes a strictly
stationary sequence. The stationary bootstrap sample is now chosen as
the block of the first $n$ elements in \eqref{eq:bp0}. In what follows, we write
$(Y_{t^\ast})_{t\ge 1}$ for the bootstrap \seq\ \eqref{eq:bp0}, indicating
that this \seq\ is nothing but the original $Y$-\seq\ sampled at the
random indices $(K_1,\ldots,K_1+L_1-1,K_2,\ldots,K_2+L_2-1,\ldots)$
with the convention that indices larger than $n$ are taken
modulo $n$.
\par
In what follows, the \pro y measure
generated by the bootstrap procedure is denoted by $P^\ast$, i.e.,
$P^\ast(\cdot)=P(\cdot \mid (X_t))$. The corresponding
expected value, variance and covariance are denoted by $E^\ast$,
$\var^\ast$ and $\cov^\ast$.
\subsection{The bootstrapped sample extremogram}
Davis et al.~\cite{davis:mikosch:cribben:2012} applied the stationary
bootstrap to the \seq\ of lagged vectors
\beao
I_t(h)= ( I_t^2,  I_t  I_{t+1},\ldots, I_t   I_{t+h})\,,\quad t=1,2,\ldots\,,
\eeao
for fixed $h\ge 0$ and showed consistency of the bootstrapped 
sample extremogram. 
In particular, they showed the following result which we cite for
further reference. A close inspection of the 
proof in \cite{davis:mikosch:cribben:2012} shows that the results
remain true if in $I_t(h)$ we replace the quantities $I_s$ by $\wt
I_s$, $s=t,\ldots,t+h$. We denote the corresponding vector by $\wt
I_t(h)$. Consider the stationary bootstrap \seq\
$(\wt I_{t^\ast}(h))$ and write
\beao
\wt \gamma_A^\ast(i)= \dfrac mn \sum_{t=1}^{n-i} \wt I_{t^\ast}
\wt I_{t^\ast+i}
,\quad i=0,\ldots,h\,. 
\eeao

\bth \label{thm:bp1} Consider an 
$\mathbb{R}^d$-valued strictly stationary
regularly varying sequence $(X_t)$ with index $\alpha>0$ and 
assume the following conditions:
\begin{enumerate}
\item[\rm 1.] The mixing condition
{\rm (M1)} and in addition 
$ \sum_{h=1}^{\infty}h\xi_h <\infty$.
\item[\rm 2.] The growth conditions $\theta=\theta_n\to 0$ and $n\theta^2/m\to \infty$.
\item[\rm 3.] The set $A$ is bounded away
  from zero, $\mu_1(\partial A)=0$ and $\mu_1(A)>0$.
\end{enumerate}
Then the following bootstrap consistency results hold for $h\ge 0$:
\beao
E^{*}\big(\wt \gamma_A^\ast (h)\big) \overset{P}{\to}
\gamma_A(h)\quad\mbox{and}\quad
 \var^{*} \big( (n/m)^{0.5} \wt \gamma_A^\ast (h) 
\big)\overset{P}{\to} \sigma_{hh}\,,
\eeao
where the covariance matrix $\Sigma_h=(\sigma_{ij})$ is given in
Lemma~\ref{thm:acf}.
Moreover, writing $d$ for any metric describing weak \con\ in
Euclidean space and $(Z_i)_{i=0,\ldots,h}$ for an $N(0,\Sigma_h)$
Gaussian vector, we also have
\beao
d\Big((n/m)^{1/2}\big(\wt
{\gamma}_A^\ast(i)-\widetilde{\gamma}_A(i)\big)_{i=0,\ldots,h}, (Z_i)_{i=0,\ldots,h}\Big)
\stp 0\,,\quad \nto \,.
\eeao
\ethe 
In what follows, we will write $d$ for any
metric describing weak \con\ in any space of interest. 

\subsection{The bootstrapped integrated \per }
Bootstrapping the \seq\ $(I_t(h))$ has the
advantage that we preserve  the neighbors $I_{t^\ast+i}$
of $I_{t^\ast}$ from the original \seq\ $( I_s)$. However, this method
depends on the lag $h$ and creates problems 
if the number of lags increases with the sample
size $n$. In what follows, we will apply the stationary bootstrap 
directly to $(I_{t})$. Then we have to re-define the 
bootstrap sample extremogram at any lag
$h< n$. Write 
\beao
\overline{I}_n= n^{-1} \sum_{t=1}^n I_t\quad\mbox{and}\quad  
\wh{I}_{t}= I_{t}-\overline{I}_n\,,\quad t\in \bbz\,,
\eeao 
and define the corresponding bootstrap sample extremogram
\beao
\wh\gamma_A^{*}(h)= \dfrac mn
\sum_{t=1}^{n-h}\wh\Imath_{t^\ast}
\wh\Imath_{(t+h)^\ast}\,,\quad h=0,\ldots,n-1\,,
\eeao 
and the bootstrap periodogram
\beao
I_{n,A}^\ast(\la)= \frac m n \Big|\sum_{t=1}^n\wh I_{t^\ast}
\ex^{-it\,\la}\Big|^2\,,\quad \la\in \Pi\,.
\eeao
Note the crucial difference: in general, 
$\Imath_{t^\ast} \Imath_{(t+h)^\ast} \ne \Imath_{t^\ast}
\Imath_{t^\ast +h} $, but, as we will see in Lemma~\ref{lem:bpclt1}, the quantities 
$\widetilde{\gamma}_A^{*}(h)$ and $\wh{\gamma}_A^{*}(h)$
are \asy ally close for fixed $h\ge 0$. 
\par
In what follows, we focus on the bootstrap for the continuous version
$J_{n,A}$ of the integrated  \per\ for a given smooth weight
\fct\ $g$; bootstrap consistency can also be shown for the  
discretized version $\wh J_{n,A}$; we omit further details.
In the definition of $J_{n,A}$ in \eqref{eq:cont}, we
simply replace $(I_t)$ by $(\wh I_{t^\ast})$, resulting in 
its bootstrap version
\beao
J^\ast_{n,A}(\la)  &=& \int_0^\la I_{n,A}^\ast(x) \,g(x)\,dx
=\psi_0 \,\wh \gamma_A^\ast(0)+2 
\sum_{h=1}^{n-1} \psi_{h}\,\wh \gamma_A^\ast(h)\,,\quad \la\in \Pi\,.
\eeao
\par
Now we can formulate a bootstrap analog of Theorem~\ref{thm:main} which shows
the consistency of the stationary bootstrap procedure.
\bth \label{thm:bpclt}
Assume the conditions of Theorem~\ref{thm:main} and \ref{thm:bp1}.
%\beao
%\label{eq:bpclt1}E^\ast\big( \widehat{\Jmath}_{nA}^{*}(\la)\big) &
%\overset{P}{\to}& J_{A}(\la) \,.
%\eeao 
Then 
\beao
d\Big((n/m)^{1/2}\big(\Jmath_{n,A}^\ast-  
E^\ast J_{n,A}^\ast\big), G\Big)\stp 0\,,\quad \nto\,,
\eeao
where the Gaussian process $G$ is defined in Theorem~\ref{thm:main}
and $d$ is any metric which describes weak \con\ in
$\bbc(\Pi)$.
\ethe
\bre\label{rem:bias}  Recall that, in general, it is not possible to replace the 
centering $EJ_{n,A}$ of $J_{n,A}$ in the 
\fct al \clt\ of Theorem~\ref{thm:main} by its limit 
$\int_0^\cdot h_A(\la)\,g(\la)\,d\la $. A similar remark applies 
to Theorem~\ref{thm:bpclt}.  Although 
$\sup_{\la\in \Pi} |E^\ast J_{n,A}^\ast(\la)-J_{n,A}(\la)|\stp 0$, under
the conditions of Theorem~\ref{thm:bpclt},
it is in general not possible to replace 
the centering $E^\ast J_{n,A}^\ast$ by $J_{n,A}$; see 
Lemma~\ref{lem:bias}. Thus, Theorem~\ref{thm:bpclt}
does not yield bootstrap consistency in a textbook sense but it rather
provides a simulation technique for the limit process $G$. In
turn, the simulation of this process makes it possible to provide
confidence bands for the goodness of fit test statistics considered
above. We will apply this simulation procedure in
Section~\ref{sec:simulation}.
\ere
\subsection{A simulation study}\label{sec:simulation}
\begin{figure}
\begin{minipage}[ht]{0.33\linewidth}  
\centering  
\includegraphics[width=\textwidth]{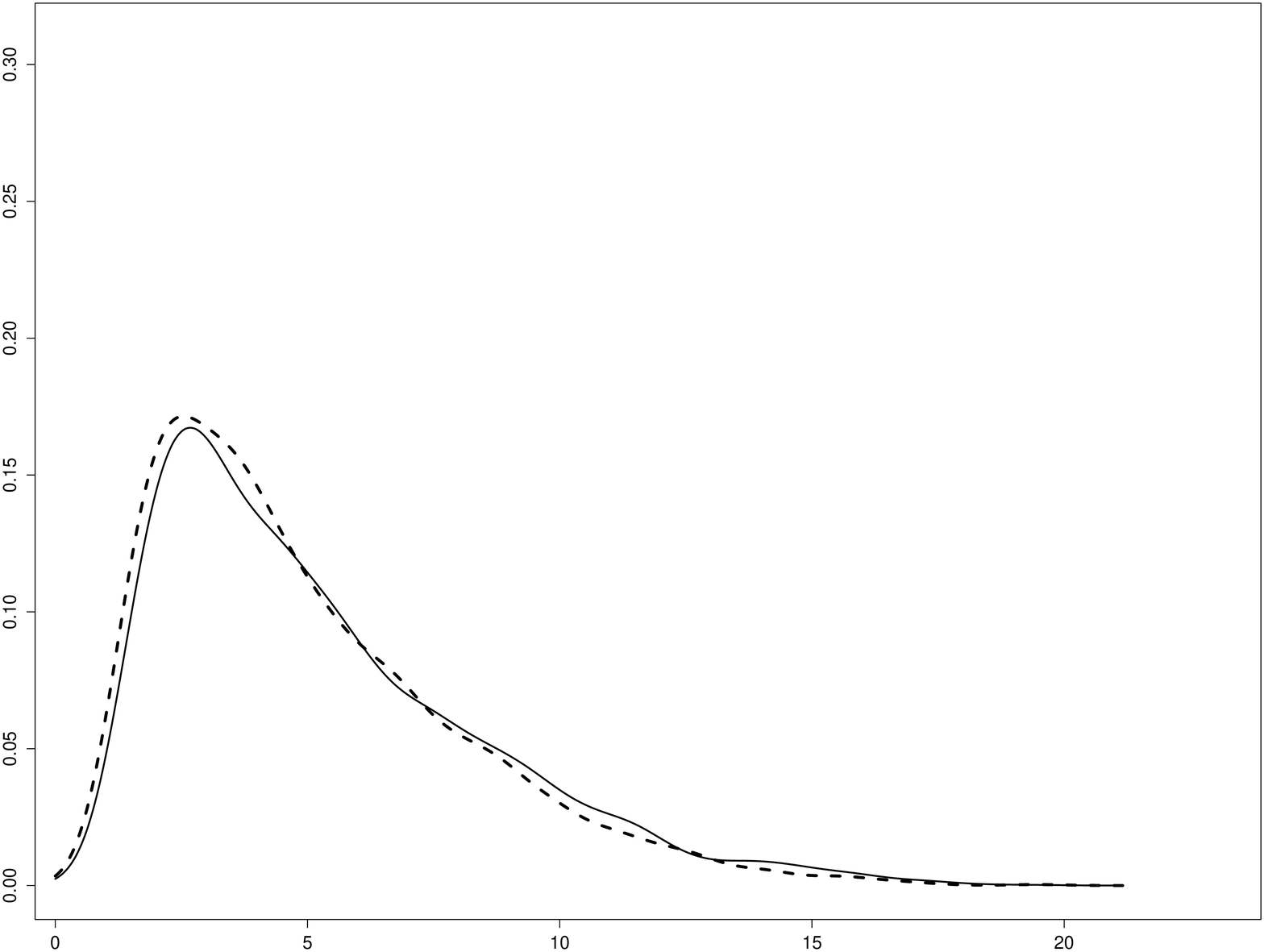}  
\end{minipage}%
\begin{minipage}[ht]{0.33\linewidth}  
\centering  
\includegraphics[width=\textwidth]{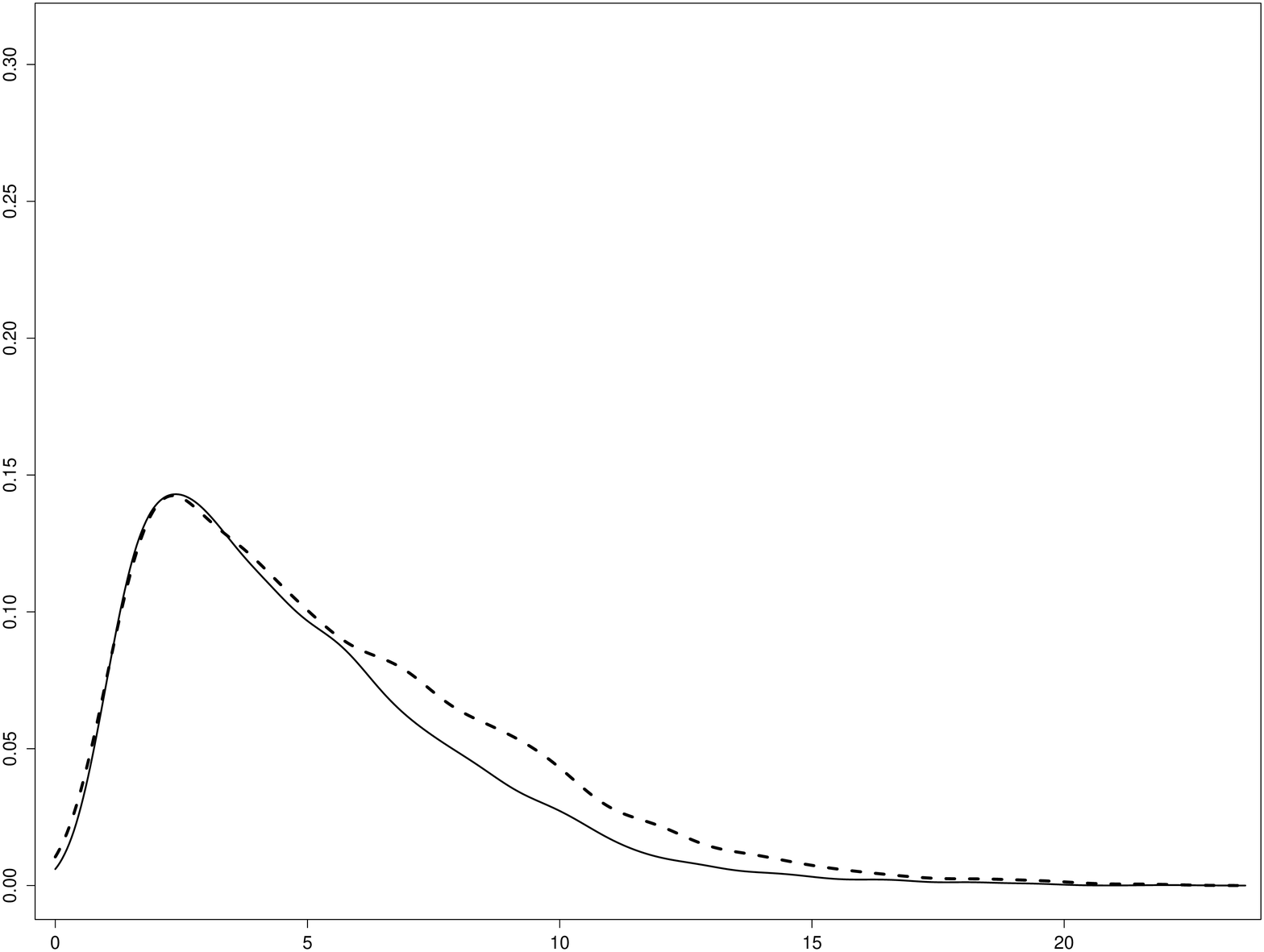}  
\end{minipage}%
\begin{minipage}[ht]{0.33\linewidth}  
\centering  
\includegraphics[width=\textwidth]{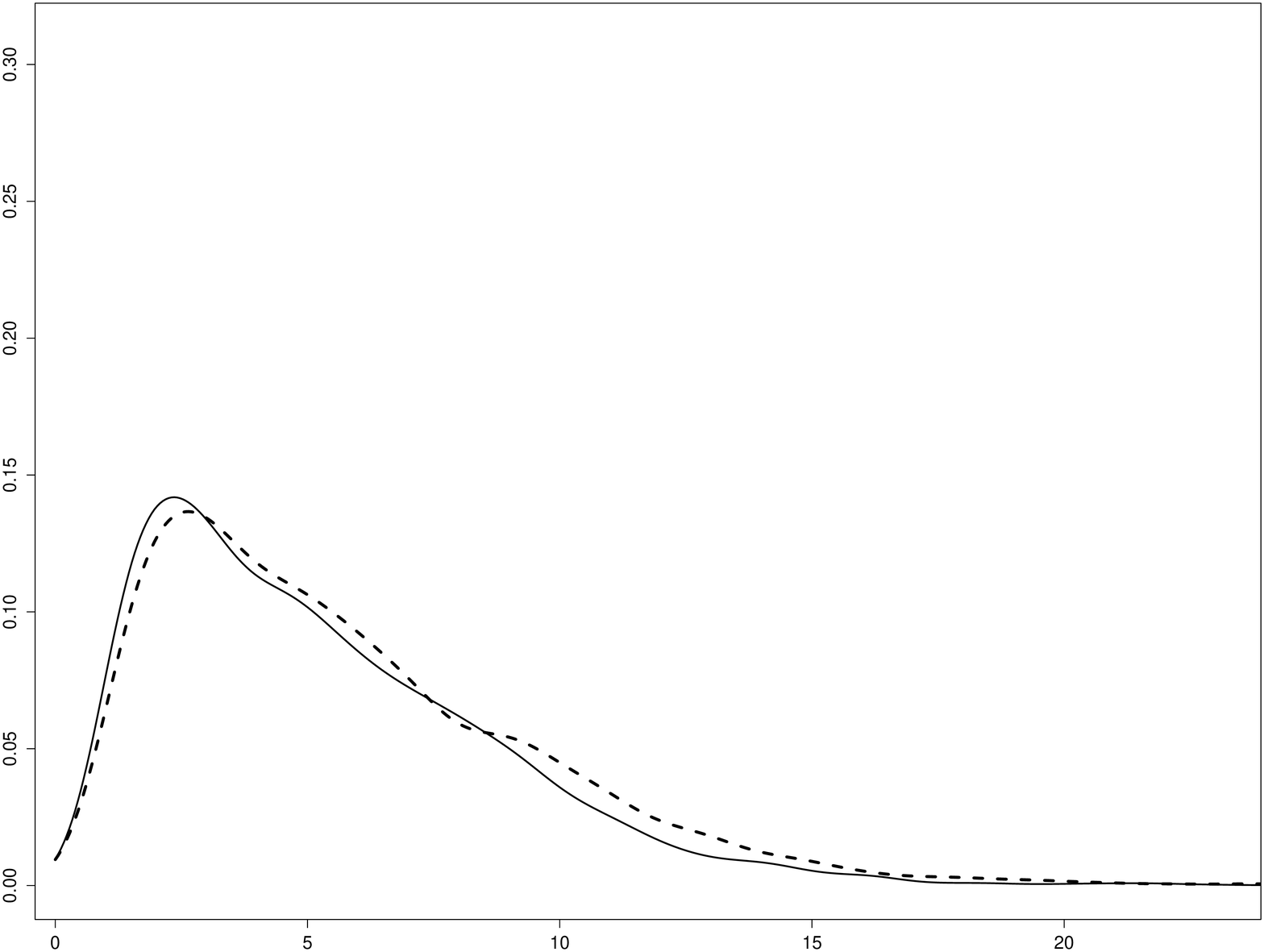}  
\end{minipage}%
\\
\begin{minipage}[htbp]{0.33\linewidth}  
\centering  
\includegraphics[width=\textwidth]{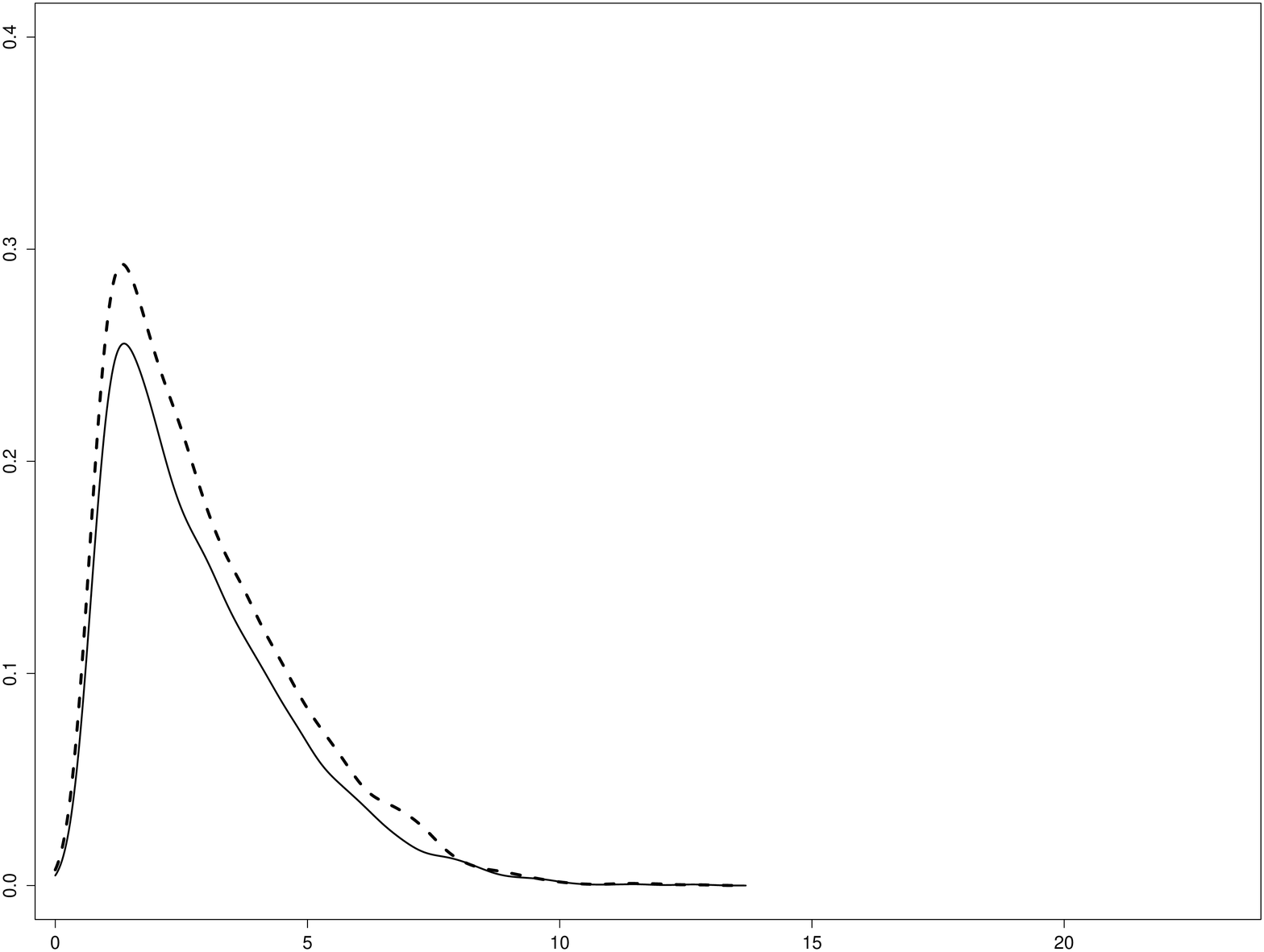}  
\end{minipage}%
\begin{minipage}[ht]{0.33\linewidth}  
\centering  
\includegraphics[width=\textwidth]{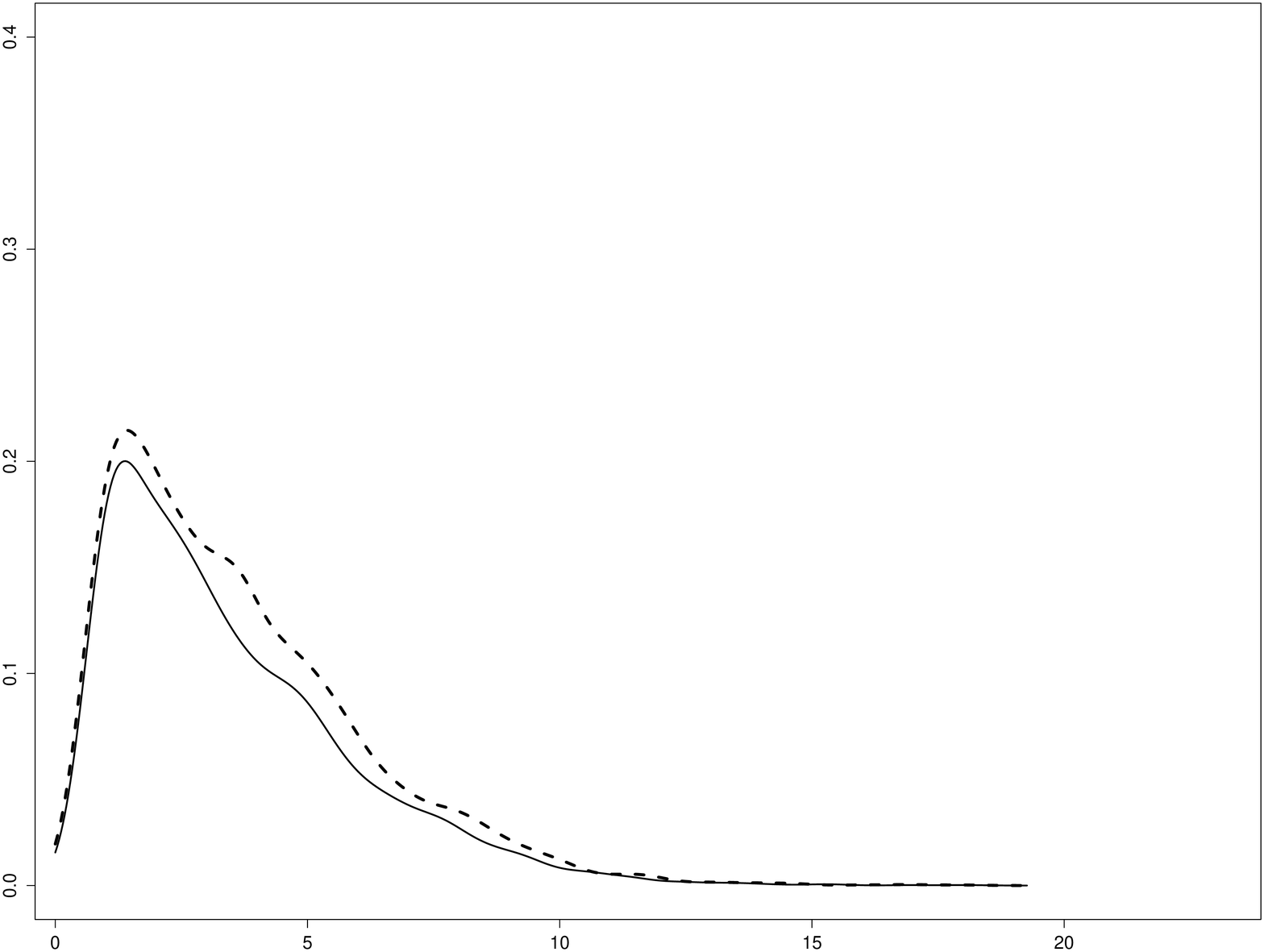}  
\end{minipage}%
\begin{minipage}[ht]{0.33\linewidth}  
\centering  
\includegraphics[width=\textwidth]{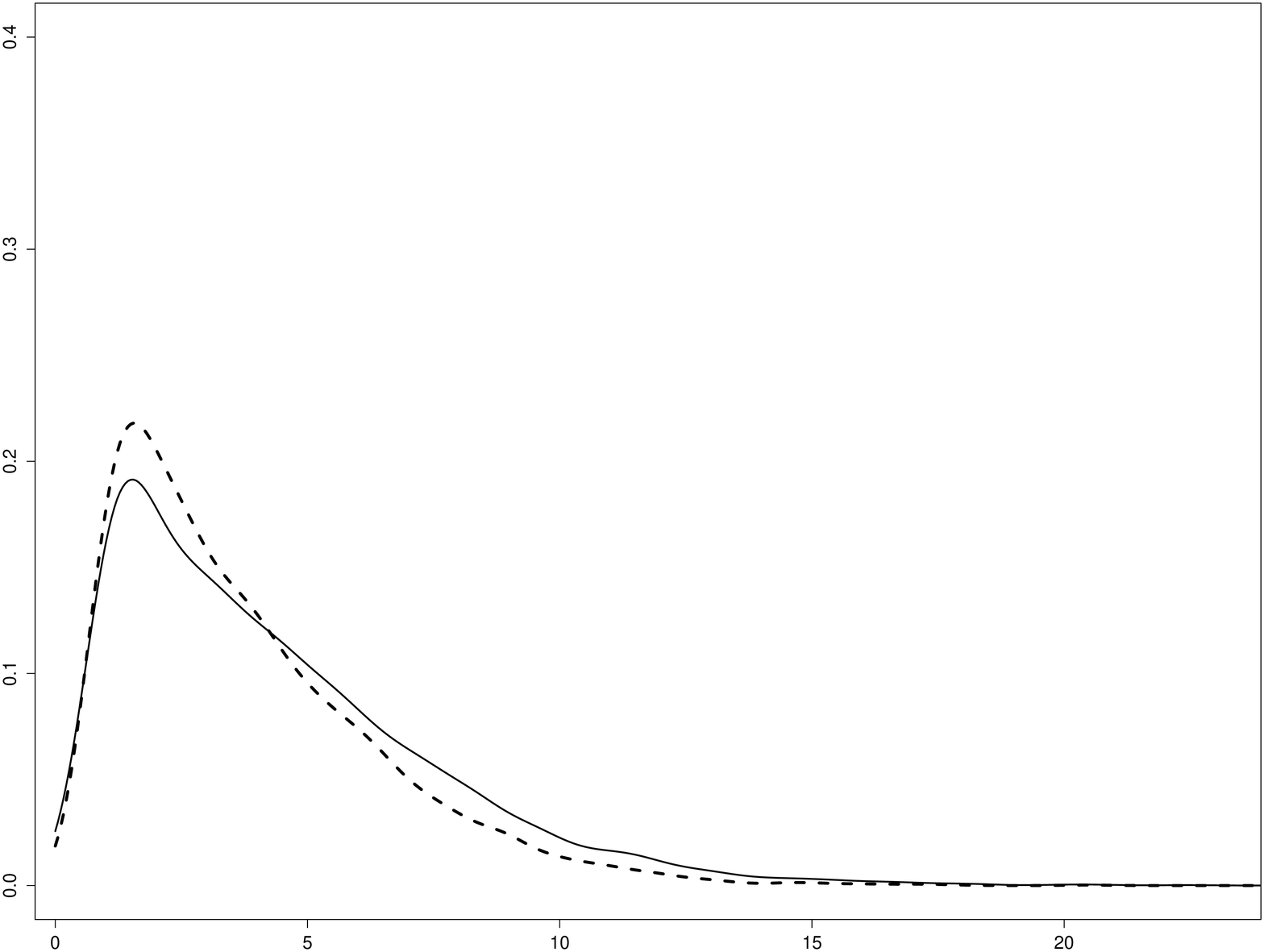}  
\end{minipage}%
\caption{\red Density of the normalized
GRS (solid
line) and its bootstrap approximation. 
The sample size is $n=2,000$.
and the thresholds $a_m$ are chosen such that
$p_0=P(X>a_m)=0.10, \, 0.05, \,0.03$ corresponding to the first,
second and third column. {\em Top:} The sample is drawn from 
the ARMA$(1,1)$ process $X_t=0.8X_{t-1}+0.1Z_{t-1}+Z_t$, where
$(Z_t)$ is iid $t$-distributed with $\alpha=3$ degrees of
freedom. 
{\em Bottom:} The sample is drawn from the
GARCH$(1,1)$ process $X_t=\sigma_tZ_t$, where
$\sigma_t^2=0.1+0.1X_{t-1}^2+0.84\sigma_{t-1}^2$ and $(Z_t)$ is iid
$t$-distributed with $4$ degrees of freedom. In this case, the index of \regvar\ for $(X_t)$ is 
$\alpha=3.49$; see Table~2 in 
Davis and Mikosch \cite{davis:mikosch:2009c}.}
\label{fig:03}
\end{figure}

\begin{figure}
\begin{minipage}[ht]{0.5\linewidth}  
\centering  
\includegraphics[width=\textwidth]{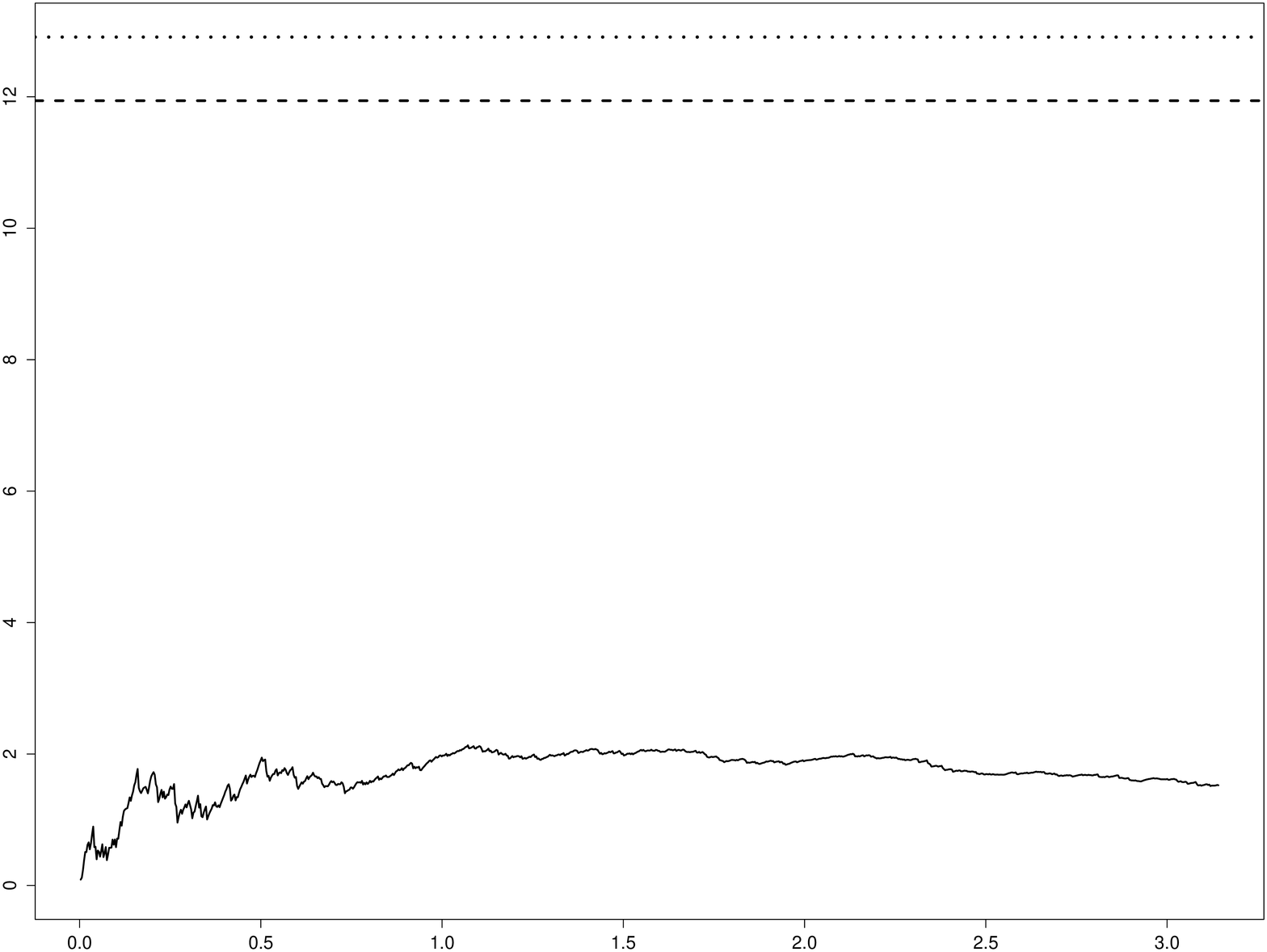}  
\end{minipage}%
\begin{minipage}[ht]{0.5\linewidth}  
\centering  
\includegraphics[width=\textwidth]{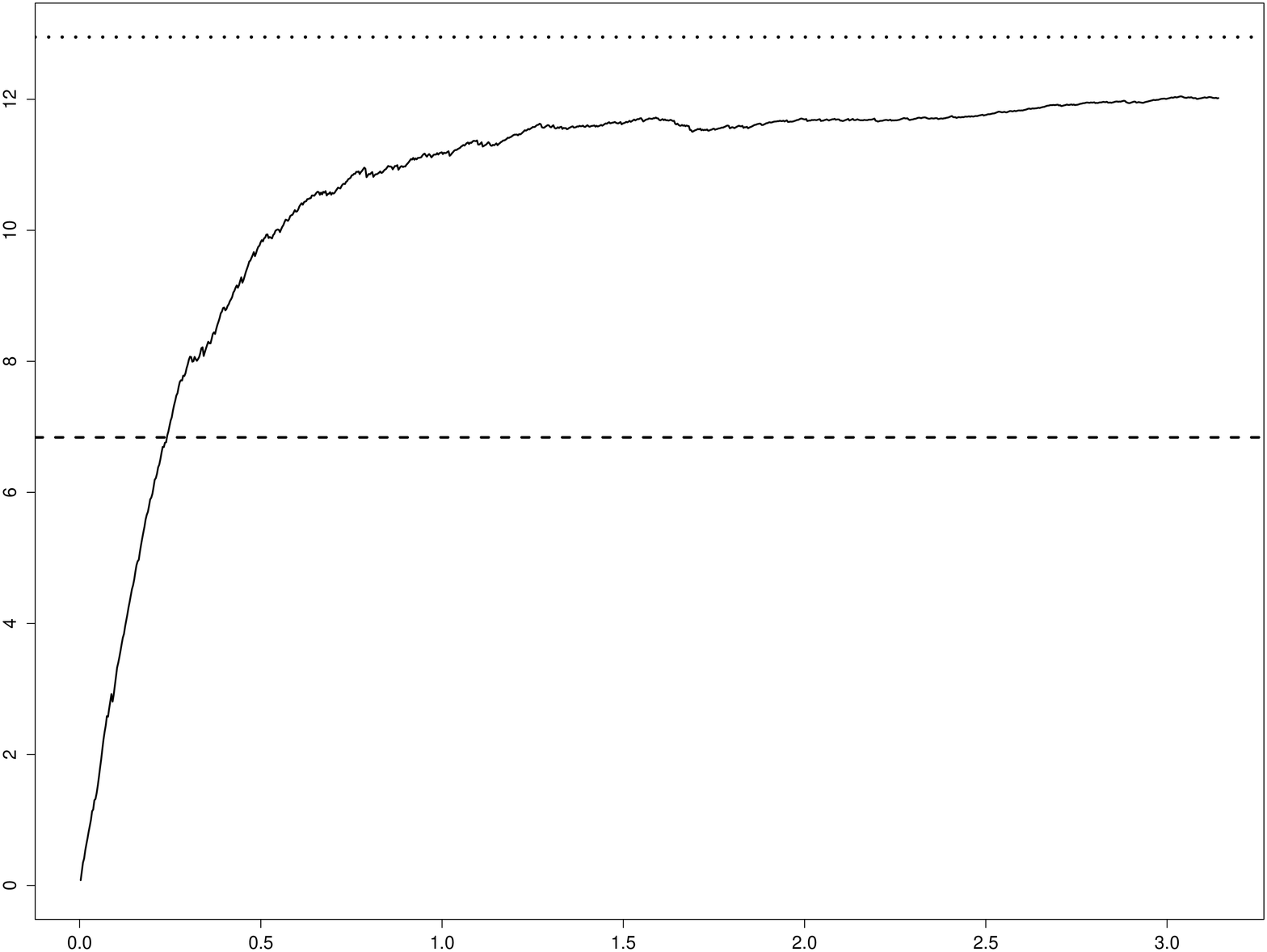}  
\end{minipage}%
\\
\begin{minipage}[ht]{0.5\linewidth}  
\centering  
\includegraphics[width=\textwidth]{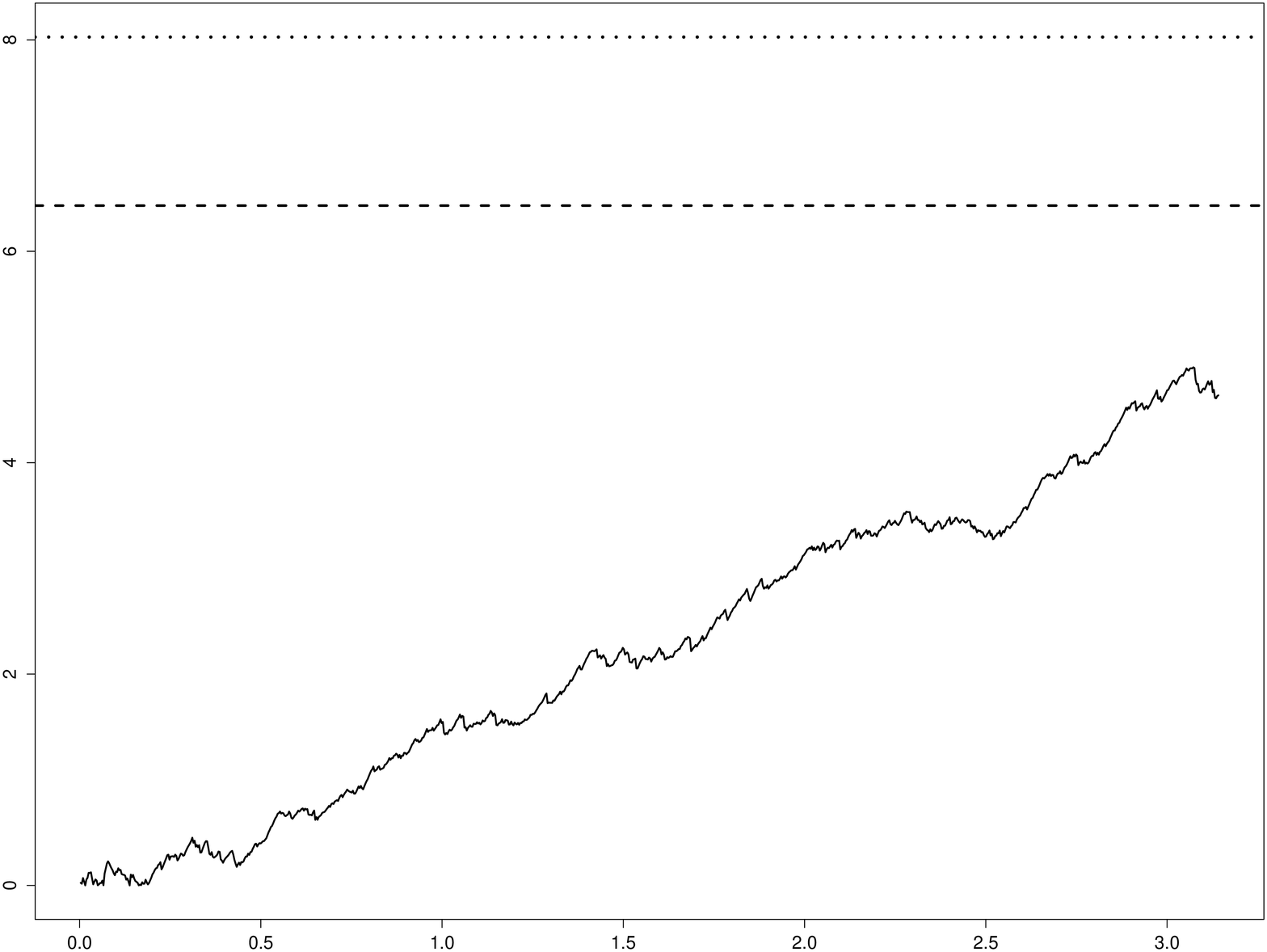}  
\end{minipage}%
\begin{minipage}[ht]{0.5\linewidth}  
\centering  
\includegraphics[width=\textwidth]{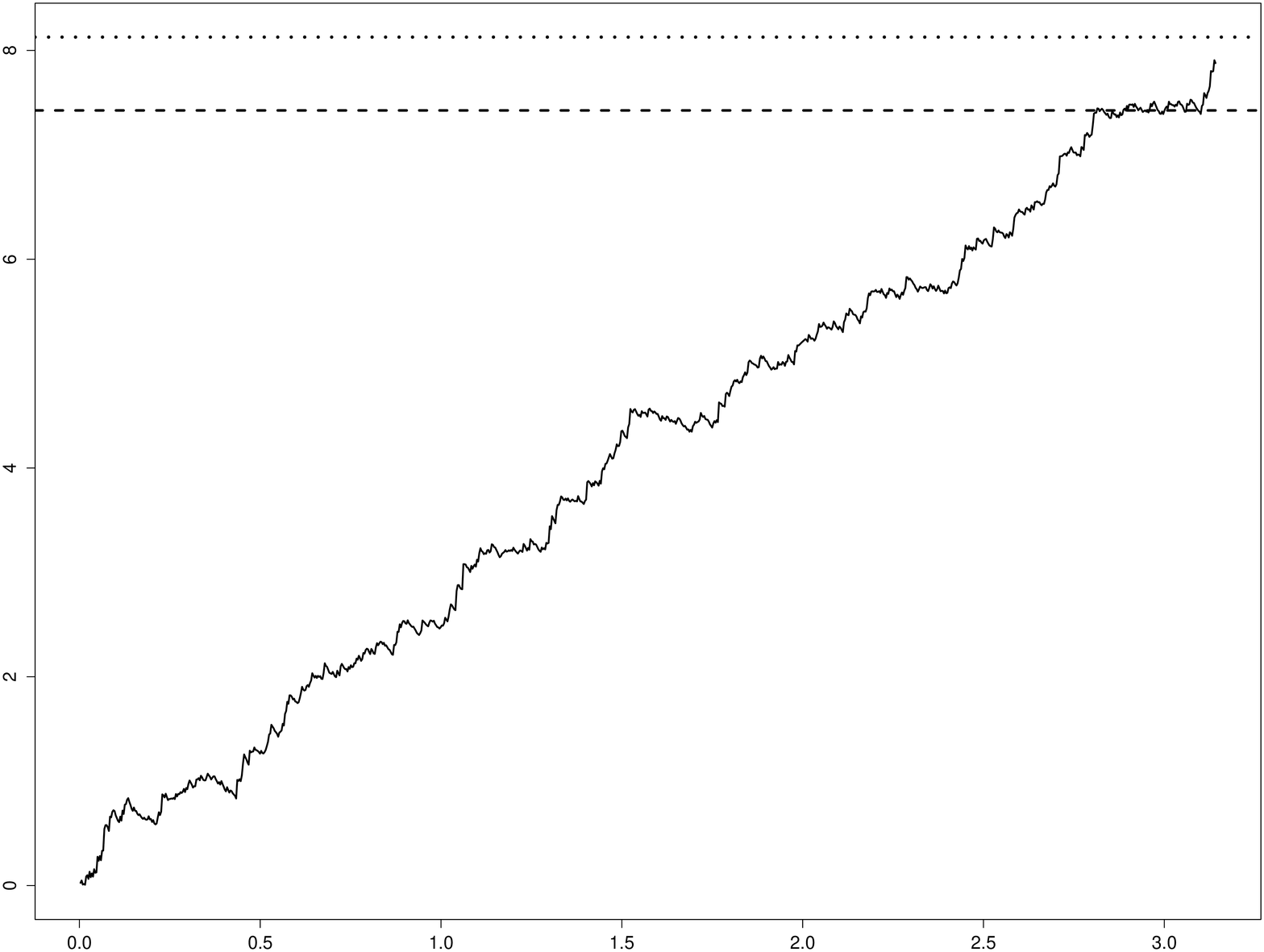}  
\end{minipage}%
\caption{\red Paths of the
integrated periodogram $(n/m)^{0.5}|J_{n,A}-EJ_{n,A}|$ with $p_0=0.05$ 
for samples of size 
$n=2,000$. {\bf Top:}  We work under the null hypothesis
of the ARMA(1,1) model
$X_t=0.8X_{t-1}+0.3Z_{t-1}+Z_t$, where $(Z_t)$ is iid $t$-distributed 
with $\alpha=3$ degrees of freedom. 
{\em Left:} The sample is drawn from the null model. 
The lower and upper dotted lines $y=11.9$ and $y=12.9$ correspond to the 
bootstrap-based and true 
95\%-quantiles of the GRS, respectively. The null hypothesis 
would be accepted.
{\em Right:} The sample is drawn from the 
ARMA(1,1) process 
$X_t=0.8X_{t-1}+0.1Z_{t-1}+Z_t$ with the same \ds\ for $(Z_t)$.
The lower dotted line $y=6.84$ is the  
bootstrap-based 95\%-quantile of the GRS. Based on it, the test would
reject the null. However, it would accept the null if one chose the
95\%-quantile of the null model.
{\bf Bottom:} We work under the null hypothesis of the 
GARCH$(1,1)$ process $X_t=\sigma_tZ_t$, where
$\sigma_t^2=10^{-7}+0.1X_{t-1}^2+0.81\sigma_{t-1}^2$ and $(Z_t)$ is iid
$t$-distributed with $4$ degrees of freedom. {\em Left:} The sample is
chosen from the null model. 
The lower and upper dotted lines $y=6.4$ and $y=8$ correspond to the 
bootstrap-based and true 95\%-quantiles of the GRS, respectively.
The null would be accepted for both quantiles. 
{\em Right:} The sample
is drawn from a GARCH$(1,1)$ process with
$\sigma_t^2=10^{-7}+0.1X_{t-1}^2+0.84\sigma_{t-1}^2$ and the same \ds\
of $(Z_t)$. The lower dotted line $y=7.4$
is the bootstrap-based 95\%-quantile of the GRS. The null would be
rejected in this case while it would be accepted if one used the 
95\%-quantile $y=8$ based on the null hypothesis. }
\label{fig:04}
\end{figure}

\begin{figure}
\begin{minipage}[ht]{0.5\linewidth}  
\centering  
\includegraphics[width=\textwidth]{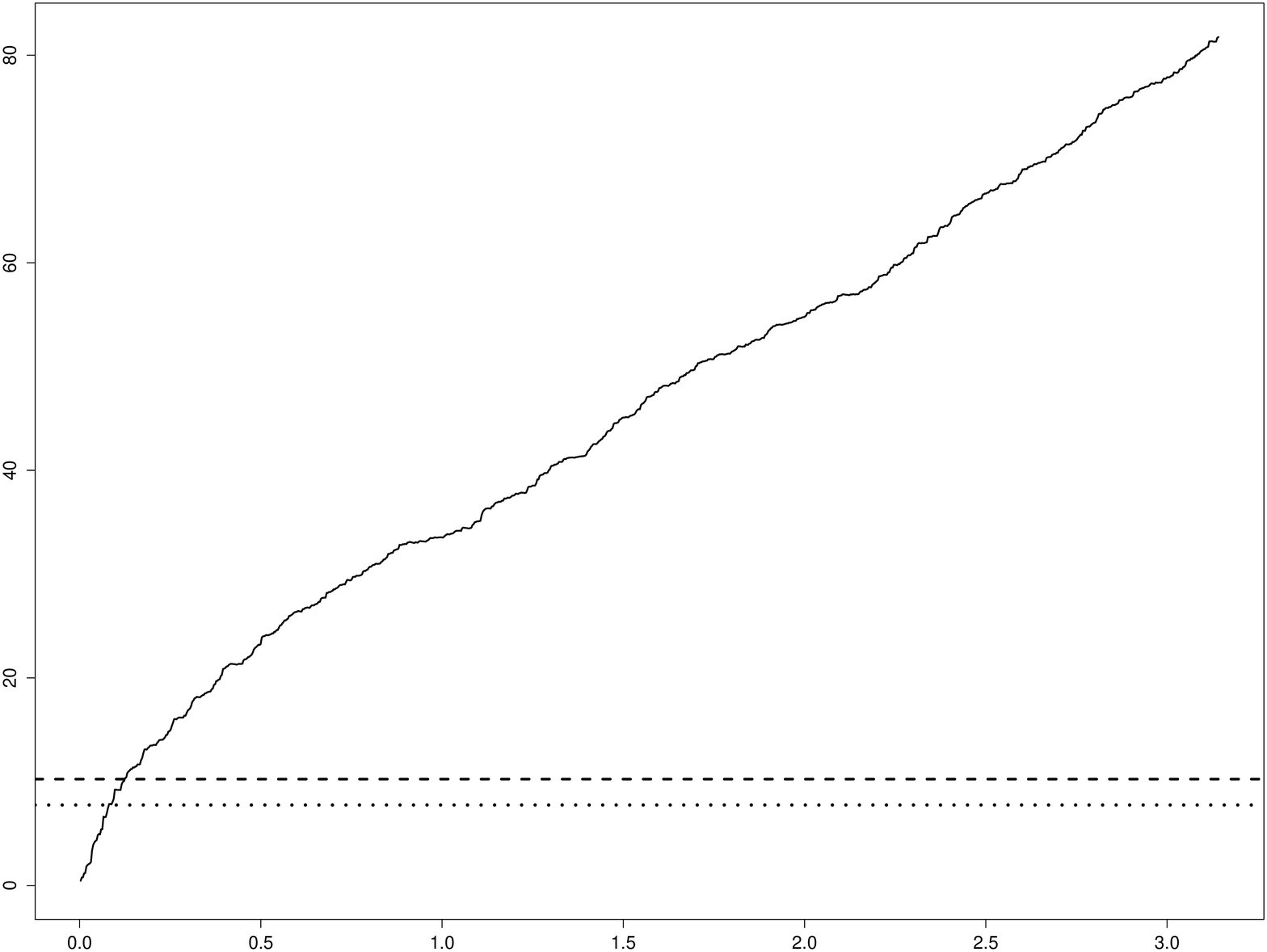}  
\end{minipage}%
\begin{minipage}[ht]{0.5\linewidth}  
\centering  
\includegraphics[width=\textwidth]{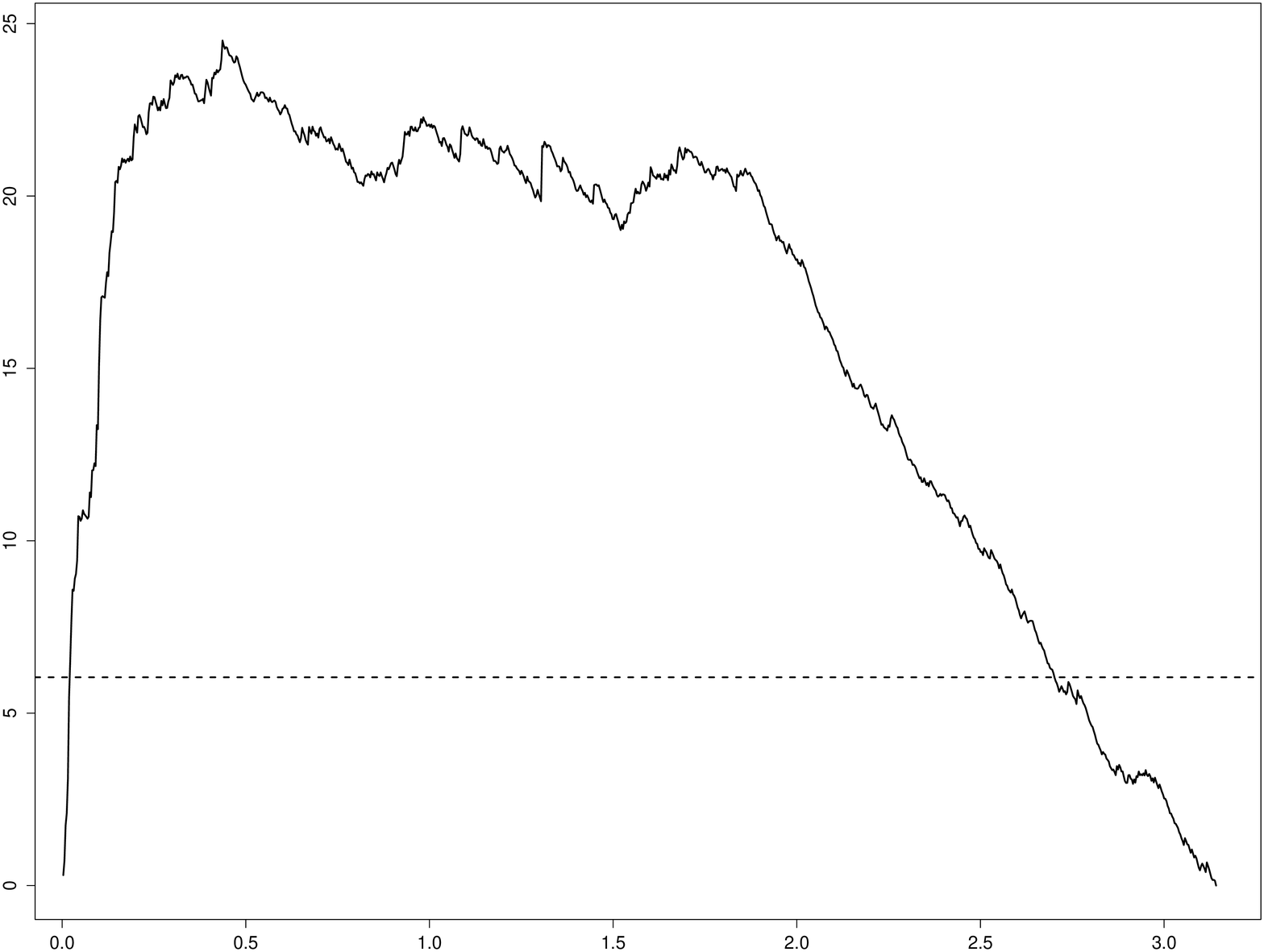}  
\end{minipage}%
\caption{\red The sample of size $n=2,000$ is drawn from a
stochastic volatility process $X_t=\sigma_tZ_t$ with
log-volatility $\log \sigma_t=0.9\log \sigma_{t-1}+\epsilon_t$
for an iid standard normal \seq\ 
$(\epsilon_t)$, $Z_t$ 
is $t$-distributed with $3.6$ degrees of freedom. 
\emph{Left:} Sample path of $(n/m)^{0.5}|J_{n,A}-EJ_{n,A}|$ with
$p_0=0.05$. 
The lower and upper dotted lines $y=7.8$ and $y=10.2$ correspond to
the true and bootstrap-based 95\%-quantiles of the 
GRS under the null hypothesis of  a
GARCH(1,1) process $\wt X_t=\wt \sigma_t\widetilde{Z}_t$ with
$\wt \sigma_t^2=6.23\times 10^{-3}+0.1\wt X_{t-1}^2+0.8\wt \sigma_{t-1}^2$ and iid
$t$-distributed $(\widetilde{Z}_t)$ with $4$ degrees of freedom. This
process has tail index $3.68$; 
see Table~1 in \cite{davis:mikosch:2009c}. The test clearly rejects the
null hypothesis.
\emph{Right:} Sample path of the integrated \per\ absolute
  value $n^{0.5}|J_{n,A}-\psi_0 \wt
\gamma_A(0)|$. The dotted line is the 95\%-quantile of the \ds\  of the
supremum of the absolute values of a Brownian bridge. The test clearly
rejects the null hypothesis that $(X_t)$ is iid.}
\label{fig:05}
\end{figure}

\begin{figure}
\begin{minipage}[ht]{0.5\linewidth}  
\centering  
\includegraphics[width=\textwidth]{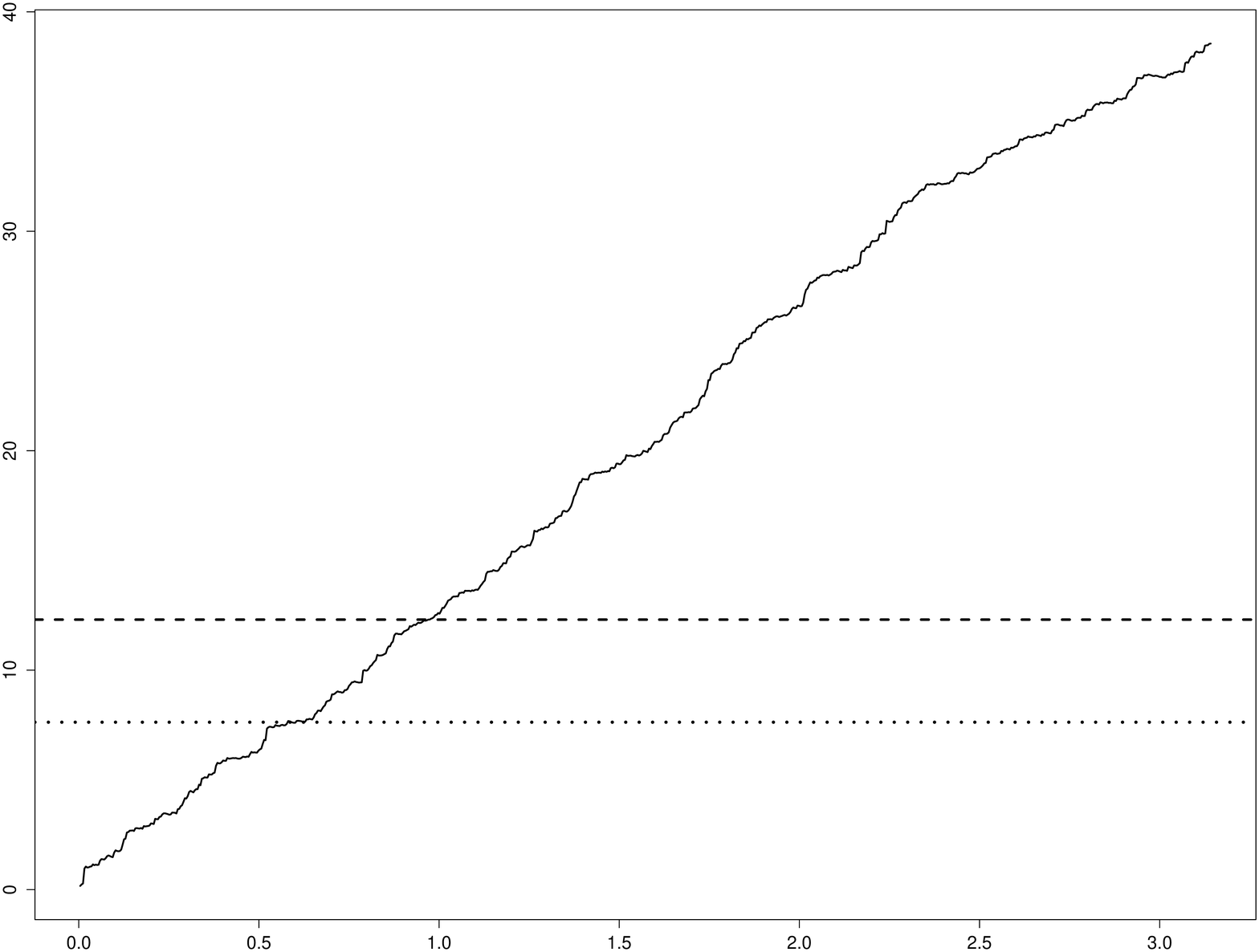}  
\end{minipage}%
\begin{minipage}[ht]{0.5\linewidth}  
\centering  
\includegraphics[width=\textwidth]{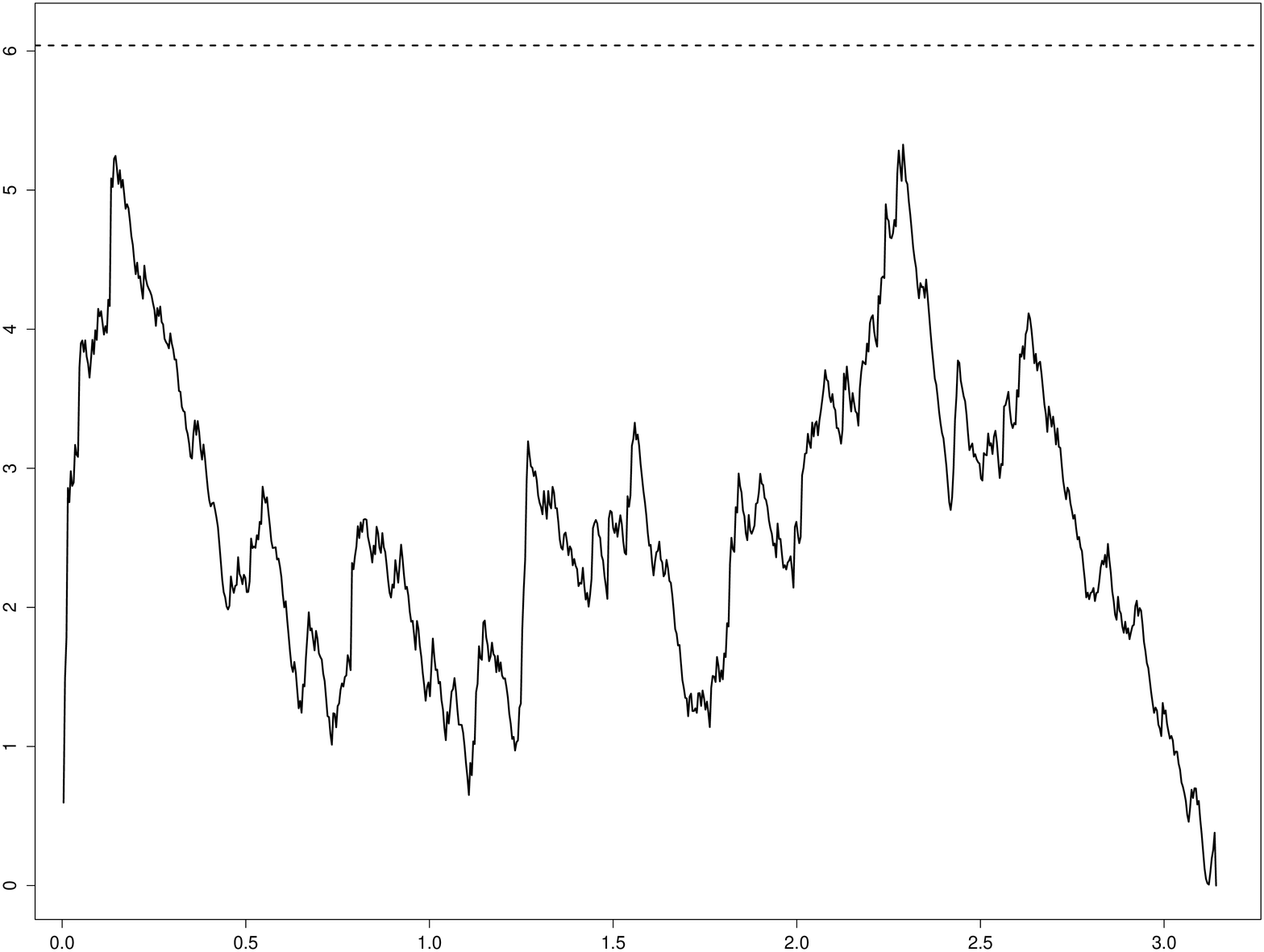}  
\end{minipage}%
\caption{\red GRTs for $1,560$ Goldman Sachs 1-minute log-returns.
{\em Left}: The integrated \per\ 
  $(n/m)^{0.5}|J_{n,A}-EJ_{n,A}|$ with $p_0=0.05$
 under the null hypothesis that the data are generated by the \garch\ model
 $\sigma_t^2=0.019+0.1X_{t-1}^2+0.87\sigma^2_{t-1}$ with iid $t$-distributed
   noise with 4 degrees of freedom.
The lower and upper dotted lines $y=7.6$ and $y=12.3$ represent 
the true and bootstrap-based 95\%-quantiles of the GRS 
under the null hypothesis. The hypothesis
of \garch\ is clearly rejected.
{\em Right:}  The integrated \per\ $n^{0.5}|J_{n,A}-\psi_0 \wt
\gamma_A(0)|$ with $p_0=0.05$ under the null hypothesis of an iid \seq . 
The dotted line represents the \asy\ 95\%-quantile based on the
approximation of the GRS by the supremum of the absolute values 
of a Brownian bridge. The null hypothesis is not rejected.}\label{fig:07}
\end{figure}

{\red 
We focus on the Grenander-Rosenblatt statistic (GRS) on the \lhs\ of 
\eqref{eq:g-r} for different \ts\ models, distinct thresholds and 
sample sizes. Under the null hypothesis of a particular \ts\ model,
one can simulate the quantiles of the GRS from the theoretical model.
In this study we also follow a different approach. First, we determine
the expected value \fct\ $EJ_{nA}$ and the threshold 
$a_m$ \st\ $p_0=P(X>a_m)=1/m$ by simulation
from the theoretical model and then we use the stationary bootstrap to
calculate the \asy\ quantiles of the GRS. 
This \ds\ is be obtained by repeated simulation of 
$(n/m)^{0.5}\sup_{x\in \Pi}| J_{n,A}^\ast(x) - E^\ast
J_{n,A}^\ast(x)|$; Theorem~\ref{thm:bpclt} provides a
justification for this approach.\footnote{
  Throughout, to exploit the power of the Fast Fourier Transform, we
  use the Riemann sum approximations to the integrated \per s. We do
  not indicate this fact in the notation.} 
In the cases when the
expected value \fct\  $EJ_{nA}$ can be replaced by its limit, 
i.e., when the bias of $J_{n,A}$ is negligible, this 
approach has the advantage that the test is 
non-parametric. An example are models satisfying the \asy\
independence condition (AI) and $n/m^2\to 0$ as $\nto$; see 
Theorem~\ref{thm:mdep2a} and the remark following it.
Of course, for an iid \seq\ 
or $\eta$-dependent \seq\ one can also use the quantiles of the 
limit \ds\ of the GRS which are known or can be simulated; see
\eqref{eq:pp1}
and \eqref{eq:ppr1}.
\par  
In what follows, we apply the Grenander-Rosenblatt test (GRT) to various
univariate (real-life or simulated) 
time series  $X_t,t=1,\ldots,n$ for different sample sizes $n$ and
thresholds $a_m$. 
We always choose $A=(1,\infty)$ and $g\equiv 1$.  Whenever we apply
the stationary bootstrap we choose the geometric parameter
$\theta=1/50$. Density plots and simulated
quantiles are derived  from $4,000$ independent repetitions,
also in the bootstrap case.
\par
In Figure~\ref{fig:03} we illustrate how the stationary bootstrap
works for different thresholds $a_m$ and sample size $n=2,000$.
We show the density of the 
normalized GRS on the \lhs\
of \eqref{eq:g-r} and its bootstrap approximation. We choose
\regvary\ ARMA$(1,1)$ and GARCH$(1,1)$ models. The densities of the
GRS and its bootstrap approximation are close
to each other. We take this fact as justification for using the
bootstrap quantiles of the GRS in the test. While the densities in
the ARMA case do not seem too sensitive to the choice of the high
threshold $a_m$, the shape of the densities change for the GARCH
model when switching from $p_0=0.10$ to $p_0=0.05$, while they look similar
for $p_0=0.05$ and $p_0=0.01$. 
%We do not address the accuracy of the
%approximation in the tails in this paper.
\par 
In Figure~\ref{fig:04} we show sample paths of the normalized and 
centered integrated periodogram $(n/m)^{0.5}|J_{n,A}-EJ_{n,A}|$
  with $p_0=0.05$
for samples of size $n=2, 000$ from ARMA(1,1) and 
GARCH(1,1) models together with
95\%-quantiles of the GRS both under the correct and under an incorrect
null hypothesis. Due to the need of centering with $EJ_{n,A}$ these
sample paths are affected both by the sample and the model. Indeed, if the
model is chosen incorrectly we will typically subtract 
the incorrect centering and calculate an incorrect threshold
$a_m$. When using both the bootstrap-based or true
95\%-quantiles of the GRS, the model is not
rejected if the sample is in agreement with the null hypothesis. However,
if the sample comes from a model whose parameters slightly deviate
from the parameters of the null hypothesis the  
incorrect expected value $EJ_{nA}$ and wrong threshold $a_m$ 
change the sample
path of the integrated \per\ in such a way that the bootstrap-based
GRT  rejects the null hypothesis while it
does not reject the null if one uses the quantiles based 
on the null hypothesis. It is advantageous to show both
95\%-quantiles: they deviate rather significantly, indicating that
we chose an incorrect null model.
\par
In Figure~\ref{fig:05} we consider a \sv\ model $X_t=\sigma_tZ_t$,
where $(\sigma_t)$ is a log-normal stationary process
independent of the iid $t$-distributed \seq\ $(Z_t)$. The
$\alpha$ degrees of freedom of the $t$-\ds\ coincide with the index
$\alpha$ of \regvar\ of $(X_t)$. The extremogram of this \seq\
vanishes at all positive lags. This fact is in agreement with the extremogram
of an iid \seq\ but is in contrast to a \garch\ process. 
Choosing $p_0=0.05$, we apply the
GRS to a \sv\ sample of size $n=2,000$ under the incorrect 
null hypothesis of a \garch\ model with tail index close to the chosen
$\alpha$. The test clearly rejects the null hypothesis. We also run a
GRT for the \sv\ sample under the null hypothesis of an iid 
$t$-distributed \seq\ with $\alpha$ degrees of freedom. We use the
approximation of the \ds\ of the GRS by the \ds\ of the supremum of a
Brownian bridge;   see Example \ref{exam:bb}. Also in this case, the null
is clearly rejected.
\par
In Figure~\ref{fig:07} we deal with a \ts\  $(X_t)$  of $1,560$
1-minute log-returns of 
Goldman Sachs stock from the period November 7-10, 2011. It has
estimated tail index $\alpha\approx 3$. Using standard software, we
fitted a \garch\ model 
\st\ $\sigma_t^2=0.019+0.1X_{t-1}^2+0.87\sigma^2_{t-1}$.
Hill and QQ plots of
the residuals of this model indicate that the noise is well 
fitted  by a $t$-\ds\ with (approximately) 4 degrees of freedom. 
The theoretical index of regular variation of this GARCH(1,1) model is
$\alpha=3.13$; 
see Table~2 in \cite{davis:mikosch:2009c}. We test the null hypothesis
of a \garch\ model with the aforementioned parameters. This hypothesis
is rejected. On the other hand, the GRT passes under the hypothesis of an
iid \seq , where we choose $a_m$ as the 95\% empirical quantile.
This means that the extremes of this data set 
are more in agreement with an
iid than with a GARCH structure. This is perhaps not
surprising in view of a  high \fre y data return series while GARCH
seems more suitable for  fitting low frequency returns. 
\par
A \garch\
model is often considered to give a  good fit to daily log-returns of stock
prices and foreign exchange (FX) rates. For example, such a judgement
may be  based on tests for zero autocorrelation of the residuals, their
absolute values and squares. We did not find evidence of GARCH
behavior in the extremes of three $5$-year \ts\ of daily Euro-USD FX rate
log-returns: from 2002 to 2006 (before the financial crisis), from
2006 to 2010 (including the financial crisis), from 2009 to 2013
(after the financial crisis); see
Figure~\ref{fig:08}. We choose different thresholds $a_m$.
When $p_0=0.05$ the null hypothesis of an iid \seq\ is
accepted for 2002-2006 and 2009-2013, but not for 2006-2010. 
The null hypothesis 
of a fitted GARCH process with $\sigma_t^2=2.37\times
10^{-7}+0.1X_{t-1}^2+0.8\sigma^2_{t-1}$ and iid $t$-distributed
noise with 4 degrees of freedom is also rejected by the GRT for 
2006-2010. 
For this latter period, the stationarity assumption may be doubted. 
We repeat the GRTs for $p_0=0.02$ in the periods 2002-2006 and
2009-2013. In the latter case the iid null hypothesis is still not
rejected while it is rejected in the former case. The abrupt change of
the behavior of the GRT may be due to the sample size (roughly 1,280
for each \ts ). For $p_0=0.02$ one would use only 2\% of the data for
the calculation of the GRT.
%{\blue We also mention the surprising fact that the residuals of 
%GARCH models fitted to the three \ts\ pass the GRT under the iid
%hypothesis. (The residuals for the second time series do not pass the GRT. ) A heuristic explanation of this phenomenon may be found in 
%\sta\ \cite{starica:2003}.}
\par
Our simulation study points at some of the problems one has to face
when using goodness of fit tests based on the extremes of a \ts .
A major problem is the choice of the threshold $a_m$. A data driven
choice would be preferable but we do not have a theoretical answer to the
problem. We propose to use graphical methods to compare the shapes 
of the extremogram and the integrated \per\ for different thresholds
and to choose a sufficiently high threshold where 
the shapes stabilize. A message from the simulations is that the
sample size $n$ should not be too small.  For example, the GRTs in
Figure~\ref{fig:08}  with  $n\approx 1,280$ give rather distinct answers when
switching from $p_0=0.05$ to $p_0=0.02$. The sample extremogram and the
integrated \per\ render meaningless for too high thresholds because
most indicator \fct s of extreme events will be zero. The
simulation study indicates that it is useful to exploit the
true quantiles of the GRS 
(obtained by simulation from a model under the null hypothesis) as
well as corresponding bootstrap-based quantile of the GRS. In particular,
when the null hypothesis is incorrect the two 95\% quantiles (say) will
typically differ, pointing at the incorrect null hypothesis.
We do not address the problem of goodness of fit tests in the case when
the null hypothesis depends on estimated parameters; the \asy\ theory
does not change under mild conditions on the \con\ rates of the estimators.

}

\begin{figure}
\begin{minipage}[ht]{0.4\linewidth}  
\centering  
\includegraphics[width=\textwidth]{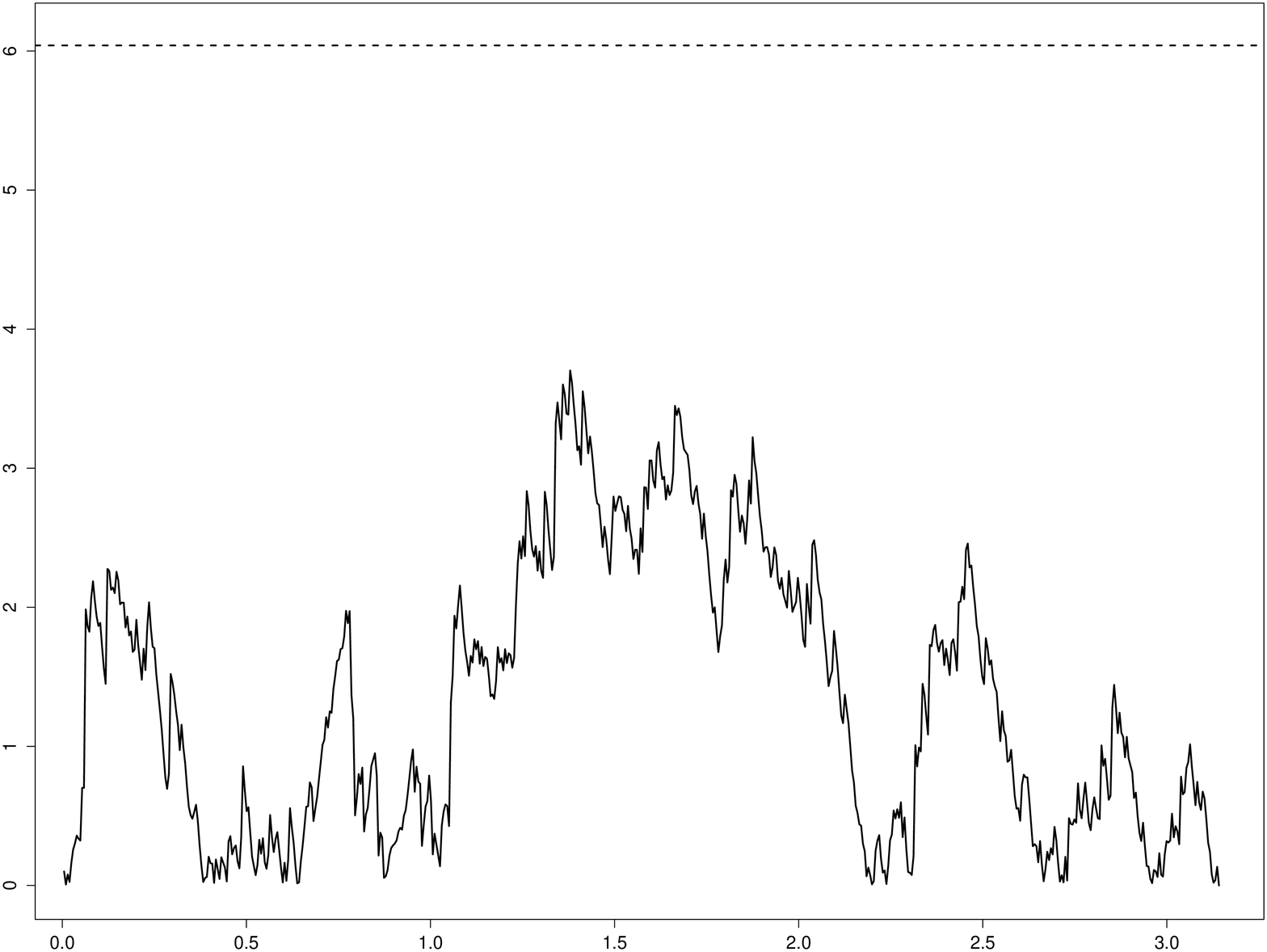}  
\end{minipage}%
\begin{minipage}[ht]{0.4\linewidth}  
\centering  
\includegraphics[width=\textwidth]{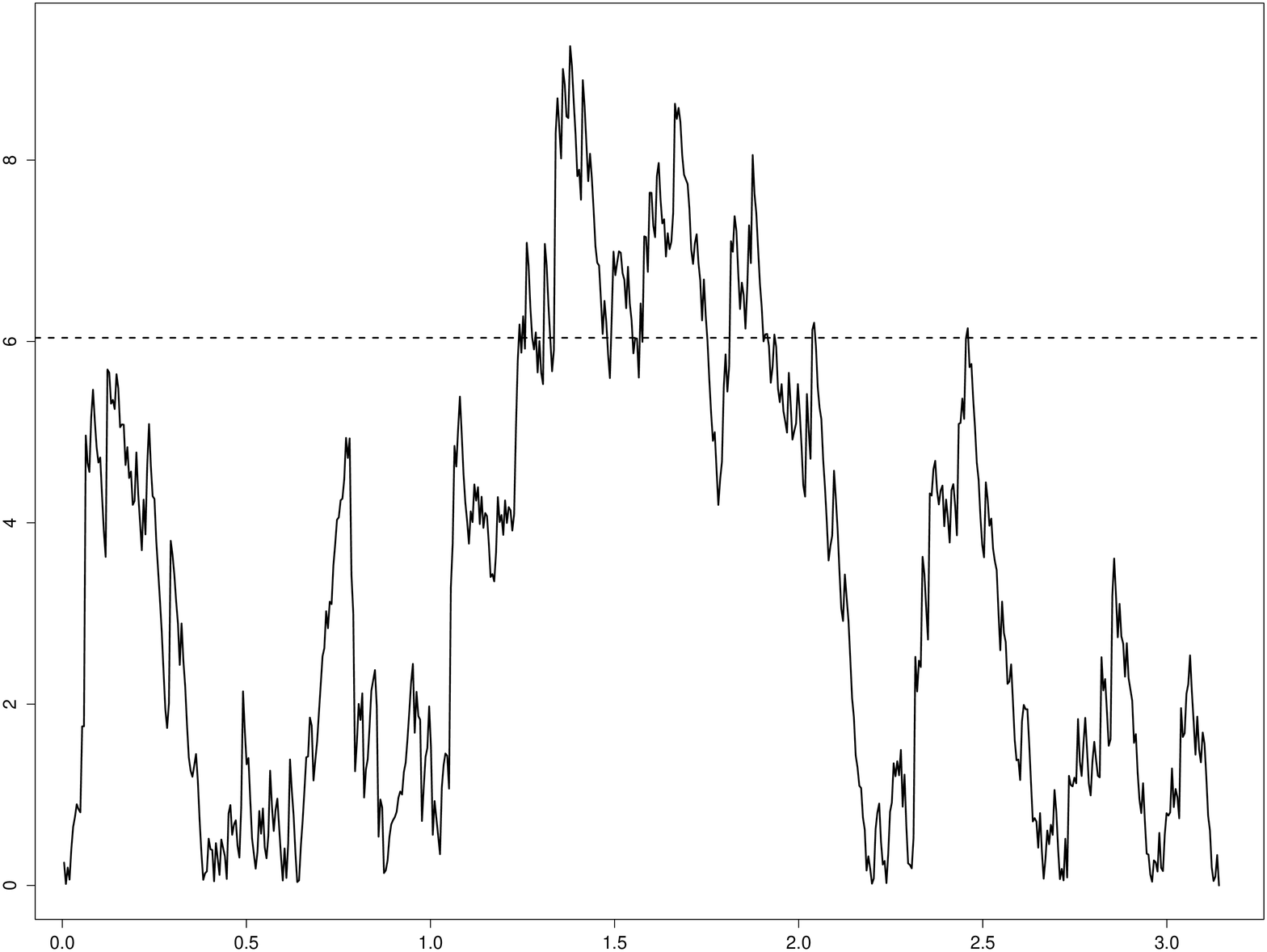}  
\end{minipage}%
\\
\begin{minipage}[ht]{0.4\linewidth}  
\centering  
\includegraphics[width=\textwidth]{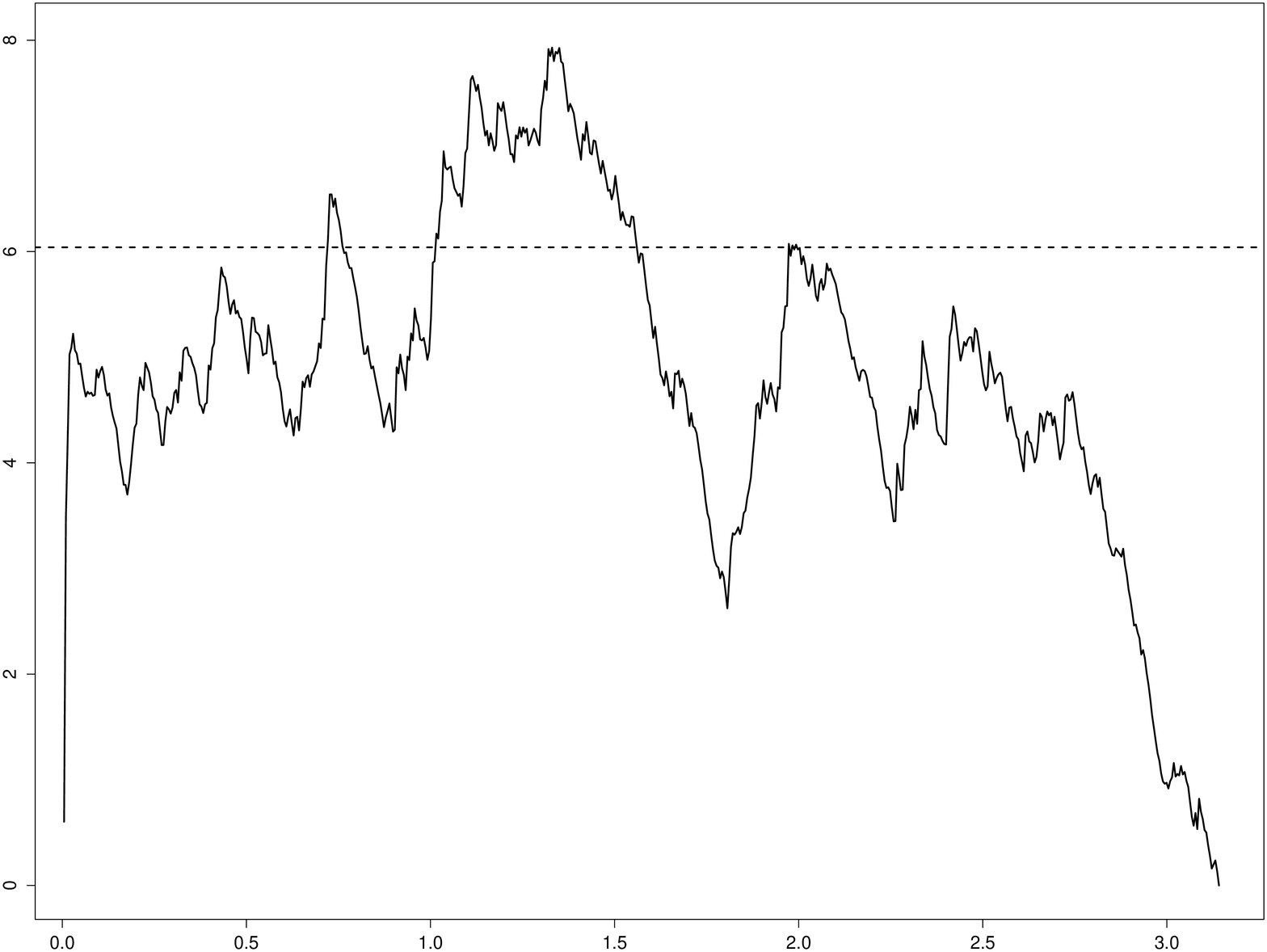}  
\end{minipage}%
\begin{minipage}[htbp]{0.4\linewidth}  
\centering  
\includegraphics[width=\textwidth]{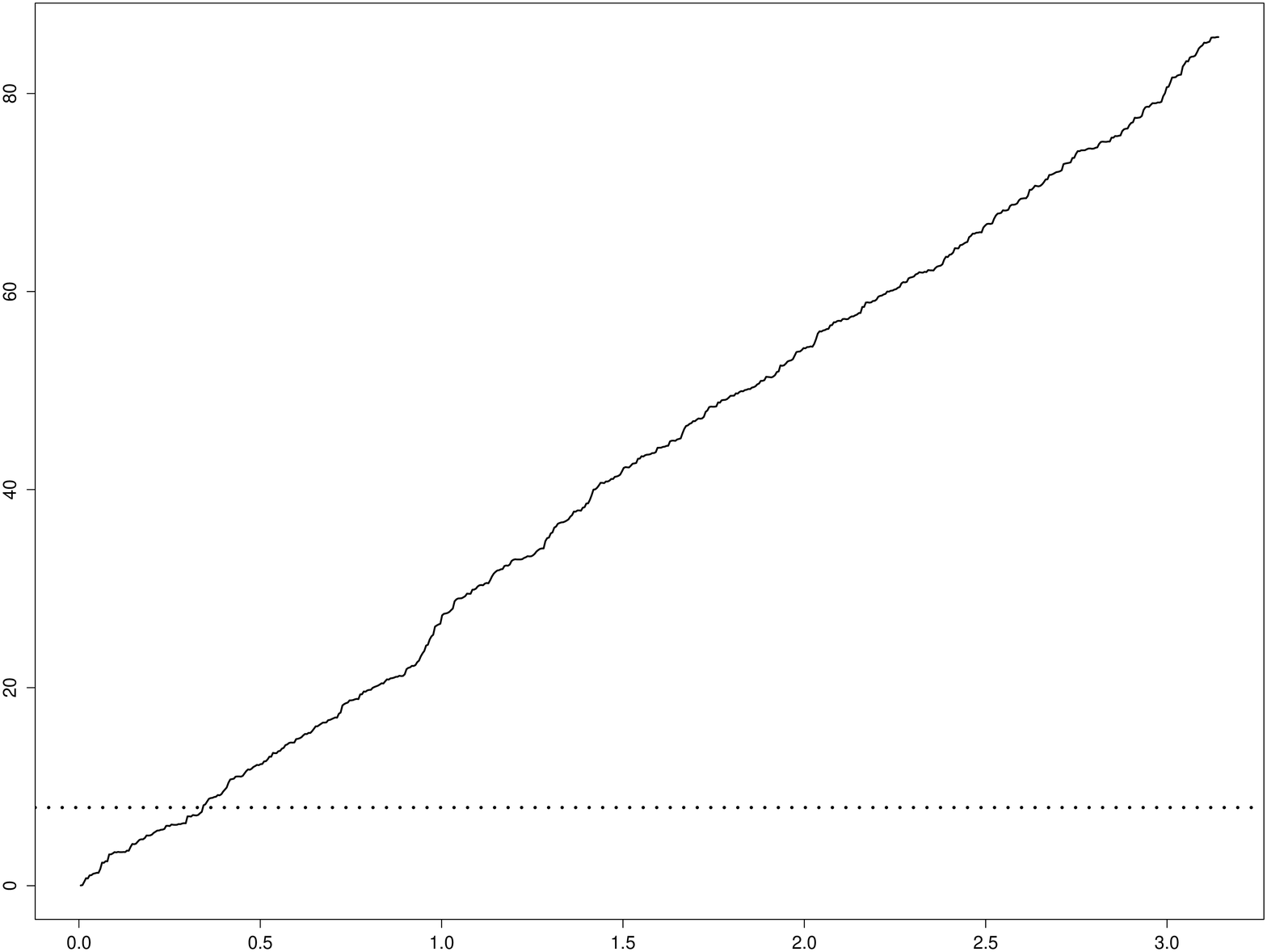}  
\end{minipage}%
\\
\begin{minipage}[ht]{0.4\linewidth}  
\centering  
\includegraphics[width=\textwidth]{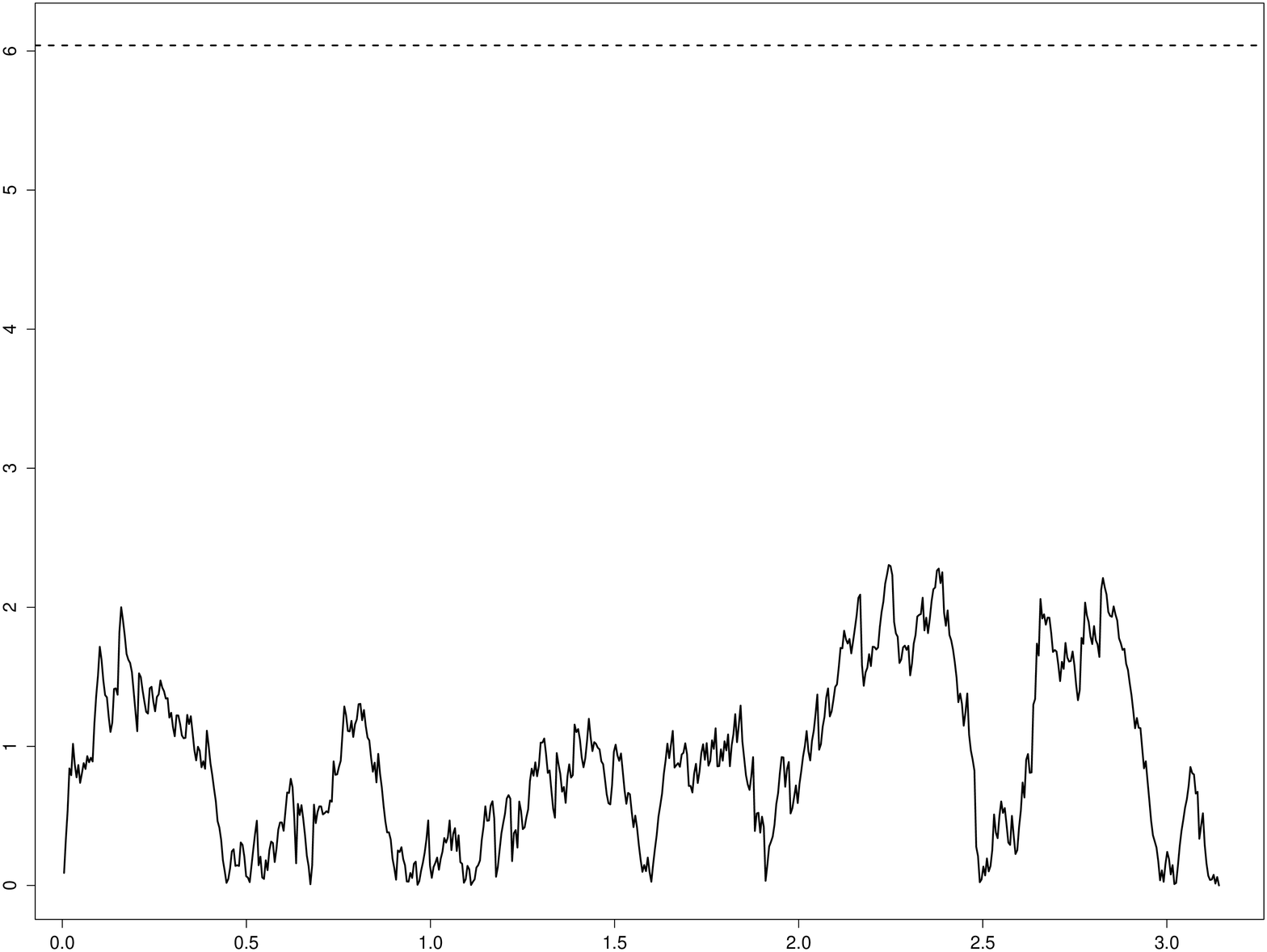}  
\end{minipage}%
\begin{minipage}[ht]{0.4\linewidth}  
\centering  
\includegraphics[width=\textwidth]{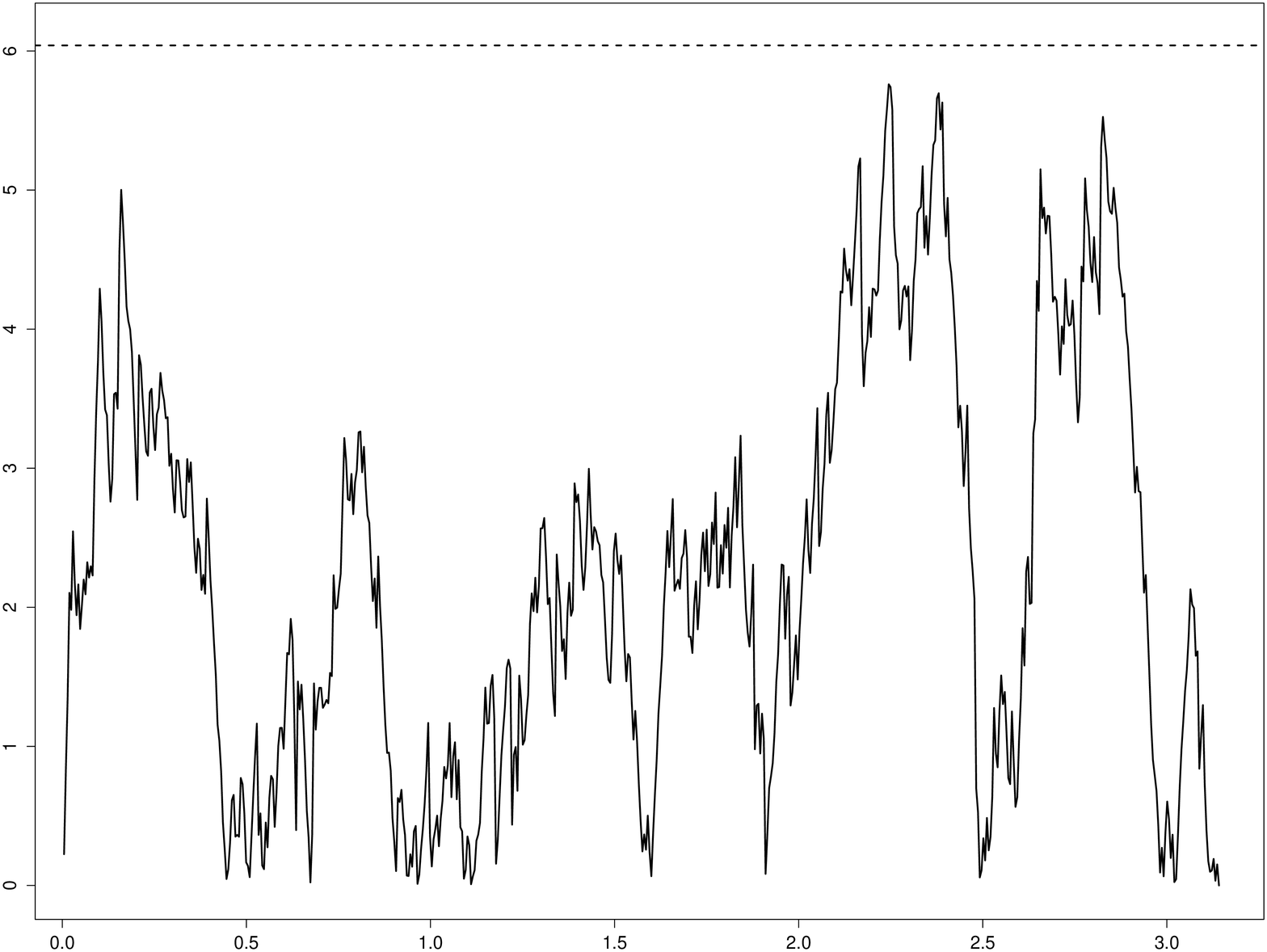}  
\end{minipage}%
\caption{\red \small
GRTs for daily Euro-USD FX rate
  log-returns 2002-2006(top, $n=1,280$), 
2006-2010 (middle, $n=1,279$), 2009-2013 (bottom, $n=1,281$).
The graphs show the integrated \per s  $n^{0.5}|J_{n,A}-\psi_0 \wt
\gamma_A(0)|$ under the null hypothesis of an iid \seq\ and 
($n/m)^{0.5}|J_{n,A}-EJ_{n,A}|$ in the case of a fitted GARCH model.
Under the iid hypothesis, the dotted lines represent the
95\%-quantile obtained from the limiting supremum of the absolute
values of a Brownian bridge. Under the GARCH hypothesis, the dotted
line represents the  bootstrap-based 
95\%-quantile of the GRS. 
{\bf Top:}  FX rate log returns 2002-2006 ($n=1,280$). We test under the iid
null hypothesis. For $p_0=0.05$ {\em (left)}, the null is not rejected. 
This is in contrast to the case $p_0=0.02$ {\em (right)} which leads to a 
clear rejection. The qualitative difference may be due to the
relatively small sample size which renders the test 
statistics meaningless. 
{\bf Middle}: FX rate log returns 2006-2010 ($n=1,279$). {\em Left.}
The iid null hypothesis with $p_0=0.05$ is rejected. {\em Right.} 
A \garch\ model with
 $\sigma_t^2=2.37\times 10^{-7}+0.1X_{t-1}^2+0.8\sigma^2_{t-1}$ and
 iid $t$-distributed noise with 4 degrees of freedom is fitted to the
 data. The null hypothesis of this GARCH is clearly rejected.
{\bf Bottom}:  FX rate log returns 2009-2013 ($n=1,281$).
The iid null hypothesis with $p_0=0.05$ {\em (left)} and $p_0=0.02$ 
{\em (right)}
is not  rejected. }
\label{fig:08}
\end{figure}

\section{Proof of Lemma~\ref{thm:app}}\label{sec:lemma}
{\bf Part 1.} Recall the series \rep s of 
$J_{n,A}(g)$ and $J_A(g)$ from  \eqref{eq:intper} and
\eqref{eq:intpera}, respectively. 
Then for every fixed $k\ge 1$, large $n$,
\beao
J_{n,A}(g)-J_A(g)&=& \Big(c_0(g) [\wt\gamma_A(0)-\gamma_A(0)]+2
  \sum_{h=1}^{k}c_h(g)\,[\widetilde{\gamma}_A(h)-\gamma_A(h)]\Big)\\ &&+
  2\sum_{h=k+1}^{n-1} c_h(g) \,[\wt \gamma_A(h)  -\gamma_A(h)]
-2 \sum_{h=n}^\infty c_h(g)\,\gamma_A(h)\\
&=& I_1(k)+I_2(k)-I_3\,.
\eeao
Then $I_3\to 0$ as $n\to\infty$ since $(\gamma_A(h))$ is summable and $EI_1(k)$ converges to zero as
$\nto$ 
due to
\regvar , for every $k$. In view of \eqref{eq:m1} in (M),
\beao
\Big|E\sum_{h=r_n+1}^{n-1}\widetilde{\gamma}_A(h) c_h (g)\Big| &= &
\Big|\frac m n \sum_{h=r_n+1}^{n-1}  (n-h)\, c_h(g)\,(p_h-p_0^2)\Big| \\ 
 &\le & c\, m\,\sum_{h=r_n+1}^{\infty}\xi_h \to 0\,, \quad n\to \infty\,, 
\eeao
and \eqref{eq:m2} in (M) implies
\beao
\lim_{k\to \infty} \limsup_{\nto} \Big|E
  \sum_{h=k+1}^{r_n}\widetilde{\gamma}_A(h) c_h(g)\Big|
 \le  c\,\lim_{k\to \infty} \limsup_{\nto} m \,\sum_{h=k+1}^{r_n}  p_h=0\,.
 \eeao
Since $\lim_{k\to \infty}\sum_{h=k+1}^\infty \gamma_A(h)=0$, we have
$\lim_{\kto}\limsup_{\nto} |EI_2(k)|=0$.
This proves Part 1.\\
{\bf Part 2.} 
It follows from Theorem 3.1
in Davis and Mikosch \cite{davis:mikosch:2009} that $\wt
\gamma_A(h)\stackrel{L^2}{\rightarrow} \gamma_A (h)$, $h\ge 1$. Hence 
$I_1(k)\stackrel{L^2}{\rightarrow} 0$ as $\nto$ for fixed $k\ge 1$.
It remains to show that 
$
\lim_{k\to \infty} \limsup_{\nto} \var(I_2(k))=0\,.
$
We have 
\beao
I_2(k)=2
\Big(\sum_{h=k+1}^{r_n}
+\sum_{h=r_n+1}^{n-1}\Big) c_h(g)
\big[\wt \gamma_A(h)-\gamma_A(h)\big]= 2 I_{21}(k)+ 2I_{22}\,. 
\end{eqnarray*}
In view of Lemma~\ref{lem:varcov} we get the bound
\begin{eqnarray*}
&&\var(I_{21}(k))
\le  \dfrac{m^2}{n}\sum_{h=k+1}^{r_n}
 \sum_{l=0}^{r_n-h}|c_h(g)c_{h+l}(g)|\times \\&& \Big(|\Gamma(0,h,0,h+l)|+
 \sum_{t=1}^{n-h-l}|\Gamma(0,h,t,t+h+l)|
  + \sum_{t=1}^{n-h}|\Gamma(0,h+l,t,t+h)| \Big)\\ 
 &= &Q_1+Q_2+Q_3\,.
\end{eqnarray*} 
Since $|c_h(g)|\le c/h$ (see \eqref{eq:fcoef}),
\begin{eqnarray*}
|Q_1| &\le& c\,
 \dfrac{m^2}{n} \sum_{h=k+1}^{r_n}|c_h(g)|\sum_{s=h}^{r_n} |c_s(g)| p_s
= c \dfrac{m^2}{n} \sum_{s=k+1}^{r_n}|c_s(g)| p_s \sum_{h=k+1}^{s}
|c_h(g)|\\
&\le & c \dfrac{m^2}{n} \sum_{s=k+1}^{r_n} p_s s^{-1}\log s \,,
\end{eqnarray*}
and the \rhs\ converges to 0 by first letting $\nto$ and then $\kto$, 
using  \eqref{eq:m2}.
Since the structures of
$Q_2$ and $Q_3$ are similar we restrict ourselves to showing 
$Q_2\to 0$ as $\nto,\kto$. We observe that
{\small \begin{eqnarray*}
|Q_2| &\le& c\frac{m^2}{n}\,\sum_{h=k+1}^{r_n} \sum_{s=h}^{r_n} 
\dfrac{1}{h s}\Big( \sum_{t=1}^{2r_n}+ \sum_{t=2r_n+1}^{n} \Big)|\Gamma(0,h,t,t+s)|\\ 
 & \le& c \frac{m \log^2 r_n}{n} m \sum_{h=k+1}^{3r_n}p_h+ c\,\frac{m
   \log ^2 r_n}{n} m\sum_{h=r_n+1}^n\xi_h +
c n^{-1}\Big(m\sum_{h=k+1}^{r_n} p_h/h\Big)^2\,.
\end{eqnarray*}}\noindent
In the last step, we used \eqref{eq:nb2}. The \rhs\ vanishes as $\nto$
and $\kto$.
Finally, we conclude that
$\lim_{k\to \infty} \limsup_{\nto}\var(I_{21}(k))=0$.
\par 
Now we turn to bounding $\var(I_{22})$.
In view of Lemma~\ref{lem:varcov} we have
\begin{eqnarray*}
\var(I_{22})
 &\le & \frac{m^2}{n}\sum_{h=r_n+1}^{n-1}
 \sum_{s=h}^{n-1}|c_h(g)c_{s}(g)| \Big(|\Gamma(0,h,0,s)|+
 \sum_{t=1}^{n-s}|\Gamma(0,h,t,t+s)|\\ 
 & &\quad  + 
\sum_{t=1}^{n-h}|\Gamma(0,s,t,t+h)| \Big) = Q_4+Q_5+Q_6\,.
\end{eqnarray*}
We have by \eqref{eq:fcoef},
\beao
Q_4&\le & c\frac{m^2}{n}\sum_{h=r_n+1}^{n-1}
 \sum_{s=h}^{n-1}|c_h(g)c_{s}(g)|\, |E\wt I_0 \wt I_s|\\
&\le & c\frac{m^2}{n}\sum_{h=r_n+1}^{n-1} h^{-2}
 \sum_{s=h}^{n-1}[(p_s-p_0^2) + p_0^2]\\
&\le & c\Big[\dfrac{m}{n r_n} m \sum_{h=r_n+1}^\infty \xi_h +
\dfrac{(p_0 m)^2}{r_n}\Big]=o(1)\,,\quad \nto\,.
\eeao
The terms $Q_5$ and $Q_6$ can be treated in a similar way; we focus on
$Q_5$. By \eqref{eq:fcoef},
\beao
Q_5&\le & c\frac{m^2}{n}\sum_{h=r_n+1}^{n-1}
 \sum_{s=h}^{h+r_n} (hs)^{-1}
 \sum_{t=1}^{r_n}|\Gamma(0,h,t,t+s)|\\
&&+c\frac{m^2}{n}\sum_{h=r_n+1}^{n-1}
  \sum_{s=h+1}^{n-1} 
 \sum_{t=r_n+1}^{n-s}
 (hs)^{-1}|\Gamma(0,h,t,t+s)|\\
&&+c\frac{m^2}{n}\sum_{h=r_n+1}^{n-1}
\sum_{s=h+r_n+1}^{n-1} 
 \sum_{t=1}^{r_n}
 (hs)^{-1}|\Gamma(0,h,t,t+s)|\\
&=& Q_{51}+Q_{52} + Q_{53}\,,
\eeao
and
\beao
Q_{51}&\le &c\frac{m^2}{n}\sum_{h=r_n+1}^{n-1}
 \sum_{s=h}^{h+r_n} (hs)^{-1} 
 \sum_{t=1}^{r_n} \big[(p_h-p_0^2) + p_0^2\big]\\
 &\le &c \Big(\dfrac{m}{n}m\sum_{h=r_n+1}^\infty\xi_h +
(mp_0)^2\dfrac{r_n }{n}\Big)\to 0\,,\quad\nto\,.   
\eeao
Next we consider  $Q_{52}$ and $Q_{53}$. By \eqref{eq:nb1}, we have
\beao
Q_{52}\le c\frac{2m^2}{n}\sum_{h=r_n+1}^{n-1}
 \sum_{s=h}^{n-1} (hs)^{-1} 
 \sum_{t=r_n+1}^n \xi_t
\le c\frac{m \log^2 n}{n} 
m \sum_{t=r_n+1}^\infty \xi_t\,.
\eeao
The right-hand side converges to zero by using the assumption $m\log^2
n/n=O(1)$ and the condition \eqref{eq:m1}. Similarly, using
\eqref{eq:nb1}, we obtain
\beao
Q_{53}\le c\frac{m }{n} m\sum_{h=r_n+1}^\infty \xi_{h}\,.
 \eeao
We conclude that $\var(I_{22})\to 0$ as $\nto$.

We proved above that $E(J_{n,A}-J_A(g))^2\to 0$, hence 
  $J_{n,A}(g)\overset{P}{\to} J_A(g)$, combined with 
  \eqref{eq:acf00},  yields $J_{n,A}^\circ (g)\overset{P}{\to}
  J_A^\circ (g)$. 

\section{Proof of
  Theorem~\ref{thm:main}}\label{sec:proof1}%\setcounter{equation}{0}
We start by proving \eqref{eq:fclt1}.
An application of the \cmt\ in $\bbc(\Pi)$ 
and Lemma~\ref{thm:acf} yield in $\bbc(\Pi)$ for every $k\ge 1$,
{\small
\beao\big(\dfrac mn\big)^{0.5}
 \big( \psi_0 \,(\wt \gamma_A(0)- E\wt \gamma_A(0))
+ 2 \sum_{h=1}^{k} \psi_h  \,(\wt \gamma_A(h)-E \wt \gamma_A(h))\big)
\std  \psi_0 Z_0 +2 \sum_{h=1}^{k} \psi_h Z_h\,.
\eeao}\noindent
Here $(Z_h)$ is a mean zero Gaussian process with covariance structure specified in Lemma~\ref{thm:acf}.
In view of Theorem 2 in 
Dehling et al.  \cite{dehling:durieu:volny:2009} relation
\eqref{eq:fclt1} will follow if we can prove the following result.
\ble \label{lem:fclt1} Assume that the conditions of Theorem~\ref{thm:main} hold.
Then for any  $\vep>0$,
\beao
\lim_{\kto}\limsup_{\nto} P\Big((n/m)^{0.5}\sup_{\la \in \Pi} \Big|\sum_{h=k+1}^{n-1} \psi_h(\la)  \,(\wt \gamma_A(h)-  E\wt \gamma_A(h))\Big |>\vep\Big)=0\,.
\eeao
\ele\noindent
\emph{Proof of Lemma \ref{lem:fclt1}.} We borrow the techniques of the
proof of Theorem~3.2 in Kl\"uppelberg and Mikosch \cite{kluppelberg:mikosch:1996}.
Without loss of generality
we assume that $k=2^a-1$ and 
$n=2^{b+1}$ where $a<b$ are integers; if 
$k$ or $n$ do not have this representation
we have to modify the proof slightly but we omit details.
For integer $q>0$ and some constant $\kappa>0$ to be chosen later, let
$\varepsilon_q=2^{-2q/\kappa}$.
We have for $\vep>0$,{\small
\begin{eqnarray*}
\lefteqn{Q=P\Big((n/m)^{0.5}\sup_{\la\in \Pi}\Big|\sum_{h=k+1}^{n-1} (\widetilde{\gamma}_A(h)-E\widetilde{\gamma}_A(h))\, \psi_h(\la) \Big|>\varepsilon \Big)}\\ 
&\le &P\Big((n/m)^{0.5} \sum_{q=a}^b\sup_{\lambda \in
  \Pi}\Big|\sum_{h=2^q}^{2^{q+1}-1} (\widetilde{\gamma}_A(h)
-E\widetilde{\gamma}_A(h))\,
\psi_h(\la) \Big|>\varepsilon )\\ 
 &\le & P\Big(\sum_{q=a}^b\varepsilon_q>\varepsilon \Big)+
P\Big( \bigcup_{q=a}^b\Big\{(n/m)^{0.5}\sup_{\la \in
  \Pi}\Big|\sum_{h=2^q}^{2^{q+1}-1}
(\widetilde{\gamma}_A(h)-E\widetilde{\gamma}_A(h))
 \psi_h(\la) \Big|>\varepsilon_q \Big\}\Big)\\
 &\le & \sum_{q=a}^bP\Big( (n/m)^{0.5}\sup_{\lambda \in
   \Pi}\Big|\sum_{h=2^q}^{2^{q+1}-1}(
 \widetilde{\gamma}_A(h)-E\widetilde{\gamma}_A(h))\, 
\psi_h(\la) \Big|>\varepsilon_q \Big)
=\sum_{q=a}^b Q_q\,.
\end{eqnarray*} }\noindent
In the last steps we used that $P(\sum_{q=a}^b\varepsilon_q>\varepsilon )$ vanishes for
fixed $\vep$ and sufficiently large
$a$. Next we will bound the expressions $Q_q$.
Write $J_{q,v}= \{(v-1)2^q+1,\ldots,v2^q\}$
and for $ j\in J_{q,v}$ and $\la\in [0,2^{-2 q}\pi]$,
\beao
Y_{qj}(\lambda)=(n/m)^{0.5}\sum_{h=2^q}^{2^{q+1}-1}(\widetilde{\gamma}_A(h)-E\widetilde{\gamma}_{A}
(h)) \psi_h (\lambda  +(j-1)\pi 2^{-2q}) \,.
\eeao
 Then
\beao
Q_q &= &P\Big((n/m)^{0.5}\max_{v=1,\ldots, 2^q} 
\max_{j\in J_{q,v}} \sup_{\lambda\in [(j-1)\pi 2^{-2q+1},j\pi
  2^{-2q+1}]}\\&&
\Big|\sum_{h=2^q}^{2^{q+1}-1} (\widetilde{\gamma}_A(h) -E\widetilde{\gamma}_A(h))\,\psi_h(\la) \Big|>\varepsilon_q \Big)\\ 
 &\le &
%\sum_{v=1,\ldots, 2^q} P\Big((n/m)^{0.5} \max_{j\in J_{qv}} \sup_{\lambda\in [(j-1)\pi 2^{-2q+1},j\pi 2^{-2q+1}]}\Big|\sum_{h=2^q}^{2^{q+1}-1}( \widetilde{\gamma}_{nA}(h)-E\widetilde{\gamma}_{nA}(h)) \psi_h(\la) \Big|>\varepsilon_q \Big)\\
%&=& 
\sum_{v=1}^{2^q} P\Big( ((n/m)^{0.5}\max_{j\in J_{q,v}}
 \sup_{\lambda\in [0, 2^{-2q+1}\pi]}\,|Y_{qj}(\la)|>\varepsilon_q  \Big) =  \sum_{v=1}^{2^q}Q_{qv}\,. 
\eeao 
We will bound each of the terms
$Q_{qv}$ by twice applying the maximal inequality 
of Theorem 10.2 in Billingsley~\cite{billingsley:1999}. For this
reason we have to control the variance of the increments of the
process 
$Y_{qj}$ both as a \fct\ of $\la$ and $j$. 
In particular, we will derive 
the following bound{\small
\beam \label{eq:maxpre}
 \frac n m\,E\Big(\sum_{h=2^q}^{2^{q+1}-1}( \widetilde{\gamma}_A(h)-E
\widetilde{\gamma}_A(h))\, d_h(\w,\la,j,j') \Big)^2  \le c\,|j-j'|^2 |\lambda-\omega|^{2\beta}\, K_{k,n}\,,
\eeam}\noindent
where $\beta$ is the H\"{o}lder coefficient of the function $g$, 
\beao
K_{k,n} \le c \Big[m \sum_{h=r_n+1}^{\infty}\xi_h + m
  \sum_{h=k+1}^{r_n} p_h + r_n/m\Big]
\eeao 
and for  $j<j'$ in $J_{q,v}$, $h\in \{2^q,\ldots,2^{q+1}-1\}$ and 
$\w <\la$ in $[0,2^{-2
  q+1}\pi]$,
\beam\label{dh}
\lefteqn{d_h(\w,\la,j,j')}\\
&=&\big(\psi_h(\lambda+(j'-1)\pi 2^{-2q+1})-\psi_h(\lambda+(j-1)\pi 2^{-2q+1})\big)\nonumber\\&&-\big(\psi_h(\omega+(j'-1)\pi 2^{-2q+1})-\psi_h(\omega+(j-1)\pi 2^{-2q+1})\big)\nonumber\\
 &=& \int_{\lambda+(j-1)\pi 2^{-2q+1}}^{\lambda+(j'-1)\pi 2^{-2q+1}}g(x)\cos(h x)dx-\int_{\omega+(j-1)\pi 2^{-2q+1}}^{\omega+(j'-1)\pi 2^{-2q+1}}g(x)\cos(h x)dx\nonumber\\ 
 & =& \int_{(j-1)\pi 2^{-2q+1}}^{(j'-1)\pi
   2^{-2q+1}}\Big(g(x+\lambda)[
 \cos(h(x+\lambda))-\cos(h(x+\omega))]\nonumber\\
&&\hspace{2.2cm}-[g(x+\lambda)-g(x+\omega)]\cos(h(x+\omega) ) \Big)\,dx\,.\nonumber
\eeam
Since $g$ is $\beta$-H\"{o}lder continuous we have{\small
\begin{eqnarray*}
 \Big|\int_{(j-1)\pi 2^{-2q+1}}^{(j'-1)\pi 2^{-2q+1}}
 [g(x+\lambda)-g(x+\omega)]\cos(h(\omega+x) )dx \Big|\le
 c(\lambda-\omega)^{\beta} (j'-j) 2^{-2q}\,.
\end{eqnarray*}}\noindent
Similarly,
\begin{eqnarray*}
\lefteqn{\Big|\int_{(j-1)\pi 2^{-2q}}^{(j'-1)\pi 2^{-2q}}g(x+\lambda)[ \cos(h(\lambda+x))-\cos(h(\omega+x))] dx \Big|}\\
&=& \Big|\int_{(j-1)\pi 2^{-2q}}^{(j'-1)\pi 2^{-2q}}g(x+\lambda)( 2\sin(h(\lambda-\omega)/2)\sin(h(\lambda+\omega+2x)/2)) dx \Big|\\
 & \le&c h(\lambda-\omega)(j'-j) 2^{-2q} \le c(\lambda-\omega)(j'-j) 2^{-q}\,.
\end{eqnarray*}
The last two inequalities yield for a constant $c$ only depending on $g$,
\beam\label{eq:101}
|d_h(\w,\la,j,j')|\le c|\lambda-\omega|^{\beta}\,|j'-j|\, 2^{-q}\,.
\eeam
Using this bound, we have
\beam
\lefteqn{\frac n m \,E\Big(\sum_{h=2^q}^{2^{q+1}-1}( \widetilde{\gamma}_A(h)-E
\widetilde{\gamma}_A(h))\, d_h(\w,\la,j,j') \Big)^2}\nonumber\\ 
% & \le & c\,\frac m n\sum_{h=2^q}^{2^{q+1}-1} \sum_{l=0}^{2^{q+1}-h-1} | d_h(\w,\la,j,j') d_{h+l}(\w,\la,j,j')| \Big|E( \widetilde{\gamma}_{nA}(h)-E
%\widetilde{\gamma}_{nA}(h))( \widetilde{\gamma}_{nA}(h+l)-E
%\widetilde{\gamma}_{nA}(h+l))\Big|\\ 
 &\le & c\,|j-j'|^{2} |\lambda-\omega|^{2\beta}\,2^{-2q}\frac nm  \sum_{h=2^q}^{2^{q+1}-1} \sum_{s=h}^{2^{q+1}-1}
\big|\cov\big(
\widetilde{\gamma}_A(h),\widetilde{\gamma}_A(s)\big)\big| \,.
\label{eq:100} 
\eeam
In what follows, it will be convenient to write $\sum_{h,l}^{(q)}
= \sum_{h=2^q}^{2^{q+1}-1} \sum_{l=0}^{2^{q+1}-h-1}$.
In view of Lemma~\ref{lem:varcov} we can bound the last term in
\eqref{eq:100}
as follows:{\small
\begin{eqnarray*}
&&\frac nm \sum_{h,l}^{(q)}
\big|\cov(\widetilde{\gamma}_A(h),
\widetilde{\gamma}_A(h+l))\big|\\
 &= & \frac{m}{n}\sum_{h,l}^{(q)}
  \Big|(n-h-l)\Gamma(0,h,0,h+l) +
 \sum_{t=1}^{n-h-l-1}(n-h-l-t)\Gamma (0,h,t,t+h+l)\\ 
 & & +\sum_{t=1}^{n-h-1} \min(n-h-l,n-h-t) \Gamma(0,h+l,t,t+h)\\&& -(n-h)(n-h-l)(p_h-p_0^2)(p_{h+l}-p_0^2) \Big|\\ 
 & \le & m \sum_{h,l}^{(q)} \Big[|\Gamma(0,h,0,h+l)|+ \sum_{t=1}^{h+r_n}|\Gamma(0,t,h,t+h+l)|+ \sum_{t=1}^{h+l+r_n}|\Gamma(0,h+l,t,t+h)|\\ 
 & &+\frac 1 n\Big|\sum_{t=h+r_n+1}^{n-h-l-1}(n-t-h-l)\Gamma(0,h,t,t+h+l)\\&&+\sum_{t=h+l+r_n+1}^{n-h-1}(n-t-h)\Gamma(0,h+l,t,t+h)\\&&-(n-h)(n-h-l)(p_h-p_0^2)(p_{h+l}-p_0^2) \Big|   \Big]\\
&=&W_1+W_2+W_3+W_4\,.
\end{eqnarray*}}\noindent
We will treat two cases of interest for the sums $\sum_{h,l}^{(q)}$: when
$2^{q+1}-1\le r_n$ and $2^q>r_n$. If
$2^q\le r_n<2^{q+1}-1$ the sums $\sum_{h,l}^{(q)}$ can be split into
two sums corresponding to $h\le r_n$ and $h>r_n$
and these can be treated in a similar fashion.
\par
We start by studying the case $2^{q+1}-1\le r_n$. Then $r_n\ge
2^{q+1}-1\ge h\ge 2^q>k$ and consequently $ 2^{q+1}-h-1\le 2^q$. Thus, 
$
W_1% m \sum_{h=2^q}^{2^{q+1}-1} \sum_{l=0}^{2^{q+1}-h-1} |\Gamma (0,h,0,h+l)| \\
\le c 2^q m \sum_{h=k+1}^{r_n}p_h\,.
$
The terms $W_2$, $W_3$ have a similar structure and can be
treated in the same way; we focus on $W_2$. 
Then we get the following bound from Lemma~\ref{lem:varcov}
\begin{eqnarray*}
W_2 % m \sum_{h,l}^{(q)}
%\sum_{t=1}^{h+r_n}|\Gamma(0,t,h,t+h+l)|
\le c\,2^{2q} \,\Big[m\sum_{h=k+1}^{r_n} p_h+m\sum_{h=r_n+1}^{2r_n}\xi_h+ (r_n/m)\Big]\,.
\end{eqnarray*}
In view of \eqref{eq:nb2}, we also have {\small
\begin{eqnarray*}
W_4 &\le & \frac{m}{n} \sum_{h,l}^{(q)}
%\Big|\sum_{t=h+r_n+1}^{n-h-l-1}(n-t-h-l)\Gamma(0,h,t,t+h+l)\\
%&&\qquad+\sum_{t=h+r_n+l+1}^{n-h-1}(n-t-h)\Gamma(0,h+l,t,t+h)
%-(n-h)(n-h-l)(p_h-p_0^2)(p_{h+l}-p_0^2) \Big|\\ 
% &\le & \frac{m}{n} \sum_{h=2^q}^{2^{q+1}-1} \sum_{l=0}^{2^{q+1}-h-1}
\Bigg[\sum_{t=h+r_n+1}^{n-h-l}(n-t-h-l)\big| \Gamma(0,h,t,t+h+l)-(p_h-p_0^2)(p_{h+l}-p_0^2)\big|\\
&&\qquad+\sum_{t=h+r_n+l+1}^{n-h}(n-t-h)\big|\Gamma(0,h+l,t,t+h)-(p_h-p_0^2)(p_{h+l}-p_0^2)\big|\\ 
 & & \qquad
+cnr_n\big|(p_h-p_0^2)(p_{h+l}-p_0^2)\big| \Bigg]\\
&\le &c 2^{2q} m\sum_{h=r_n+1}^n\xi_h + c\,\frac{r_n}{m}\Big( m\sum_{h=k+1}^{r_n}p_h\Big)^2\,.
\end{eqnarray*}}
Next we assume that $2^q>r_n$. By \eqref{eq:nb1} and \eqref{eq:nb2},
{\small \begin{eqnarray*}
W_1%&=& m \sum_{h=2^q}^{2^{q+1}-1} \sum_{l=0}^{2^{q+1}-h-1} |\Gamma
%(0,h,0,h+l)|\\ 
 &\le & m \sum_{h=2^q}^{2^{q+1}-1} \sum_{l=0}^{r_n} |\Gamma
 (0,0,h,h+l)-(p_0-p_0^2)(p_l-p_0^2)|+ 2^q m \sum_{l=0}^{r_n} (p_0-p_0^2)(p_l-p_0^2)\\
 && +m \sum_{h=2^q}^{2^{q+1}-1} \sum_{l=r_n+1}^{2^{q+1}-h-1} |\Gamma
 (0,0,h,h+l)|\\ 
 & & \\
&\le& c\,2^q m\sum_{h=r_n+1}^{\infty}\xi_h +  \frac{2^q r_n}{m}(mp_0)^2
\,.
\end{eqnarray*} }
We again focus on  $W_2$; $W_3$ can be treated in a  similar way. {\small
\begin{eqnarray*}
W_2 & \le& m\,\sum_{h=2^q}^{2^{q+1}-1}\Big(
\sum_{l=1}^{ r_n}
\Big( \sum_{t=1}^{r_n}+
 \sum_{t=r_n+1}^h 
+ \sum_{t=h+1}^{h+r_n}\Big)+
 \sum_{l=r_n+1}^{2^{q+1}-1}
\sum_{t=1}^{h+r_n}
\Big)\Big)|\Gamma(0,t,h,t+h+l)|\\
 &\le &
c \, 2^{2q} m \sum_{h=r_n+1}^{\infty}\xi_h  +c\,2^{2q} \frac{r_n}{m}(mp_0)^2%+
%c 2^{2q} m \sum_{h=r_n+1}^{\infty}\xi_h+ c 2^{2q}m \sum_{t=r_n+1}^{\infty}\xi_t%+c 2^{2q}m \sum_{h=r_n+1}^{\infty}\xi_h\,.
\end{eqnarray*}}\noindent
To obtain the bounds for $W_4$ we use \eqref{eq:nb2}: {\small 
\begin{eqnarray*}
W_4 &\le & \frac{m}{n} \sum_{h=2^q}^{2^{q+1}-1} \sum_{l=0}^{2^{q+1}-h-1}
\Bigg[\sum_{t=h+r_n+1}^{n-h-l}(n-t-h-l)\\&&\times \big| \Gamma(0,h,t,t+h+l)-(p_h-p_0^2)(p_{h+l}-p_0^2)\big|\\
&&+\sum_{t=h+r_n+l+1}^{n-h}(n-t-h)\big|\Gamma(0,h+l,t,t+h)-(p_h-p_0^2)(p_{h+l}-p_0^2)\big|\\ 
 & & 
+cn2^q\big|(p_h-p_0^2)(p_{h+l}-p_0^2)\big|\Bigg]\\ 
 & \le& c 2^{2q} m \sum_{t=r_n+1}^{\infty} \xi_t + c(2^q/m) \Bigg( m \sum_{h=r_n+1}^{\infty}\xi_h\Bigg)^2\,.
\end{eqnarray*}}
Collecting the bounds for $W_i$, $i\le 4$, and using \eqref{eq:100}, we finally proved 
\eqref{eq:maxpre}.
\par
{\red Using this bound,
we can apply 
the maximal inequality of
Theorem 10.2 in Billingsley~\cite{billingsley:1999} \wrt\ the variable
$\lambda \le 2^{-2q}\pi$ and for fixed $j,j'$}:
\begin{eqnarray*}
P(\max_{0\le \lambda \le
  2^{-2q}\pi}|Y_j(\lambda)-Y_{j'}(\lambda)|>\varepsilon_q)&\le& c
\varepsilon_q^{-2}  (2^{-2q}\pi)^{2\beta}\,(j-j')^2\,K_{k,n}\\
&\le & c\, 2^{4q (1-\beta+ \kappa^{-1})}\,((j-j')2^{-2q})^2\,K_{k,n}\,.
\end{eqnarray*}
Another application of this maximal inequality to  
$\max_{0\le \lambda \le 2^{-2q}\pi} | Y_j(\lambda)|$
\wrt\ the variable $j\in J_{q,v}$
yields
\begin{eqnarray*}
Q_{qv}=P\Big(\max_{j\in \{(v-1)2^q+1,\ldots,v2^q \}} \max_{0\le
  \lambda \le 2^{-2q}\pi} | Y_j(\lambda)| >\varepsilon_q \Big) \le c
 2^{4q (2^{-1}-\beta+\kappa^{-1})} K_{k,n}\,.
\end{eqnarray*}
Then we also have
\begin{eqnarray*}
Q_q\le \sum_{v=1}^{2^q}Q_{qv} \le c\, 2^{4q (3/4-\beta+\kappa^{-1})} K_{k,n}\,.
\end{eqnarray*}
The \rhs\ converges to zero as $q\to\infty$ provided $\beta\in
(3/4,1]$ and $\kappa$ is chosen sufficiently large.
Therefore we conclude for every $\vep>0$, 
\beam\label{eq:Q}
Q\le \sum_{q=a}^b Q_q \le c K_{k,n} \sum_{q=a}^{\infty}2^{4q (3/4-\beta+\kappa^{-1})}\,.
\eeam
The \rhs\ converges to zero by first letting $\nto$ and then
$k\to\infty$. This concludes the proof of \eqref{eq:fclt1}.
\par
Next we turn to the proof of \eqref{eq:fclt2}. It will follow from
\eqref{eq:fclt1} once we prove the  following lemma.
\ble \label{lem:fclt2}
Assume that the conditions of Theorem~\ref{thm:main} hold. Then
for any $\varepsilon>0$, as $\nto$,
\begin{eqnarray*}
P\Big(
(n/m)^{0.5}\sup_{\la\in\Pi}\big|\big(\widehat{J}_{n,A}(\la)-E\widehat{J}_{n,A}(\la)\big)-\big(J_{n,A}(\la)-EJ_{n,A}(\la)\big)\big|>\varepsilon
\Big)\to 0\,.
\end{eqnarray*}
\ele 
\noindent
\emph{Proof of Lemma \ref{lem:fclt2}:} For any fixed $k\ge 1$ we have
\begin{eqnarray*}
 \lefteqn{P\Big( (n/m)^{0.5}\sup_{\la\in\Pi
  }\big|\big(\widehat{J}_{n,A}(\la)-E\widehat{J}_{n,A}(\la)\big)-\big(J_{n,A}(\la)-EJ_{n,A}(\la)\big)\big|>\varepsilon
\Big)}\\ 
 &\le  & P\Big( (n/m)^{0.5}\sup_{\la\in
  \Pi}\Big|\sum_{h=0}^k
\big(\widetilde{\gamma}_A(h)-E\widetilde{\gamma}_A(h)\big)\big(\psi_h(\la)-\widehat{\psi}_h(\la)\big)\Big|>\varepsilon/3
\Big)\\
 & & +P\Big( (n/m)^{0.5}\sup_{\la\in\Pi
 }\Big|\sum_{h=k+1}^{n}
\big(\widetilde{\gamma}_A(h)-E\widetilde{\gamma}_A(h)\big)\psi_h(\la)\Big|>\varepsilon/3
\Big)\\
&&+P\Big( (n/m)^{0.5}\sup_{\la\in\Pi
 }\Big|\sum_{h=k+1}^{n}
\big(\widetilde{\gamma}_A(h)-E\widetilde{\gamma}_A(h)\big)\widehat{\psi}_h(\la)\Big|>\varepsilon/3)\\
 &= & V_1+V_2+V_3\,.
\end{eqnarray*} 
An application of Chebyshev's and H\"older's inequalities yields,
\beao
V_1 &\le & 9\vep^{-2}\frac{n}{m}E\sup_{\la\in
  \Pi}\Big|\sum_{h=0}^k
\big(\widetilde{\gamma}_A(h)-E\widetilde{\gamma}_A(h)\big)\big(\psi_h(\la)-\widehat{\psi}_h(\la)\big)\Big|^2\nonumber\\
&\le &c\, \dfrac nm \,E\sup_{\la\in
  \Pi}\sum_{h=0}^k \big(\widetilde{\gamma}_A(h)-E\widetilde{\gamma}_A(h)   \big)^2
\big|\psi_h(\la)-\widehat{\psi}_h(\la)\big|\,
\sum_{s=0}^k \big|\psi_s(\la)-\widehat{\psi}_s(\la)\big|\\
&\le &c\, k\dfrac nm \,\sum_{h=0}^k \var\big(\widetilde{\gamma}_A(h)\big)
\sup_{x\in
  \Pi}\big|\psi_h(x)-\widehat{\psi}_h(x)\big|
 \,.
\eeao
Next we will study $\sup_{\la\in
  \Pi}|\psi_h(\la)-\widehat{\psi}_h(\la)|$. Trivially, for $x\in\Pi$,
\beao
\Big|\int_{\w_n(x_n)}^x \cos (h\lambda) \,g(\la)\,d\la\Big|\le
c/n\,, 
\eeao
where the constant $c$ only depends on $g$. We also have for the
\fre ies $x\in \Pi$,{\small
\beam
\lefteqn{|\psi_h(\w_n(x_n))-\widehat{\psi}_h(\w_n(x_n))|}\nonumber\\
 &= & \Big|\sum_{i=1}^{x_n}\Big( \int_{\w_n(i-1)}^{\w_n(i)}\cos(h\lambda)\,g(\lambda)\,d\lambda - \w_n(1) \cos(h\w_n(i))\,g(\w_n(i))\Big)\Big|\nonumber\\ 
 & \le& \sum_{i=1}^{x_n}\Big|  \int_{\w_n(i-1)}^{\w_n(i)}\cos(h\lambda)\,(g(\lambda)-g(\w_n(i)))\,d\lambda\Big|\label{eq:96}\\
&& +\Big|\sum_{i=1}^{x_n} g(\w_n(i))\Big(\frac{\sin(h\w_n(i))-\sin(h\w_n(i-1))}{h}-\w_n(1)\cos(h\w_n(i)) \Big)\Big|\,.\nonumber\\\label{eq:97}
\eeam}\noindent
Since $g$ is $\beta$-H\"{o}lder continuous there exists 
a constant $c>0$ such that 
\beao
|g(\lambda)-g(\w_n(i))|\le c n^{-\beta}\,,\quad 
\la\in [\w_n(i-1),\w_n(i)]\,.
\eeao 
Therefore the term in \eqref{eq:96} is bounded by 
$c\,n^{-\beta}$. 
A Taylor expansion as $z\to 0$ yields
$\sin(z)=z-z^3/3!+o(z^3)$. Then we have for $h\le n$, {\small
\begin{eqnarray*}
&&\Big|\frac{\sin(h\w_n(i))-\sin(h\w_n(i-1))}{h}-\w_n(1)\cos(h\w_n(i))\Big|\\
 &= &\Big|2h^{-1}\sin(h\w_n(0.5))\cos(h\theta(i+0.5))-\w_n(1)\cos(h\w_n(i))\Big|\\ 
 & =&\Big|2h^{-1}(\sin(h\w_n(0.5))-h\w_n(0.5)) \cos(h\theta(i+0.5))\\&&+\w_n(1)(\cos(h\w_n(i+0.5))-\cos(h\w_n(i)))\Big|\\ 
 &\le & c(h\w_n(1))^3+
 \w_n(1)\Big|2\sin(h\w_n(0.25)\sin(h\w_n(i+0.25)) \Big|
 \le  c\, (h^3 n^{-3}+h n^{-2})\,.
\end{eqnarray*}}\noindent
Consequently, we have the bound  $c
(k/n)(1+k^2/n)$ for \eqref{eq:97} uniformly for $x\in\Pi$ and $h\le k$,
%\begin{eqnarray*}
%\sum_{i=1}^{x_n}\Big|
%g(\theta_n(i))\Big(\frac{\sin(h\theta_n(i))-\sin(h\theta_n(i-1))}{h}-\theta_n(1)\cos(h\theta_n(i))
%\Big)\Big|\le c r_n(n^{-1}+r_n^2/n^2)\,.
%\end{eqnarray*}
Thus, uniformly for $h\le k$,
\begin{eqnarray*}
\sup_{x\in\Pi} |\widetilde{\psi}_h(x)-\widehat{\psi}_h(x)| \le c\big[
n^{-\beta}+  (k/n)(1+k^2/n)]\,.
\end{eqnarray*}
As we have shown in Lemma~\ref{thm:acf}, $(n/m)\,\sum_{h=0}^k
\var\big(\widetilde{\gamma}_A(h)\big) \le c\,k$; see also 
Davis and Mikosch \cite{davis:mikosch:2009}, Lemma~5.2. Thus, as $\nto$,
\begin{eqnarray*}
V_1\le c\big[
k^2n^{-\beta}+  (k^3/n)(1+k^2/n)]\to 0\,.
\end{eqnarray*}
\par
It follows from Lemma \ref{lem:fclt1} that $\lim_{k\to \infty}
\limsup_{\nto}V_2=0$. We adapt the proof of  Lemma \ref{lem:fclt1}
for the case $V_3$. Abusing notation, 
consider
\begin{eqnarray*}
d_h(\omega,\lambda,j,j')
&=&\big(\widehat{\psi}_h(\lambda+(j'-1)\pi 2^{-2q+1})-\widehat{\psi}_h(\lambda+(j-1)\pi
2^{-2q+1})\big)\\
&&-\big(\widehat{\psi}_h(\omega+(j'-1)\pi 2^{-2q+1})-\widehat{\psi}_h(\omega+(j-1)\pi 2^{-2q+1})\big)\,.
\end{eqnarray*}
Recall that
we assume $n=2^b$ for some integer $b$ and $x_n=[nx/(2\pi)]$. 
Therefore for $\la\in\Pi$ and
integer $j$,
\beao
(\la+ (j-1) \pi 2^{-2q+1})_n&=& [n\la/(2\pi)+(j-1) 2^{-2q+b} ]\\&=&[n\la/(2\pi)]+(j-1) 2^{-2q+b}\\&=&\la_n+(j-1)  2^{-2q+b}\,.
\eeao
Thus we can write{\small
\begin{eqnarray*}
&&d_h(\omega,\lambda,j,j')\\&=
 &\dfrac{2\pi}n \sum_{i=\lambda_n+(j-1) 2^{b-2q}}^{\lambda_n+(j'-1)
   2^{b-2q}} g(\w_n(i))\cos (h\w_n(i))\\&&
  -\dfrac{2\pi} n\sum_{i=\omega_n+(j-1)2^{b-2q}}^{\omega_n+(j'-1)
   2^{b-2q}} g(\w_n(i))\cos (h\w_n(i))\\ 
 & =& \dfrac{2\pi} n\sum_{i=(j-1)2^{b-2q}}^{(j'-1) 2^{b-2q}}
 \Big(g(\w_n(\lambda_n+i))[\cos
 (h\w_n(\lambda_n+i))-\cos(h\w_n(\omega_n+i))]\\ 
 &
 &\hspace{2.2cm}-[g(\w_n(\lambda_n+i))-g(\w_n(\omega_n+i))]\cos(h\w_n(\omega_n+i)) \Big) =T_1+T_2\,.
\end{eqnarray*}}\noindent
Calculation yields 
\begin{eqnarray*}
|T_1|& \le&
c|\w_n(\lambda_n)-\w_n(\omega_n)|\big|(j'-j)2^{-2q}\big|2^q
\le c\big|(\lambda_n-\omega_n)/n\big|\big|(j'-j)2^{-2q}\big| 2^q\,,\\
|T_2|&\le & 
c|\w_n(\lambda_n)-\w_n(\omega_n)|^{\beta}\big|(j'-j)2^{-2q}\big|\le
c\big|(\lambda_n-\omega_n)/n\big|^{\beta}\big|(j'-j)2^{-2q}\big|\,.
\end{eqnarray*}
Combining these bounds, we have, 
\begin{eqnarray*}
|d_h(\omega,\lambda,j,j')|\le c|(\lambda_n-\omega_n)/n|^{\beta}\big|(j'-j)2^{-2q}\big|2^q\,.
\end{eqnarray*}
In the remaining argument we can follow the proof of Lemma
\ref{lem:fclt1}; the only difference is that we have to replace the
supremum over $\la,\omega \in [0, j 2^{-2q+1}]$ 
by the corresponding quantities $\la_n/n,\omega_n/n \in [0, j 2^{-2q+1}]$. 
This proves\\$\lim_{k\to \infty}\limsup_{\nto}V_3=0$ and concludes the
proof of the lemma.
\par
The proofs of \eqref{eq:fclt01} and \eqref{eq:fclt02} are 
completely analogous. 
Instead of the relations \eqref{eq:acf11} one has to use \eqref{eq:opu2}.

\section{Proof of Theorem~\ref{thm:mdep2}}%\setcounter{equation}{0}
\label{sec:proof2}
We adapt the proof
of Theorem~\ref{thm:main}. We need to prove that 
\begin{eqnarray*}
n \sum_{h,l}^{(q)}
|\cov(\widetilde{\gamma}_A(h),\widetilde{\gamma}_A(h+l))| \le c2^q\,,
\end{eqnarray*}
where $\sum_{h,l}^{(q)}$ is defined in the proof of
Lemma~\ref{lem:fclt1}. Here $h>\eta$. {\small
\begin{eqnarray*}
\lefteqn{n\sum_{h,l}^{(q)}
\big|\cov(\widetilde{\gamma}_A(h),
\widetilde{\gamma}_A(h+l))\big|}\\
 &= & \frac{m^2}{n}\sum_{h,l}^{(q)}
  \Big|(n-h-l)\Gamma(0,h,0,h+l) +
 \sum_{t=1}^{n-h-l-1}(n-h-l-t)\Gamma (0,h,t,t+h+l)\\ 
 & & +\sum_{t=1}^{n-h-1} \min(n-h-l,n-h-t) \Gamma(0,h+l,t,t+h) \Big|\\ 
 &\le &  m^2 \sum_{h=2^q}^{2^{q+1}-1} |\Gamma(0,0,h,h)| + m^2
 \sum_{h=2^q}^{2^{q+1}-1} \sum_{l=1}^{2^{q+1}-h-1} |\Gamma(0,0,h,h+l)|\\
 & &  + m^2  \sum_{h,l}^{(q)} \sum_{t=1}^{n-h-l-1}
 |\Gamma(0,h,t,t+h+l)|+ m^2  \sum_{h,l}^{(q)} \sum_{t=1}^{n-h-l-1}
 |\Gamma(0,h+l,t,t+h)|\\ 
 &=  & m^2 \sum_{h=2^q}^{2^{q+1}-1} |\Gamma(0,0,h,h)|+m^2
 \sum_{h=2^q}^{2^{q+1}-1} \sum_{l=1}^{\eta}|\Gamma(0,0,h,h+l)|\\ 
 & & + m^2   \sum_{h=2^q}^{2^{q+1}-1} \sum_{l=1}^{\eta} \sum_{t=1}^{\eta}
 |\Gamma(0,h,t,t+h+l)|+m^2   \sum_{h=2^q}^{2^{q+1}-1}
 \sum_{l=1}^{\eta} \sum_{t=1}^{\eta}|\Gamma(0,h+l,t,t+h)|\\
 &\le & c2^q
\end{eqnarray*} }
In the above calculation, we use the facts that for $s\le t\le u\le v$,
$\Gamma(s,t,u,v)=0$ where $t-s>\eta$ or $v-u>\eta$. 

In the remaining argument we can follow the proof of
Theorem~\ref{thm:main}; instead of Lemma~\ref{thm:acf} we use the central limit
theory of Lemma~\ref{thm:mdep1}. \qed

\section{Proof of Theorem~\ref{thm:bpclt}}\label{sec:proof3} 
We will mimic the proof of Theorem~\ref{thm:main}. 
We start by proving a result for the bootstrapped sample extremogram
$\wh \gamma_A^\ast$
analogous to Theorem~\ref{thm:bp1}.
\ble \label{lem:bpclt1}
Under the conditions and with the notation of Theorem~\ref{thm:bp1},
for $h\ge 0$, 
\beao
d\Big((n/m)^{0.5} \big(
\wh{\gamma}_A^{*}(i)-E^{*}\wh{\gamma}_A^{*}(i)
\big)_{i=0,\ldots,h},(Z_i)_{i=0,\ldots,h}\Big)  \stp 0\,,\quad \nto\,.
\eeao
\ele
\begin{proof}
We start by observing (see Lemma~\ref{lem:moment})
that for $h\ge 0$ 
\beao
E^{*}\wt{\gamma}_A^{*}(h)&=& \dfrac{m}n \,(n-h)E^\ast 
\wt I_{1^\ast}\wt I_{1^\ast+h}\\&= &(1-h/n) \,\Big[\wt \gamma_A
(h)+\frac mn \sum_{t=n-h+1}^n \wt I_t\wt  I_{t+h}\Big] \,,
\\
E^{*}\wh{\gamma}_A^{*}(h)&=& \dfrac{m}n \,(n-h)E^\ast 
\wh I_{1^\ast}\wh I_{(1+h)^\ast}\\&=&
 (1-h/n) \,(1-\theta)^h\Big[\wh \gamma_A (h)
+\frac mn\sum_{t=n-h+1}^n \wh I_t\wh  I_{t+h}\Big]
\,,
\eeao
where we interpret indices larger than $n$ modulo $n$,
and therefore
\beam\label{eq:a3}
(n/m)^{0.5}\big[ {(1-\theta)^h}
E^{*}\wt{\gamma}_A^{*}(h)-E^{*}\wh{\gamma}_A^{*}(h)\big]=O_P(m^{-1})\stp 0\,,
\eeam
where we used that $\ov I_n^2-p_0^2=O_P(1/\sqrt{mn})$.
By virtue of Theorem~\ref{thm:bp1} 
it suffices to show that for any $\vep>0$ and $h\ge 0$, as $\nto$, 
\beao
P^\ast
\Big((n/m)^{0.5}\Big|{(1-\theta)^h}\big(\widetilde{\gamma}_A^{*}(h)-E^{*}\widetilde{\gamma}_A^{*}(h)\big)-\big(\widehat{\gamma}_A^{*}(h)-E^{*}\widehat{\gamma}_A^{*}(h)\big)\Big|>\vep\Big)\stp
0\,.
\eeao  
Markov's inequality ensures that it suffices to prove
that 
\begin{eqnarray*}
\frac{n}{m}
\var^\ast\big({
  (1-\theta)^h}\widetilde{\gamma}_A^{*}(h)-\widehat{\gamma}_A^{*}(h)\big) \ \stp 0\,,\quad \nto\,.
\end{eqnarray*}
We observe that {\small
\beao\lefteqn{
\dfrac nm \var^\ast \big({(1-\theta)^h}\wt{\gamma}_A^{*}(h)-\wh{\gamma}_A^{*}(h)
\big)}\\
&=& m \big(1-\frac h n\big)  \var^\ast
\Big(\widehat{\Imath}_{1^{*}}\widehat{\Imath}_{(1+h)^{*}}-{(1-\theta)^h}\wt
I_{1^{*}}\wt I_{1^{*}+h}
 \Big)\\
&&+ 2 m \sum_{s=1}^{n-h-1} \big(1-\frac{h+s}n\big)\times \\&&\cov^\ast \Big(\wh
I_{1^\ast}\widehat{\Imath}_{(1+h)^{*}}-{(1-\theta)^h}\wt
I_{1^{*}}\wt I_{1^{*}+h}, \wh I_{(1+s)^\ast}\widehat{\Imath}_{(1+s+h)^{*}}-{ (1-\theta)^h}\wt
I_{(1+s)^\ast}\wt I_{(1+s)^\ast+h}\Big)\\
 &= &  m \big(1-\frac h n \big)  \var^\ast
\Big(\widehat{\Imath}_{1^{*}}\widehat{\Imath}_{(1+h)^{*}}-{(1-\theta)^h}\wt
I_{1^{*}}\wt I_{1^{*}+h}
 \Big)\\
 & & + 2  m \sum_{s=1}^{n-h-1} \big(1-\frac{h+s}n\big)\,\Big[\cov^\ast (\wh
I_{1^\ast}\widehat{\Imath}_{(1+h)^{*}}, \wh
I_{(1+s)^\ast}\widehat{\Imath}_{(1+s+h)^{*}})\\&&-{(1-\theta)^h}
\cov^\ast (\wh I_{1^\ast}\widehat{\Imath}_{(1+h)^{*}}, \wt
I_{(1+s)^\ast}\wt I_{(1+s)^\ast+h})\\&&
-{(1-\theta)^h}\cov^\ast (\wt
I_{1^{*}}\wt I_{1^{*}+h}, \wh
I_{(1+s)^\ast}\widehat{\Imath}_{(1+s+h)^{*}}) \\&&
+{(1-\theta)^{2h}}\cov^\ast (\wt
I_{1^{*}}\wt I_{1^{*}+h},\wt I_{(1+s)^\ast}\wt I_{(1+s)^\ast +h})\Big] 
  = Q_1+Q_2\,.
\end{eqnarray*} }\noindent
We will show that the \rhs\ converges to zero in $P$-\pro y,
where we focus on $Q_2$ and omit the details for $Q_1$.
We start by looking at the summands in $Q_2$ for fixed $s\le h$, using
the structure of the covariances in
Lemma~\ref{lem:moment}. The expressions for
the covariances in Lemma~\ref{lem:moment} contain terms with
normalization $n^{-2}$. For example,  by \eqref{eq:mo30}  a
corresponding term in $Q_2$ is of the order  
\beao
m \Big(n^{-1}\sum_{i=1}^n \wt I_i\wt I_{i+h}\Big)^2=
m^{-1} \Big(\frac m n \sum_{i=1}^n \wt I_i \wt I_{i+h} \Big)^2=O_P(m^{-1})\,,
\eeao
since $\frac m n \sum_{i=1}^n \wt I_i \wt I_{i+h}\stp \gamma_A(h)$;
see Lemma~\ref{thm:acf}. In the latter sums, the $\wt I_i$'s can be
exchanged by the $I_i$'s or the $\wh I_i$'s. Therefore all other 
terms in $Q_2$ with normalization $mn^{-2}$ converge to zero in $P$-\pro y.
Another appeal to  Lemma~\ref{lem:moment} shows that
it remains to consider those expressions in $Q_2$ that are
normalized by $m n^{-1}$ again for fixed $s\le h$. From \eqref{eq:mo32} and
\eqref{eq:mo4} we see that, on one hand,  we have to deal with the
differences
{\small
\beam\label{eq:r0}
(1-\theta)^{s+h}\; \,\frac mn
\sum_{i=1}^n\wh I_i\wh I_{i+s}\wh I_{i+h}\wh I_{i+s+h}-
(1-\theta)^{s+2h}\frac mn
\sum_{i=1}^n\wt I_i\wt I_{i+h}\wh I_{i+s}\wh I_{i+s+h}\,,
\eeam}\noindent
but both sums are consistent estimators of
$\lim_{\nto} m P(a_m^{-1}X_0\in A,a_m^{-1}X_s\in A,a_m^{-1}X_h\in
A,a_m^{-1}
X_{s+h}\in A)$ (see \cite{davis:mikosch:2009}, Theorem~3.1).
Therefore \eqref{eq:r0} converges to zero in
$P$-\pro y. On the other hand, in view of  \eqref{eq:mo30} and 
\eqref{eq:mo31} we have to deal with the differences, for~$s\le h$, 
\beao
{(1-\theta)^{s+2h}} \frac mn
 \sum_{i=1}^n\widetilde{\Imath}_i
\widetilde{\Imath}_{i+s}\widetilde{\Imath}_{i+h}\widetilde{\Imath}_{i+s+h}-
{(1-\theta)^{2h}} \;   \frac mn \sum_{i=1}^n
 \widehat{\Imath}_i \widehat{\Imath}_{i+h}
 \widetilde{\Imath}_{i+s}\widetilde{\Imath}_{i+s+h}\,,
\eeao
which again converge to zero in $P$-\pro y.
These arguments finish the proof for $s\le h$.
\par 
An inspection of the covariances in 
Lemma~\ref{lem:moment} shows that for $s>h$ all expressions 
with normalization $n^{-2}$ do not depend on $s$. The corresponding 
aggregated terms in $Q_2$ are then given by {\small
\beao
&&2m \sum_{s=h+1}^{n-h-1}\big(1-\frac{h+s} n\big)
\Big[-(1-\theta)^{s+h}
  \Big(n^{-1}\sum_{i=1}^n \widehat{\Imath}_i\widehat{\Imath}_{i+h}
  \Big)^2\\&&
+(1-\theta)^{s+h} \Big(n^{-1}\sum_{i=1}^n\widehat{\Imath}_i\widehat{\Imath}_{i+h} \Big)
 \Big(n^{-1}\sum_{i=1}^n\widetilde{\Imath}_i\widetilde{\Imath}_{i+h}
 \Big)\Big)
\\&&+
{(1-\theta)^{s+2h}} \; 
 \Big(n^{-1}\sum_{i=1}^n\widehat{\Imath}_i\widehat{\Imath}_{i+h} \Big)
 \Big(n^{-1}\sum_{i=1}^n\widetilde{\Imath}_i\widetilde{\Imath}_{i+h} \Big)
-
{(1-\theta)^{s+2h}}\;
  \Big(n^{-1}\sum_{i=1}^n \wt I_i\wt I_{i+h} \Big)^2
\Big]\\
&=&- 2m^{-1}\Big(\frac mn
\sum_{i=1}^n\widehat{\Imath}_i\widehat{\Imath}_{i+h} 
-\frac  mn
\sum_{i=1}^n\wt{\Imath}_i\wt{\Imath}_{i+h}\Big) \Big(\frac mn \sum_{i=1}^n\widehat{\Imath}_i\widehat{\Imath}_{i+h} \Big) 
\sum_{s=h+1}^{n-h-1}\big(1-\frac{h+s} n\big)\; (1-\theta)^{s+h}\\
&&- 2m^{-1}\Big(\frac mn
\sum_{i=1}^n\wt{\Imath}_i\wt{\Imath}_{i+h} 
-\frac  mn
\sum_{i=1}^n\wh{\Imath}_i\wh{\Imath}_{i+h}\Big) \Big(\frac mn \sum_{i=1}^n\wt{\Imath}_i\wt{\Imath}_{i+h} \Big) 
\sum_{s=h+1}^{n-h-1}\big(1-\frac{h+s} n\big)\; {(1-\theta)^{s+2h}}\\
&=& O_P( 1/ (\theta \sqrt{mn}))= o_P(1)\,.
\eeao  }\noindent
In the last step we used \eqref{eq:a3} and
the assumption $n\theta^2/m\to\infty$.
Finally, we deal with the remaining terms in  $Q_2$. In
view of Lemma~\ref{lem:moment} they are given by {\small
\beao
&&2m \sum_{s=h+1}^{n-h-1}\big(1-\frac{h+s} n\big)
\Big[(1-\theta)^{s+h}\;n^{-1}\sum_{i=1}^n \widehat{\Imath}_i\widehat{\Imath}_{i+s}
  \widehat{\Imath}_{i+h} \widehat{\Imath}_{i+s+h}\\&&
-(1-\theta)^{s+h} \; n^{-1} \sum_{i=1}^n
 \widehat{\Imath}_i \widehat{\Imath}_{i+h}
 \widetilde{\Imath}_{i+s}\widetilde{\Imath}_{i+s+h}\\&&
-{(1-\theta)^{s+2h}} \;   n^{-1} \sum_{i=1}^n
 \widetilde{\Imath}_i \widetilde{\Imath}_{i+h}
 \widehat{\Imath}_{i+s}\widehat{\Imath}_{i+s+h}
+{(1-\theta)^{s+2h}}\;n^{-1}\sum_{i=1}^n  \wt I_i\wt{\Imath}_{i+s}
  \wt{\Imath}_{i+h} \wt{\Imath}_{i+s+h}
\Big]\\
&=&2m \sum_{s=h+1}^{n-h-1}\big(1-\frac{h+s} n\big)
(1-\theta)^{s+h} n^ {-1} \sum_{i=1}^n\
 \wh{\Imath}_i \wh{\Imath}_{i+h}
 \big(
 \wh{\Imath}_{i+s} \wh{\Imath}_{i+s+h}-
 \wt{\Imath}_{i+s}\wt{\Imath}_{i+s+h}\big)\\
&+&2m \sum_{s=h+1}^{n-h-1}\big(1-\frac{h+s} n\big)
(1-\theta)^{s+2h}\;n^{-1}\sum_{i=1}^n \wt{\Imath}_i\wt{\Imath}_{i+h}\big(
  \wt{\Imath}_{i+s} \wt{\Imath}_{i+s+h}-\wh{\Imath}_{i+s} 
\wh{\Imath}_{i+s+h}\big)=J_0\,.
\eeao}\noindent
Using the assumption $n\theta^2 /m\to\infty$,
we have 
\beao 
E|J_0|&\le& c\, m\,\sum_{s=h+1}^{n-h-1}(1-\theta)^{s+h} E|\wh I_0 \wh
I_h-\wt I_0\wt I_h|\\&\le&  c\, m  E|p_0-\ov
I_n|\sum_{s=h+1}^{n-h-1}(1-\theta)^{s+h}\le c (m/n)^{0.5}\theta^{-1}=o(1)\,.
\eeao
This finishes the proof of the lemma.\qed
\end{proof}
We conclude from Lemma~\ref{lem:bpclt1} that for any $k\ge 1$, as $\nto$,
{\small \beao
&&d \Big((n/m)^{0.5}\Big(
\psi_0 \,\big(\wh \gamma_A^\ast(0)-E^\ast\wh \gamma_A^\ast(0)\big) +2 
\sum_{h=1}^{k} \psi_{h}\,\big(\wh \gamma_A^\ast(h)
-E^\ast \wh \gamma_A^\ast(h)\big)
\Big),\\&&\quad
\psi_0 \,Z_0+2 
\sum_{h=1}^{k} \psi_{h}\,Z_h
\Big)\stp 0\,,
\eeao }\noindent
where the dependence structure of $(Z_h)$ is 
defined in Lemma~\ref{thm:acf}.
\par
The proof of the theorem is finished   by the following result
which parallels 
Lemma~\ref{lem:fclt1}.
\ble \label{lem:bpclt3}
Assume the conditions of Theorem~\ref{thm:bpclt}. Then the
following relation holds for $\delta>0$
\beam
\label{eq:bpr1}\red
\lim_{k\to \infty} \limsup_{\nto}P\Big( (n/m)^{0.5}
\sup_{\lambda\in \Pi}
\Big|\sum_{h=k+1}^{n-1}\psi_h(\lambda)\,
\big(\wh{\gamma}_A^{*}(h)-E^{*}\wh{\gamma}_A^{*}(h)\big)\Big|>\delta\Big)
 = 0\,.\nonumber\\
\eeam 
\ele

\begin{proof}
We follow the lines of the proof of Lemma~\ref{lem:fclt1} and use the
same notation.
We again assume without loss of generality that $k=2^a-1$ and
$n=2^{b+1}$ for integers  $a<b$, $a$ chosen sufficiently large,  and
we write $\varepsilon_q=2^{-2q/\kappa}$ for $\kappa>0$ to be chosen
later. Then, for large $a$ depending on $\vep>0$, the steps of the
proof 
lead to the inequality (cf. \eqref{eq:Q})
\beao
Q^{*} & =& P^{*}\Big((n/m)^{0.5} \sup_{\lambda\in \Pi}
\Big|\sum_{h=k+1}^{n-1} \big(\widehat{\gamma}_A^{*}(h)-E^{*}\widehat{\gamma}_A^{*}(h)\big)\,
\psi_h(\lambda)\Big|>\varepsilon \Big)\\ 
&\le &
c\sum_{q=a}^b 2^{4q(0.75-\beta+\kappa^{-1})}K_{q}\,,
\eeao
where $\beta\in (3/4,1]$ is the H\"older coefficient of the \fct\ $g$,
the number $\kappa>0$ can be chosen arbitrarily large and
\beao
  K_{q}= \frac nm \sum_{h=2^q}^{2^{q+1}-1}
  \sum_{s=h}^{2^{q+1}-1}|\cov^\ast(\wh\gamma_A^\ast(h),\wh \gamma_A^\ast(s)|\,.
\eeao
By the Cauchy-Schwarz inequality, for $s,h\in [2^q,2^{q+1})$ and $h\le
s$,
\beao
(n/m)^2|\cov^\ast(\wh\gamma_A^\ast(h),\wh \gamma_A^\ast(s)|^2\le 
(n/m) \var^\ast(\wh\gamma_A^\ast(h))\;(n/m) \var^\ast(\wh\gamma_A^\ast(s))\,.
\eeao
We will show that 
\beam\label{eq:opi}
(n/m) E \var^\ast(\wh\gamma_A^\ast(h))\le c
\eeam
 for
some constant $c$, uniformly for $k\le h\le n$ and $n$. Then 
\beao
EQ^\ast \le c\, \sum_{q=a}^b 2^{4q(3/4-\beta+\kappa^{-1})} \le c\,
\sum_{q=a}^ \infty 2^{4q(3/4-\beta+\kappa^{-1})}\,.
\eeao
The \rhs\ converges since $\beta \in (3/4,1]$ and $\kappa$ can be
chosen arbitrarily large. Moreover, the \rhs\ converges to zero as $k\to\infty$.
\par
Thus it remains to show \eqref{eq:opi}.
In view of Lemma~\ref{lem:moment} we have {\small
\begin{eqnarray*}
\lefteqn{(n/m)E\var^{*}(\widehat{\gamma}_A^{*}(h))}\\
  & =& (m/n)
 \Big[(n-h)E\var^{*}(\widehat{\Imath}_{1^{*}}\widehat{\Imath}_{(1+h)^{*}})\\&&+2
 \sum_{t=1}^{n-h-1}(n-h-t)E\cov^{*}(\widehat{\Imath}_{1^{*}}\widehat{\Imath}_{(1+h)^{*}},\widehat{\Imath}_{(1+t)^{*}}\widehat{\Imath}_{(1+t+h)^{*}})
 \Big]\\
 & =&
\Big[ m (1-h/n)(1-\theta)^{2h}
\Big[E(\widehat{\Imath}_1\widehat{\Imath}_{1+h})^2-E
 \Big(n^{-1}\sum_{i=1}^n\widehat{\Imath}_i\widehat{\Imath}_{i+h}
 \Big)^2\Big]\Big]\\ 
 & & +2 m \sum_{t=1}^{n-h-1}
 (1-(h+t)/n)\Big[ n^{-1}\sum_{i=1}^nE\widehat{\Imath}_i\widehat{\Imath}_{i+h}\widehat{\Imath}_{i+t}\widehat{\Imath}_{i+t+h} \Big](1-\theta)^{t+h}\\
 & & +2m
 \sum_{t=1}^{\min(h-1,n-h-1)}(1-(h+t)/n) E\Big(n^{-1}\sum_{i=1}^n
 \widehat{\Imath}_i\widehat{\Imath}_{i+t}\Big)^2
 ((1-\theta)^{2t}-(1-\theta)^{t+h})\\ 
 & & -2m\sum_{t=1}^{\min(h-1,n-h-1)}(1-(h+t)/n)
E \Big(n^{-1}\sum_{i=1}^n\widehat{\Imath}_i\widehat{\Imath}_{i+h}
 \Big)^2(1-\theta)^{2h}\\
 & & -2m \sum_{t=h}^{n-h-1}(1-(h+t)/n)
E \Big(n^{-1}\sum_{i=1}^n\widehat{\Imath}_i\widehat{\Imath}_{i+h}
 \Big)^2(1-\theta)^{t+h}\\
&\le &
m\,E(\widehat{\Imath}_1\widehat{\Imath}_{1+h})^2 
  +2 m \sum_{t=1}^{n-h-1} (1-(h+t)/n)\Big(n^{-1}\sum_{i=1}^nE\widehat{\Imath}_i\widehat{\Imath}_{i+h}\widehat{\Imath}_{i+t}\widehat{\Imath}_{i+t+h}\Big) (1-\theta)^{t+h}\\
 &&+ 2m
 \sum_{t=1}^{\min(h-1,n-h-1)}(1-(h+t)/n) E\Big(n^{-1}\sum_{i=1}^n
 \widehat{\Imath}_i\widehat{\Imath}_{i+t}\Big)^2
 (1-\theta)^{2t}\\
 & =& V_1+V_2+V_3\,.
\end{eqnarray*} }
We observe that, for some constant $c_0>0$,
\beao
V_1&\le & m\,E(\widehat{\Imath}_1\widehat{\Imath}_{1+h})^2\le c\,m\, 
\big[EI_1I_{1+h} + (E\ov I_n)^2\Big]\le c m\,p_0\le c_0\,.
\eeao
For $V_2$, we observe that for $i\le n$,
\beao
&&m \theta^{-1} \big|E\big[\widehat{\Imath}_i\widehat{\Imath}_{i+h}\widehat{\Imath}_{i+t}\widehat{\Imath}_{i+t+h}-
\wt{\Imath}_i\wt{\Imath}_{i+h}\wt{\Imath}_{i+t}\wt{\Imath}_{i+t+h}\big]\big|\\
&\le& c\,  m \theta^{-1}E\big|\ov I_n-p_0\big|= O( \sqrt{m/n} \theta^{-1})=o(1)\,,
\eeao
by virtue of the condition $n\theta^2/m\to\infty$. Therefore, for
showing that $|V_2|\le c$ uniformly for $h,n$,  it suffices to show that
$|\wt V_2|\le c$, where $\wt V_2$ is obtained from $V_2$ by replacing 
the $\wh I_t$'s by the corresponding $\wt I_t$'s. 
Taking into account $E\wt{\Imath}_1\wt{\Imath}_{1+t}= p_t - p_0^2$ and the Cauchy-Schwarz
inequality,
we have for a fixed integer $M>0$, 
\beao
|\wt V_2|& \le &c\,m \sum_{t=1}^{n-h-1}\Big|
n^{-1}\sum_{i=1}^nE\wt{\Imath}_i\wt{\Imath}_{i+h}\wt{\Imath}_{i+t}\wt{\Imath}_{i+t+h}\Big|\\
&= &c\,m \sum_{t=1}^{n-h-1}|
E\wt{\Imath}_1\wt{\Imath}_{1+h}\wt{\Imath}_{1+t}\wt{\Imath}_{1+t+h}| \\
&\le & (m p_0) M+ c\,m \sum_{t=M+1}^{r_n} (p_t+p_0^2)+c\,
m\,\sum_{t=r_n+1}^\infty \xi_t\le c\,,
\eeao
in view of condition (M) and \regvar . A similar argument as for
$V_2$ shows that one may replace the $\wh I_t$'s in $V_3$ by the
corresponding $\wt I_t$'s. We denote the resulting quantity by $\wt
V_3$.
Then we have{\small
\beao
\wt V_3&\le &m \sum_{t=1}^n (1-\theta)^t E\Big(n^{-1}\sum_{i=1}^{n-t}
 \wt{\Imath}_i\wt{\Imath}_{i+t}+ n^{-1}\sum_{i=n-t+1}^{n} \wt I_i \wt
 I_{i+t-n}\Big)^2\\
&\le& c\,m \sum_{t=1}^n (1-\theta)^t E\Big(n^{-1}\sum_{i=1}^{n-t}
 \wt{\Imath}_i\wt{\Imath}_{i+t}\Big)^2+
c\,m \sum_{t=1}^n (1-\theta)^t E\Big(n^{-1}\sum_{i=n-t+1}^{n}
 \wt{\Imath}_i\wt{\Imath}_{i+t-n}\Big)^2\\
&=& \wt V_{31}+\wt V_{32}\,.
\eeao}\noindent
We will only deal with $\wt V_{31}$, the other term can be bounded in
a similar way. We observe that for fixed $M>1$, using condition (M),{\small 
\beao
\wt V_{31}&\le& c\,\frac m n  \sum_{t=1}^n (1-\theta)^t \Big(
 E(\wt I_1\wt I_{1+t})^2 + 2 \sum_{s=1}^{n-t-1}
 |E\wt{\Imath}_1\wt{\Imath}_{1+t}
 \wt{\Imath}_{1+s}\wt{\Imath}_{1+s+t}|\Big)\\
&\le & o(1) + c \frac m n  \sum_{t=1}^n (1-\theta)^t 
  \sum_{s=1}^{n-t-1}
 |E\wt{\Imath}_1\wt{\Imath}_{1+t}
 \wt{\Imath}_{1+s}\wt{\Imath}_{1+s+t}|\\
& \le & o(1)+ c \frac m n  \sum_{t=1}^n (1-\theta)^t\sum_{s=M+1}^{r_n}
(p_s+p_0^2)+  c \frac m n  \sum_{t=1}^n (1-\theta)^t\sum_{r_n+1\le
  s\le n-t-1, s\le t}  \xi_s \\
&&+ c \frac m n  \sum_{t=1}^n (1-\theta)^t\sum_{r_n+1\le
  s\le n-t-1, s> t} \big( |E\wt{\Imath}_1\wt{\Imath}_{1+t}
 \wt{\Imath}_{1+s}\wt{\Imath}_{1+s+t}-(p_t-p_0^2)^2| + (p_t-p_0^2)^2\big)\,.
\eeao}\noindent
In view of condition (M), the first two terms on the \rhs\ are
negligible as $\nto$. The third term is bounded by
\beao
 c \frac m n  \sum_{t=1}^n (1-\theta)^t\sum_{r_n+1\le
  s\le n-t-1, s> t} \xi_{s-t} +
 c m  \sum_{t=1}^n (1-\theta)^t 
(p_t-p_0^2)^2  \,.
\eeao
Multiple use of (M) again shows that the \rhs\ is negligible.
This proves \eqref{eq:opi}. \qed
\end{proof}

\ble \label{lem:moment} Under the conditions of Theorem~\ref{thm:bpclt}
the following relations hold for $s,h\ge 0$ :\footnote{If indices in
  the sums below exceed the value $n$ they are interpreted in the
  circular sense, i.e., ${\rm mod} \; n$.}{\small 
\begin{eqnarray}
\label{eq:mo1}
&&E^{*}\widehat{\Imath}_{1^{*}}  = 0\,,\\
\label{eq:mo2}&&E^{*}\widehat{\Imath}_{1^{*}}\widehat{\Imath}_{(1+h)^{*}}
 = (1-\theta)^h \;n^{-1} \sum_{i=1}^n
 \widehat{\Imath}_i\widehat{\Imath}_{i+h}\,,\quad 
E^{*}\widetilde{\Imath}_{1^{*}}\widetilde{\Imath}_{1^{*}+h}
 =n^{-1} \sum_{i=1}^n
 \wt{\Imath}_i\wt{\Imath}_{i+h}\,,\\
\label{eq:mo30}
 & &\cov^{*}(\widetilde{\Imath}_{1^{*}}\widetilde{\Imath}_{1^{*}+h},\widetilde{\Imath}_{(1+s)^{*}}\widetilde{\Imath}_{(1+s)^{*}+h}
 )\\\nonumber&=& (1-\theta)^s\Big(n^{-1}
 \sum_{i=1}^n\widetilde{\Imath}_i\widetilde{\Imath}_{i+s}\widetilde{\Imath}_{i+h}\widetilde{\Imath}_{i+s+h}-\Big(n^{-1}
\sum_{i=1}^n\widetilde{\Imath}_i\widetilde{\Imath}_{i+h} \Big)^2 \Big)\,, \\
\label{eq:mo31}&&\cov^{*}\big(
\widehat{\Imath}_{1^{*}}\widehat{\Imath}_{(1+h)^{*}},\widetilde{\Imath}_{(1+s)^{*}}\widetilde{\Imath}_{(1+s)^{*}+h}\big)
\\
 \nonumber&=& (1-\theta)^{\max(s,h)} \; \Big(  n^{-1} \sum_{i=1}^n
 \widehat{\Imath}_i \widehat{\Imath}_{i+h}
 \widetilde{\Imath}_{i+s}\widetilde{\Imath}_{i+s+h}- \Big(n^{-1}\sum_{i=1}^n\widehat{\Imath}_i\widehat{\Imath}_{i+h} \Big)
 \Big(n^{-1}\sum_{i=1}^n\widetilde{\Imath}_i\widetilde{\Imath}_{i+h} \Big)\Big)\,,\\
\label{eq:mo32}
 &
 &\cov^\ast\big(\wt I_{1^{*}} \wt I_{1^\ast +h},\widehat{\Imath}_{(1+s)^{*}}\widehat{\Imath}_{(1+s+h)^{*}}
\big)\\
 \nonumber&=& (1-\theta)^{s+h} \; \Big(  n^{-1} \sum_{i=1}^n
 \widetilde{\Imath}_i \widetilde{\Imath}_{i+h}
 \widehat{\Imath}_{i+s}\widehat{\Imath}_{i+s+h}- \Big(n^{-1}\sum_{i=1}^n\widehat{\Imath}_i\widehat{\Imath}_{i+h} \Big)
 \Big(n^{-1}\sum_{i=1}^n\widetilde{\Imath}_i\widetilde{\Imath}_{i+h} \Big)\Big)\,,\\
\label{eq:mo4}&&\cov^{*}
\big(\widehat{\Imath}_{1^{*}}\widehat{\Imath}_{(1+h)^{*}},
\widehat{\Imath}_{(1+s)^{*}}\widehat{\Imath}_{(1+s+h)^{*}}\big)\\ 
\nonumber &&=
\begin{cases}
(1-\theta)^{s+h}\;\Big[
n^{-1}
\sum_{i=1}^n\wh I_i\wh I_{i+s}\wh I_{i+h}\wh I_{i+s+h}
-\Big(n^{-1}\sum_{i=1}^n\widehat{\Imath}_i\widehat{\Imath}_{i+s}\Big)^2\Big]
+\\
\Big(n^{-1}\sum_{i=1}^n\widehat{\Imath}_i\widehat{\Imath}_{i+s}
(1-\theta)^{s}\Big)^2
-\Big(n^{-1}\sum_{i=1}^n\widehat{\Imath}_i\widehat{\Imath}_{i+h}(1-\theta)^h \Big)^2\,, \quad s< h\,,\\
(1-\theta)^{s+h}\;\Big( n^{-1}\sum_{i=1}^n \widehat{\Imath}_i\widehat{\Imath}_{i+s}
  \widehat{\Imath}_{i+h} \widehat{\Imath}_{i+s+h}-
  \Big(n^{-1}\sum_{i=1}^n \widehat{\Imath}_i\widehat{\Imath}_{i+h} \Big)^2\Big)\,, \quad s\ge h\,.
\end{cases}
\end{eqnarray} }
\ele
\begin{proof}
Relations \eqref{eq:mo1} and \eqref{eq:mo2} follow 
from the defining properties of the
stationary bootstrap; see 
Politis and Romano~\cite{politis:romano:1994}.
\par
We will only show that
\eqref{eq:mo4} holds; 
\eqref{eq:mo30}--\eqref{eq:mo32} can be proved in a similar (and
even simpler) way but we omit further details.
First assume $s<h$. Recall $L_1$ from the construction of the stationary
bootstrap scheme. Consider the following decomposition
\beao
\lefteqn{E^{*}[\widehat{\Imath}_{1^{*}}\widehat{\Imath}_{(1+h)^{*}}\widehat{\Imath}_{(1+s)^{*}}\widehat{\Imath}_{(1+s+h)^{*}}]}\\
&=&E^{*}[\widehat{\Imath}_{1^{*}}\widehat{\Imath}_{(1+h)^{*}}\widehat{\Imath}_{(1+s)^{*}}\widehat{\Imath}_{(1+s+h)^{*}}\mid
L_1\le s]\; P(L_1\le s)
\\&&+ E^{*}[\widehat{\Imath}_{1^{*}}\widehat{\Imath}_{(1+h)^{*}}\widehat{\Imath}_{(1+s)^{*}}\widehat{\Imath}_{(1+s+h)^{*}}\mid
s<L_1\le h]\;P(s<L_1\le h)
\\&&+E^{*}[\wh{\Imath}_{1^{*}}\wh{\Imath}_{(1+h)^{*}}\wh{\Imath}_{(1+s)^{*}}\wh{\Imath}_{(1+s+h)^{*}}\mid
h<L_1\le
s+h]\;P(h<L_1\le s+h)
\\&&+
E^{*}[\widehat{\Imath}_{1^{*}}\widehat{\Imath}_{(1+h)^{*}}\widehat{\Imath}_{(1+s)^{*}}\widehat{\Imath}_{(1+s+h)^{*}}\mid
L_1>s+h]\;P(L_1>s+h)\\
&=&Q_1+Q_2+Q_3+Q_4\,.
\eeao
We start with $Q_1$. For $L_1\le s<h$, $\wh I_{1^\ast}$ is independent of\\
$(\wh I_{(1+h)^\ast}, \wh I_{(1+s)^\ast},\wh I_{(1+s+h)^\ast})$, given
$(X_t)$, but $E^\ast \wh I_{1^\ast}=0$ by \eqref{eq:mo1} and therefore $Q_1=0$.
Similarly, for $h<L_1\le s+h$, $\wh I_{(1+s+h)^\ast}$ is independent
of
$(\wh{\Imath}_{1^{*}},\wh{\Imath}_{(1+h)^{*}},\wh{\Imath}_{(1+s)^{*}})$,
given  $(X_t)$, and since $E^{*}\widehat{\Imath}_{(1+s+h)^{*}}=0$,
$Q_3=0$.
Each of the values $i=1,\ldots,n$ has the same chance
to be chosen by the bootstrap, i.e., 
$P^\ast ( \wh I^\ast_1= \wh I_i)=n^{-1}$ for $i=1,\ldots,n$. Thus, for
$L_1>s+h$ and the chosen $i$,
the natural ordering 
$(1^\ast,(1+h)^\ast,(1+s)^\ast,(1+s+h)^\ast)=(i,i+h,i+s,i+s+h)$ 
is preserved and therefore
\beao
Q_4
&=&n^{-1}\sum_{i=1}^n \widehat{\Imath}_i\widehat{\Imath}_{i+s}
  \widehat{\Imath}_{i+h} \widehat{\Imath}_{i+s+h} \;P(L_1>s+h)\\
&=&n^{-1}\sum_{i=1}^n \widehat{\Imath}_i\widehat{\Imath}_{i+s}
  \widehat{\Imath}_{i+h} \widehat{\Imath}_{i+s+h}\;(1-\theta)^{s+h}\,.
\eeao
By a similar argument, \eqref{eq:mo2} and using stationarity, we have
\beao
Q_2&=&
E^{*}[\widehat{\Imath}_{1^{*}}\widehat{\Imath}_{(1+s)^{*}} \mid
s<L_1\le h ]\,
E^{*}[\widehat{\Imath}_{(1+h)^\ast}\widehat{\Imath}_{(1+h+s)^{*}}]\;P(s<L_1\le h) \,\\ & =&
n^{-2}\Big(\sum_{i=1}^n\widehat{\Imath}_i\widehat{\Imath}_{i+s} \Big)^2(1-\theta)^s\;\Big((1-\theta)^s-(1-\theta)^h)\Big)\,.
\eeao
Combining the above expressions and taking into account
\eqref{eq:mo2}, we arrive at \eqref{eq:mo4} for $s<h$.
\par
We proceed with the case $s> h$. Then we have the
corresponding decomposition 
\beao
\lefteqn{E^{*}[\widehat{\Imath}_{1^{*}}\widehat{\Imath}_{(1+h)^{*}}\widehat{\Imath}_{(1+s)^{*}}\widehat{\Imath}_{(1+s+h)^{*}}]}\\
&=&E^{*}[\widehat{\Imath}_{1^{*}}\widehat{\Imath}_{(1+h)^{*}}\widehat{\Imath}_{(1+s)^{*}}\widehat{\Imath}_{(1+s+h)^{*}}\mid
L_1\le h]\; P(L_1\le h)
\\&&+ E^{*}[\widehat{\Imath}_{1^{*}}\widehat{\Imath}_{(1+h)^{*}}\widehat{\Imath}_{(1+s)^{*}}\widehat{\Imath}_{(1+s+h)^{*}}\mid
h<L_1\le s]\;P(h<L_1\le s)
\\&&+E^{*}[\wh{\Imath}_{1^{*}}\wh{\Imath}_{(1+h)^{*}}\wh{\Imath}_{(1+s)^{*}}\wh{\Imath}_{(1+s+h)^{*}}\mid
s<L_1\le
s+h]\;P(s<L_1\le s+h)
\\&&+
E^{*}[\widehat{\Imath}_{1^{*}}\widehat{\Imath}_{(1+h)^{*}}\widehat{\Imath}_{(1+s)^{*}}\widehat{\Imath}_{(1+s+h)^{*}}\mid
L_1>s+h]\;P(L_1>s+h)\\
&=&Q_1'+Q_2'+Q_3'+Q_4'\,.
\eeao
We observe that the \lhs\ is symmetric in $h,s$ and therefore the same
arguments as above show that
$Q_1'=Q_3'=0$, $Q_4=Q_4'$ and 
\beao
Q_2'&=&
E^{*}[\widehat{\Imath}_{1^{*}}\widehat{\Imath}_{(1+h)^{*}}\mid
h<L_1\le s ] \,
E^{*}[\widehat{\Imath}_{(1+s)^{*}}\widehat{\Imath}_{(1+s+h)^{*}}]\;P(h<L_1\le
s)\\
 &= & n^{-2}\Big(\sum_{i=1}^n\widehat{\Imath}_i\widehat{\Imath}_{i+h}
 \Big)^2(1-\theta)^h\; \Big((1-\theta)^h- (1-\theta)^s\Big)
\end{eqnarray*}
The case $h=s$ can be considered as a degenerate case, where $Q_2'=0$.
This  completes the proof of \eqref{eq:mo4}. \qed
\end{proof}

We conclude with a short 
discussion of the bias problem of the bootstrapped integrated \per\
mentioned in Remark~\ref{rem:bias}.
\ble\label{lem:bias}
Assume the conditions of Theorem~\ref{thm:bpclt} and the additional 
condition $\sup_{x\in\Pi}|\psi_h(x)|\le c/h$ for $h\ge 1$ and a
constant $c$.
Then the
following relation holds as $\nto$,{\small
\beam\label{eq:mm}
&&\big(\dfrac nm\big)^{0.5}
\sup_{\lambda\in \Pi}
\Big|\psi_0(\la)\big( E^{*}\wh{\gamma}_A^{*}(0)- \wt \gamma_{A}(0)\big) +2 \sum_{h=1}^{n-1}\psi_h(\lambda)\,
\big(E^{*}\wh{\gamma}_A^{*}(h)- (1-\theta)^h\wt \gamma_{A}(h)\big)\Big|\nonumber\\&&\stp
0\,.
\eeam}\noindent
\ele
\begin{proof} We observe that for $h\ge 0$,
\beam\label{uui}
E^{*}\wh{\gamma}_A^{*}(h)- (1-\theta)^h\wt \gamma_{A}(h)&=&
(1-\theta)^h \big[(\wh \gamma_{A}(h)+\wh\gamma_{A}(n-h))-\wt
\gamma_{A}(h) \big]\nonumber\\
&=&(1-\theta)^h \big[\wt \gamma_{A}(n-h) - m(p_0- \ov I_n)^2\big]\,.
\eeam
For fixed $h$ we have $(n/m)^{0.5} m(p_0- \ov I_n)^2\stp 0$ as $\nto$ and
\beao
(n/m)^{0.5}  E|\wt \gamma_{A}(n-h)|\le c\,(m/n)^{0.5} h p_0 \to
0\,,\quad \nto\,.
\eeao
Therefore it suffices to show that
\beao
\lim_{k\to\infty}\limsup_{\nto} P\Big(\sup_{\lambda\in \Pi}
\Big| \sum_{h=k+1}^{n-1}\psi_h(\lambda)\,
\big(E^{*}\wh{\gamma}_A^{*}(h)- (1-\theta)^h\wt
\gamma_{A}(h)\big)\Big|>\delta\Big)\,,\quad \delta>0\,.
\eeao
Keeping in mind \eqref{uui},
we have 
\beao
(n/m)^{0.5} m(p_0- \ov I_n)^2\sup_{\lambda\in \Pi}
\Big|\sum_{h=k+1}^{n-1}\psi_h(\lambda) (1-\theta)^h\Big|=O_P(1/(\theta
\sqrt{mn}))=o_P(1)\,,
\eeao
where we used $\theta^2 n/m\to \infty$,
and
\beao\lefteqn{
(n/m)^{0.5}
\sup_{\lambda\in \Pi}
\Big|\sum_{h=k+1}^{n-1}\psi_h(\lambda)\,(1-\theta)^h\wt
\gamma_{A}(n-h)\Big|}\\&\le &
(n/m)^{0.5}
\sup_{\lambda\in \Pi}
\Big|\sum_{h=1}^{n-k-1}\psi_{n-h}(\lambda)\,(1-\theta)^{n-h}[\wt
\gamma_{A}(h)- m (1-h/n)
(p_h-p_0^2)]\Big|\\&&+(n/m)^{0.5}\sup_{\lambda\in
  \Pi}\Big|\sum_{h=1}^{n-k-1}\psi_{n-h}(\lambda)\,(1-\theta)^{n-h }m (1-h/n)
(p_h-p_0^2)\Big|\\
&=&I_1+I_2\,.
\eeao
Under the assumption $\sup_{x\in\Pi}|\psi_h(x)|\le c/h$ uniformly for
$h\ge 1$, we have for small $\vep>0$,
\beao
I_2&\le & (m/n)^{0.5} c \sum_{h=1}^{\infty}  \xi_h\to 0\,,\quad \nto\,.
\eeao 
Now we can adapt the proof 
of Lemma~\ref{lem:fclt1} to prove that
\beao
\lim_{k\to\infty}\limsup_{\nto} P(I_1>\delta)=0\,,\quad \delta>0\,.
\eeao
This proves \eqref{eq:mm}.\qed
\end{proof}
However, under the assumptions of Theorem~\ref{thm:main}, it is in general not possible to replace the quantities 
$(1-\theta)^h \wt \gamma_{A}(h)$ in \eqref{eq:mm} by $\wt
\gamma_{A}(h)$, i.e., in general we do not have the relation 
$(n/m)^{0.5} (E^\ast J^\ast_{n,A}- J_{n,A})\stp 0$. Indeed,
taking into account  \eqref{eq:mm} and assuming $\eta$-dependence for
$(X_t)$, we have $E\wt \gamma_A(h)=0$ for $h>\eta$ and
\beao
&&(n/m)^{0.5} (E^\ast J^\ast_{n,A}- J_{n,A})\\&=&2(n/m)^{0.5}
\sum_{h=1}^{n-1}\psi_h(\lambda)\,[(1-\theta)^h-1]\wt
\gamma_{A}(h)+o_P(1)\\&=&2(n/m)^{0.5}\sum_{h=1}^{n-1}\psi_h(\lambda)\,[(1-\theta)^h-1]\big(\wt
\gamma_{A}(h)-E \wt \gamma_A(h)\big)\\
&& + 
2(n/m)^{0.5}\sum_{h=1}^{\eta}\psi_h(\lambda)\,[(1-\theta)^h-1]
(1-h/n) m(p_h-p_0^2)+
o_P(1)\,.
\eeao 
An argument similar to the proof of Theorem~\ref{thm:main} shows that
the first term on the \rhs\ is stochastically bounded, while the
second term may diverge (for example, if $\gamma_A(\eta)>0$ and
$\psi_\eta\ne 0$)
since it is of the order $\theta (n/m)^{0.5}$
which converges to infinity in view of the assumption $\theta^2
n/m\to\infty$ which is vital for the proof of the consistency of 
the stationary bootstrap.

% \section{Acknowledgement}
% This paper was written when Yuwei Zhao was a PhD student at the
% Department of Mathematics of the University of Copenhagen. He would
% like to thank his Department for generous financial support. Major
% parts of this paper were written in 2013 when
% both authors spent sabbatical periods
% at the Department of Statistics at Columbia University and the Forschungsinstitut
% f\"ur Mathematik of ETH Z\"urich. Both authors take pleasure to thank 
% their hosts for financial support and a constructive scientific atmosphere.

\end{document}